# Generalized Gamma Convolutions, Dirichlet means, Thorin measures, with explicit examples*

**Lancelot F. James**[†]

*The Hong Kong University of Science and Technology
Dept. of Information Systems, Business Statistics and Operations Management
Clear Water Bay, Kowloon, Hong Kong SAR e-mail:* `lancelot@ust.hk`

**Bernard Roynette**

*Institut Elie Cartan, Université Henri Poincaré,
B.P. 239, 54506 Vandoeuvre les Nancy Cedex, France*

**Marc Yor**

*Laboratoire de Probabilités et Modèles Aléatoires,
Université Paris VI et VII and
Institut Universitaire de France,
4 place Jussieu - Case 188
F - 75252 Paris Cedex 05,France*

**Abstract:**

- In Section 1, we present a number of classical results concerning the Generalized Gamma Convolution ( : GGC) variables, their Wiener-Gamma representations, and relation with the Dirichlet processes.

- To a GGC variable, one may associate a unique Thorin measure. Let $G$ a positive r.v. and $\Gamma_t(G)$ (resp. $\Gamma_t(1/G)$) the Generalized Gamma Convolution with Thorin measure $t$-times the law of $G$ (resp. the law of $1/G$). In Section 2, we compare the laws of $\Gamma_t(G)$ and $\Gamma_t(1/G)$.

- In Section 3, we present some old and some new examples of GGC variables, among which the lengths of excursions of Bessel processes straddling an independent exponential time.

Received August 2007.

---

*This is an original survey paper
[†]Supported in part by grants HIA05/06.BM03, RGC-HKUST 6159/02P, DAG04/05.BM56 and RGC-HKUST 600907 of the HKSAR.

**Contents**

**Introduction**

**1.** This survey is concerned with the study of a rich and interesting class of infinitely divisible laws on $\mathbb{R}_+$ called **the generalized gamma convolutions** (GGC), a class introduced by O. L. Thorin in 1977 (see [49]) and then studied thoroughly by L. Bondesson [5]; both the lectures notes by Bondesson and the book by Steutel and Van Harn [48] contain many results on this class of laws. We shall also discuss their close connections to a class of random variables known as Dirichlet means whose study was initiated by Cifarelli and Regazzini [10; 11].

We shall often make, throughout this paper, the common abuse of language which consists of talking about a random variable instead of its law; thus, we shall use slightly incorrect terms such as GGC variables, and so on . . .

In order to introduce the family of GGC variables as naturally as possible, let us consider the 3 sets of r.v.'s (or of laws): $\mathcal{G}, \mathcal{S}, \mathcal{J}$ defined as follows:

a) $\mathcal{J}$ is the family of **infinitely divisible** r.v.'s (taking values in $\mathbb{R}_+$)
b) $\mathcal{S}$ is the set of **self-decomposable** r.v.'s, valued in $\mathbb{R}_+$
c) $\mathcal{G}$ is the set of positive **GGC variables**

The Laplace transforms of these variables satisfy:

$$E(e^{-\lambda X}) = \exp\left\{-a\lambda - \int_0^\infty (1 - e^{-\lambda x})\nu(dx)\right\} \qquad (\lambda \geq 0)$$

and

a') If $X \in \mathcal{J}$, then $a \geq 0$ and $\nu(dx)$ is a Lévy measure, i.e.:

$$\int_0^\infty (1 \wedge x)\nu(dx) < \infty$$

b') If $X \in \mathcal{S}$, then $a \geq 0$ and $\nu(dx) = \frac{dx}{x} h(x)$, with $h$ decreasing
c') If $X \in \mathcal{G}$, then $X \in \mathcal{S}$ and moreover:

$$h(x) = \int_0^\infty e^{-xy}\mu(dy)$$

for a $\sigma$-finite and positive measure $\mu$; $\mu(dy)$ is called the Thorin measure of $X$.

Thus, we have:

$$\mathcal{G} \subset \mathcal{S} \subset \mathcal{J}$$

and, of course, the inclusions are strict.

Another definition of a GGC variable is sometimes given in the literature as the limit (in law) of sums of independent gamma variables (with different parameters).

A comparative discussion of different (but, in the end, equivalent) properties of elements of $\mathcal{S}$ is made in Jeanblanc, Pitman, and Yor [32]. Similarly, it is the aim of this survey to gather and to compare different properties for $\mathcal{G}$. This, and more generally, the entire survey is motivated by the fact that, recently, some new elements of $\mathcal{G}$ were discovered; likewise, there is some recurring interest in the Gamma process (see, e.g. the Festschrift Volume for Dilip Madan [55]).

**2.** This survey paper consists of three parts:

• In the first part, the results being discussed are classical; they are about the relationships between different families of r.v.'s and/or processes, namely: GGC r.v.'s, Wiener-Gamma integrals, Dirichlet means, compound Poisson processes, mixing of Gamma variables, Poisson point processes, GGC subordinators, and so on. These results are detailed here in order to ease up the reading of this paper for probabilists coming from diverse horizons.

• In the second part, we discuss a notion of duality for the GGC r.v.'s; in particular, when one knows the density, or the Laplace transform of a GGC r.v. $\Gamma$, this notion of duality allows to compute explicitely the density, or the Laplace transform of the "dual GGC variable". We use this principle to compute the Thorin measure of a Pareto r.v. and of a power of a gamma r.v.

• The third part consists essentially in presenting the explicit computations of densities and of Laplace transforms of some particular GGC r.v.'s. In the main, these r.v.'s originate from the study of the length of excursion which straddles an independent exponential time for a recurrent Bessel process [4]. It is noteworthy that, in this third part, the notion of duality presented in the second part allows to obtain very easily explicit formulae for the density and the Laplace transform of a large number of GGC variables.

Finally, in the Appendix, we describe an interpolation principle between the gamma subordinator $(\gamma_t, \ t \geq 0)$ and a family of GGC subordinators.

## 1. Classical results on GGC r.v.'s

**A writing convention.**
   **1.** Each time we write an equality in law between r.v.'s and that on one or the other side of this equality, several r.v.'s occur, we always assume that these r.v.'s are independent, without mentioning it systematically.
   **2.** It will be convenient, in some instances, to speak of a r.v. instead of its law and vice versa. We hope that no confusion will ensue.

### *1.1. The Gamma process*

It is a subordinator - i.e. a Lévy process with increasing paths and càdlàg trajectories -
   All through this paper, the reference process is the standard Gamma process $(\gamma_t, \ t \geq 0)$, which is a subordinator without drift, and with Lévy measure $\frac{dx}{x} e^{-x}$ $(x > 0)$. Thus, its Lévy-Khintchine representation writes (see [2]):

$$E[e^{-\lambda \gamma_t}] \ = \ \exp\Big\{ -t \int_0^\infty (1 - e^{-\lambda x}) \frac{dx}{x} e^{-x} \Big\} \qquad (\lambda, t \geq 0) \qquad (1)$$

$$= \ \exp\{ -t \log(1 + \lambda) \} = \frac{1}{(1+\lambda)^t} \qquad\qquad\qquad (2)$$

where formula (2) is obtained from (1) and from the elementary Frullani formula (see [34], p. 6):

$$\int_0^\infty (1 - e^{-\lambda x}) \frac{dx}{x} \int_0^\infty e^{-xz} \nu(dz) = \int_0^\infty \log\Big(1 + \frac{\lambda}{z}\Big) \nu(dz) \qquad (3)$$

where, in general, $\nu$ denotes a positive measure on $\mathbb{R}_+$, which is $\sigma$-finite. Here, $\nu(dz) = \delta_1(dz)$ is the Dirac measure at 1, but formula (3) shall be useful in the sequel. We note that Frullani's formula (3) may be easily obtained by observing that the two sides of this formula take the value 0 as $\lambda = 0$ and have the same derivative with respect to $\lambda$.

For each $t > 0$ fixed, $\gamma_t$ follows a gamma law with parameter $t$:

$$P(\gamma_t \in da) = \frac{e^{-a}}{\Gamma(t)} a^{t-1} da \qquad (a \geq 0) \tag{4}$$

This process $(\gamma_t, t \geq 0)$ enjoys a large number of remarkable properties which make it a "worthy companion" of Brownian motion. In particular, Emery and Yor [18] establish a parallel between Brownian motion and its bridges on one hand, and the Gamma process and its bridges on the other hand. See also Vershik, Yor and Tsilevich [50] and Yor [55] for a survey of many remarkable properties of the gamma process.

### 1.2. Wiener-Gamma integrals and GGC variables

**1.2.a** Many times throughout this work, we shall use the properties of the integrals:

$$\widetilde{\Gamma}(h) := \int_0^\infty h(s) \, d\gamma_s \tag{5}$$

where $h : \mathbb{R}_+ \longrightarrow \mathbb{R}_+$ is a Borel function such that:

$$\int_0^\infty \log\big(1 + h(u)\big) du < \infty \tag{6}$$

Under this hypothesis (6), $\widetilde{\Gamma}(h)$ is finite a.s. (see Proposition 1.1 below). Of course, since the trajectories of the process are a.s. increasing, the integral featured in (1.5) may be defined in a path-wise manner, as a usual Stieltjes integral. We call $\widetilde{\Gamma}(h)$ a Wiener-Gamma integral, in analogy with the Wiener integrals $\int_0^\infty f(u) \, dB_u, f \in L^2(\mathbb{R}_+, du)$ which constitute the Gaussian space generated by Brownian motion $(B_u, u \geq 0)$.

Thus, the family of r.v.'s $\{\widetilde{\Gamma}(h)$ with $h$ such that $\int_0^\infty \log\big(1 + h(u)\big) du < \infty\}$ constitutes the analogue for the process $(\gamma_t, t \geq 0)$ of the first Wiener chaos for a Brownian motion $(B_u, u \geq 0)$.

Let us assume for a moment that, in (1.5), the function $h$ is constant on a finite number of intervals, i.e.:

$$h(s) = \sum_{i=0}^{n-1} h_i 1_{]s_i, s_{i+1}]}(s) \qquad (h_i \geq 0)$$

for a subdivision $\quad 0 = s_0 < s_1 < \cdots < s_n < \infty$. Then:

$$\widetilde{\Gamma}(h) = \sum_{i=0}^{n-1} h_i(\gamma_{s_{i+1}} - \gamma_{s_i})$$

Thus, $\widetilde{\Gamma}(h)$ is a linear combination of independent gamma r.v.'s and we obtain:

$$
\begin{aligned}
E\big(e^{-\lambda\widetilde{\Gamma}(h)}\big) &= \prod_{i=0}^{n-1} E\big(e^{-\lambda h_i(\gamma_{s_{i+1}}-\gamma_{s_i})}\big) \qquad (\lambda \geq 0) \\
&= \prod_{i=0}^{n-1} \frac{1}{(1+\lambda h_i)^{s_{i+1}-s_i}} \qquad \big(\text{from } (1.2)\big) \\
&= \exp\left(-\sum_{i=0}^{n-1}(s_{i+1}-s_i)\int_0^\infty (1-e^{-\lambda h_i x})\frac{dx}{x}e^{-x}\right) \\
&\qquad \big(\text{from Frullani's formula } (1.3)\big) \\
&= \exp\left(-\int_0^\infty (1-e^{-\lambda y})\frac{dy}{y}\int_0^\infty e^{-\frac{y}{h(s)}}ds\right) \\
&\qquad (\text{after making the change of variable } h_i x = y) \\
&= \exp\left(-\int_0^\infty (1-e^{-\lambda y})\frac{dy}{y}\int_0^\infty e^{-yx}\mu_h(dx)\right)
\end{aligned}
$$

where $\mu_h$ is the image of Lebesgue's measure on $\mathbb{R}_+$ by the application $s \to \frac{1}{h(s)}$. This latter formula justifies the following definition:

**1.2.b Definition 1.0**: Following ([5], p. 29), we say that a positive r.v. $\Gamma$ is a generalized gamma convolution (GGC) - without translation term - if there exists a positive Radon measure $\mu$ on $]0,\infty[$ such that:

$$
\begin{aligned}
E[e^{-\lambda\Gamma}] &= \exp\Big\{-\int_0^\infty (1-e^{-\lambda x})\frac{dx}{x}\int_0^\infty e^{-xz}\mu(dz)\Big\} & (7) \\
&= \exp\Big\{-\int_0^\infty \log\Big(1+\frac{\lambda}{z}\Big)\mu(dz)\Big\} & (8)
\end{aligned}
$$

with: $\quad \int_{]0,1]} |\log x|\,\mu(dx) < \infty \text{ and } \int_{[1,\infty[} \frac{\mu(dx)}{x} < \infty \qquad (9)$

The measure $\mu$ is called Thorin's measure associated with $\Gamma$. Thus, from the Lévy-Khintchine formula, a GGC r.v. is infinitely divisible. In fact, since its Lévy density $l_\Gamma(x) = \frac{1}{x}\int_0^\infty e^{-xz}\mu(dz)$ satisfies: $x \longrightarrow x\,l(x)$ is decreasing, then $\Gamma$ is a self decomposable r.v. (see, e.g. [37]). Such a self-decomposable r.v. $\Gamma$, assumed to be non degenerate, admits a density $f_\Gamma$ such that $f_\Gamma(x) > 0$ for every $x > 0$ (see [45] p.404 ). The study of GGC variables was initiated by O. Thorin in a series of papers (see for instance [49]).

**1.2.c** GGC variables and Wiener-Gamma representations.
The following Proposition is classical. The reader may refer to Lijoi and Regazzini [36].

**Proposition 1.1** *The class of positive GGC variables coincides with the class of Wiener-Gamma integrals. More precisely:*

**1.** If $\widetilde{\Gamma}(h) = \int_0^\infty h(u)\,d\gamma_u$, then:

$$E(e^{-\lambda \widetilde{\Gamma}(h)}) = \exp\Big\{-\int_0^\infty \log\Big(1+\frac{\lambda}{x}\Big)\mu_h(dx)\Big\} \tag{10}$$

where $\mu_h$ denotes the image of Lebesgue's measure on $\mathbb{R}_+$ under the application: $s \longrightarrow \frac{1}{h(s)}$. In other terms:

$$\int_0^\infty e^{-\frac{x}{h(s)}}\,ds = \int_0^\infty e^{-xz}\mu_h(dz) \qquad (x>0) \tag{11}$$

We note that, in (11), $h$ may vanish on some measurable set.

**2.** Let $\Gamma$ denote a GGC r.v. with Thorin measure $\mu$. Let $F_\mu(x) := \int_{]0,x]} \mu(dy)$ ($x \geq 0$) and denote by $F_\mu^{-1}$ its right continuous inverse, in the sense of the composition of functions.

Then:

$$\Gamma \stackrel{(law)}{=} \widetilde{\Gamma}(h), \ \text{with} \ h(u) = \frac{1}{F_\mu^{-1}(u)} \tag{12}$$

**Proof of Proposition 1.1** Let $\widetilde{\Gamma}(h) := \int_0^\infty h(u)\,d\gamma_u$. It is easily obtained, by approximation of $\widetilde{\Gamma}_h$ by Riemann sums, using also the fact that the Lévy measure of the process $(\gamma_t,\ t\geq 0)$ equals $\frac{dx}{x}e^{-x}$, that:

$$\begin{aligned}
E(e^{-\lambda \widetilde{\Gamma}(h)}) &= \exp\Big\{-\int_0^\infty du \int_0^\infty (1-e^{-\lambda x h(u)})\frac{dx}{x}e^{-x}\Big\} \\
&= \exp\Big\{-\int_0^\infty (1-e^{-\lambda y})\frac{dy}{y}\int_0^\infty e^{-\frac{y}{h(u)}}\,du\Big\} \tag{13}
\end{aligned}$$

after making the change of variable $xh(u) = y$. We observe, from (11) and (12), the equivalence of the conditions:

$$\int_0^\infty \log(1+h(u))\,du < \infty \Leftrightarrow \int \log\Big(1+\frac{1}{x}\Big)\mu_h(dx) < \infty$$

$$\Leftrightarrow \int_{]0,1]} |\log x|\,\mu_h(dx) < \infty \ \text{and} \ \int_{[1,\infty[} \frac{\mu_h(dx)}{x} < \infty$$

**Remark 1.2**

**1.** Formula (13) may be obtained in a slightly different manner: the process $(\gamma_s - \gamma_{s^-} := e_s,\ s \geq 0)$ of jumps of the subordinator $(\gamma_t,\ t\geq 0)$ is a Poisson point process whose intensity measure $n$ equals the Lévy measure of $(\gamma_t,\ t\geq 0)$ (see [2]):

$$n(dx) = \frac{1}{x}e^{-x}dx \tag{14}$$

Thus, from the exponential formula for Poisson point processes ([43], p. 476):

$$E\Big[\exp-\lambda \sum_{0<s\leq t} f(s,e_s)\Big] = \exp\Big\{-\int_0^t ds \int_0^\infty (1-e^{-\lambda f(s,u)})n(du)\Big\}$$

one gets:

$$\begin{aligned}
E[\exp-\lambda \int_0^\infty h(s)d\gamma_s] &= E(\exp-\lambda \sum_{s>0} h(s)(\Delta\gamma_s)) \\
&= \exp\left\{-\int_0^\infty ds \int_0^\infty (1-e^{-\lambda h(s)u})\frac{e^{-u}}{u}du\right\},
\end{aligned}$$

which agrees with the expression in (13).

**2.** Let $h, k : \mathbb{R}_+ \longrightarrow \mathbb{R}_+$ two Borel functions which satisfy (6) and assume that $\widetilde{\Gamma}(h) \stackrel{(law)}{=} \widetilde{\Gamma}(k)$. Relation (13) and the uniqueness of the Lévy measure in the Lévy-Khintchine representation imply: the images by $h$ and $k$ of Lebesgue's measure on $\mathbb{R}_+$ are identical. Thus, choosing for $k$ the increasing rearrangement $h^*$, (resp.: the decreasing rearrangement $h_*$) of $h$, we obtain that there exists essentially a unique increasing function $h^*$, (resp. a unique decreasing function $h_*$) such that:

$$\widetilde{\Gamma}(h) \stackrel{(law)}{=} \widetilde{\Gamma}(h^*) = \widetilde{\Gamma}(h_*) \tag{15}$$

We recall that the function $h^*$ (resp. the function $h_*$) is the unique (equivalence class of) increasing (resp. decreasing) function such that for every $a \geq 0$:

$$\text{meas}\,(x\,;\,h(x) < a) = \text{meas}\,(x,\,h^*(x) < a) = \text{meas}\,(x\,;\,h_*(x) < a)$$

where meas indicates Lebesgue's measure on $\mathbb{R}_+$.

From Proposition 1.1, $\widetilde{\Gamma}(h)$ is a GGC r.v. and we shall denote by $\mu_h$ the Thorin measure associated with $\widetilde{\Gamma}(h)$.

### 1.3. m-Wiener-Gamma integrals, $(m, G)$ GGC r.v.

In this work, we shall often consider a GGC r.v. whose associated Thorin measure has finite total mass. Thus, we shall now particularize Proposition 1.1 in this case.

**1.3.a** Let $m > 0$ and $h : [0, m] \to \mathbb{R}_+$ a Borel function such that:

$$\int_0^m \log(1+h(u))du < \infty \tag{16}$$

We call $m$-Wiener integral of $h$ the r.v.:

$$\widetilde{\Gamma}_m(h) := \int_0^m h(u)d\gamma_u \tag{17}$$

Since $(\gamma_u,\,0 \leq u \leq m)$ and $(\gamma_m - \gamma_{(m-u)-},\,0 \leq u \leq m)$ have the same law, we deduce from (17), after making the change of variable $m - u = v$:

$$\widetilde{\Gamma}_m(h) \stackrel{(law)}{=} \int_0^m h(u)d\gamma_u \stackrel{(law)}{=} \int_0^m h(m-u)d\gamma_u \tag{18}$$

We note, in relation with (15) above, that if $h$ is increasing, resp. decreasing, then $u \to h(m-u)$ is decreasing, resp. increasing.

**1.3.b** Let $m > 0$ and $G$ be a positive r.v. such that:

$$E\big(\log^+(1/G)\big) < \infty \tag{19}$$

We say that a positive r.v. $\Gamma$ is a $(m,G)$ GGC if:

$$E(e^{-\lambda \Gamma}) = \exp\Big\{ -m \int_0^\infty (1 - e^{-\lambda x}) \frac{dx}{x} E(e^{-xG}) \Big\} \tag{20}$$

Of course, from (7), a $(m,G)$ GGC r.v. is a GGC r.v. whose Thorin measure $\mu$ equals:

$$\mu(dx) = m\, P_G(dx) \tag{21}$$

where $P_G$ denotes the law of $G$ and we have:

$$\mu(\,]0, \infty[\,) = m \tag{22}$$

Under (21), it is clear that:

$$\int_{]0,1]} |\log x|\, \mu(dx) < \infty \text{ and } \int_{[1,\infty[} \frac{\mu(dx)}{x} < \infty$$
$$\Leftrightarrow E\big(\log^+(1/G)\big) < \infty \tag{23}$$

We denote by $\Gamma_m(G)$ the (law of the) r.v. $\Gamma$ defined by (20). Hence:

$$E[e^{-\lambda \Gamma_m(G)}] = \exp\Big\{ -m \int_0^\infty (1 - e^{-\lambda x}) \frac{dx}{x} E(e^{-xG}) \Big\} \tag{24}$$

$$= \exp\Big\{ -m\, E\Big(\log\Big(1 + \frac{\lambda}{G}\Big)\Big) \Big\} \tag{25}$$

**1.3.c** $G$ or $1/G$ ? How to choose ?

We have used the notation $\Gamma_m(G)$ due to the relation (24). Of course, the relation (25) invites, on the contrary, to adopt the notation $\Gamma_m(1/G)$. However, we shall not adopt this latter notation as the notation $\Gamma_m(G)$ is used by L. Bondesson [5] who has contributed in an essential manner to the study of the GGC variables.

**1.3.d** Proposition 1.1, when the Thorin measure has a finite total mass $m$, admits the following translation.

**Proposition 1.3** *A r.v. is a $(m,G)$ GGC if and only if it is a $m$-Wiener integral. More precisely:*

**1)** If $\qquad \Gamma_m(G) \stackrel{(law)}{=} \widetilde{\Gamma}_m(h)$, then

$$h(u) = \frac{1}{F_G^{-1}(u/m)} \qquad (u \in [0,m]) \tag{26}$$

where $F_G^{-1}$ denotes the right continuous inverse, in the sense of composition of functions, of $F_G$, the cumulative distribution function of $G$.

**2)** If $\Gamma_m(h) \stackrel{(law)}{=} \Gamma_m(G)$, then

$$G \stackrel{(law)}{=} \frac{1}{h(U_m)} \tag{27}$$

where $U_m$ denotes a uniform r.v. on $[0, m]$.

**1.3.e** Some classical results

We gather here some results which are due to L. Bondesson [5] and which we shall use in the sequel. Let $m > 0$ and $G$ satisfy (19). Then, denoting $f_{\Gamma_m(G)}$ the density of $\Gamma_m(G)$:

$$\bullet \qquad f_{\Gamma_m(G)}(x) = \frac{x^{m-1}}{\Gamma(m)} g(x) \qquad (x > 0) \tag{28}$$

where $g$ is a completely monotone function ([5], p. 49). Moreover ([5], p. 50) $g$ admits a limit on the right of 0 and:

$$g(0_+) = \exp\{m\, E\,(\log G)\} \tag{29}$$

We note that, since by hypothesis $E(\log^+(1/G)) < \infty$, $g(0_+)$ is finite if and only if $E(\log^+(G)) < \infty$. In Section 2 of this work (see Theorem 2.1) we give an explicit form of $g$ when $E(|\log G|) < \infty$.

• $m$ may be determined from the knowledge of $f_{\Gamma_m(G)}$:

$$m = \sup\left\{\alpha \geq 0\; ;\; \lim_{x \downarrow 0_+} \frac{f_{\Gamma_m(G)}(x)}{x^{\alpha-1}} = 0\right\} \tag{30}$$

(see [5], p. 51).

### 1.4. m-Wiener Gamma integrals, m Dirichlet means, Gamma (m) mixtures

**1.4.a** The preceding discussion leads us to introduce, for every $m > 0$, the Dirichlet process with parameter $m$, denoted: $(D_u^{(m)},\ 0 \leq u \leq m)$ and defined as:

$$(D_u^{(m)},\ 0 \leq u \leq m) = \left(\frac{\gamma_u}{\gamma_m},\ 0 \leq u \leq m\right) \tag{31}$$

It is well known, and it is an easy consequence of the properties of the "beta-gamma algebra" that this process $(D_u^{(m)},\ 0 \leq u \leq m)$ is independent from the r.v. $\gamma_m$, hence from $(\gamma_v,\ v \geq m)$. Indeed, if $\gamma_a$ and $\gamma_b$ are two independent gamma variables with parameter $a$, resp. $b$, the basic "beta-gamma algebra" states that:

$$\left(\frac{\gamma_a}{\gamma_a + \gamma_b},\ \gamma_a + \gamma_b\right) \stackrel{(\text{law})}{=} (\beta_{a,b}, \gamma_{a+b})$$

where $\beta_{a,b}$ and $\gamma_{a+b}$ are independent and are respectively beta $(a,b)$ and gamma $(a+b)$ distributed. In particular, $\frac{\gamma_a}{\gamma_a+\gamma_b}$ is independent from $\gamma_a + \gamma_b$. Thus, for every $u \leq m$:

$$\frac{\gamma_u}{\gamma_m} = \frac{\gamma_u}{\gamma_u + (\gamma_m - \gamma_u)}$$

is independent from $\gamma_u + (\gamma_m - \gamma_u) = \gamma_m$

This allows to write, for $h$ which satisfies (16):

$$\int_0^m h(u)\, d\gamma_u \stackrel{(law)}{=} \gamma_m \cdot \int_0^m h(u)\, d_u(D_u^{(m)}) \tag{32}$$

Thus, from Proposition 1.3, we may write for $G$ which satisfies (19):

$$\Gamma_m(G) \stackrel{(law)}{=} \gamma_m \cdot D_m(G) \tag{33}$$

with

$$D_m(G) := \int_0^m \frac{1}{F_G^{-1}(u/m)}\, d_u(D_u^{(m)}) \tag{34}$$

It follows that for every $(m, G)$ GGC r.v., the r.v. $\Gamma_m(G)$ is a Gamma $(m)$ mixture, i.e. it may be written as:

$$\Gamma_m(G) \stackrel{(law)}{=} \gamma_m \cdot Z \tag{35}$$

where $Z$ is a positive r.v.

In general, a relationship of the kind:

$$X \cdot Z \stackrel{(law)}{=} X \cdot Z'$$

(with $X$ and $Z$ independent, and $X$ and $Z'$ independent) does not allow, "via simplification" to conclude that $Z \stackrel{(law)}{=} Z'$. However, when $X$ is a gamma variable, this "simplification" is licit. More precisely:

**1.4.b** The relation (35) determines the law of $Z$. Indeed, let $Z$ and $Z'$ two positive r.v.'s such that:

$$\gamma_m \cdot Z \stackrel{(law)}{=} \gamma_m Z'$$

Then, for every $s \in \mathbb{R}$:

$$E[\gamma_m^{is}]\, E[Z^{is}] = E[\gamma_m^{is}]\, E[Z'^{is}]$$

Hence:

$$E(Z^{is}) = E(Z'^{is}) \quad \text{and} \quad Z \stackrel{(law)}{=} Z'$$

**1.4.c Remark 1.4**
(We shall not use the present Remark in the sequel of this paper). We come back to the notation of point 1 of Remark 1.2 and we denote: $(J_1^{(m)} \geq J_2^{(m)} \geq$ $\cdots$) the sequel of the lengths of jumps of the process $(\gamma_u, \ u \leq m)$ ranked in decreasing order. It is not difficult to see that since the intensity measure of the Poisson point process $((s, e_s), \ s \geq 0)$ is $ds \frac{e^{-x}}{x} dx$, then the jump times $(U_1^{(m)}, U_2^{(m)}, \cdots)$ constitute a sequence of i.i.d r.v.'s with uniform law on $[0, m]$ which is independent from the sequence $(J_k^{(m)}, \ k \geq 1)$. Thus:

$$\int_0^m h(u) d\gamma_u \stackrel{(law)}{=} \int_0^m \frac{1}{F_G^{-1}(u/m)} \, d\gamma_u \quad \text{(from (26))}$$
$$= \sum_{k \geq 1} \frac{1}{F_G^{-1}\left(\frac{U_k^{(m)}}{m}\right)} J_k^{(m)} = \sum_{k \geq 1} J_k^{(m)} \frac{1}{G_k} \quad (36)$$

where $\left(\frac{1}{G_k}, \ k \geq 1\right)$ is the sequence of i.i.d r.v.'s with common law $1/G$ and is independent (as a sequence) from the r.v.'s $(J_k^{(m)}, \ k \geq 1)$. Indeed:

$$P\left[\frac{1}{F_G^{-1}\left(\frac{U_k^{(m)}}{m}\right)} \leq x\right] = P\left[F_G^{-1}\left(\frac{U_k^{(m)}}{m}\right) \geq \frac{1}{x}\right]$$
$$= P\left[\frac{U_k^{(m)}}{m} \geq F_G\left(\frac{1}{x}\right)\right] = 1 - F_G\left(\frac{1}{x}\right)$$
$$= P\left(\frac{1}{G} \leq x\right)$$

since $\frac{U_k^{(m)}}{m}$ is uniform on $[0, 1]$. We deduce from (36) that:

$$\int_0^m h(u) d_u(D_u^{(m)}) = \frac{1}{\gamma_m} \cdot \int_0^m h(u) \, d\gamma_u = \sum_{k \geq 1} \frac{J_k^{(m)}}{\gamma_m} \cdot \frac{1}{G_k} \quad (37)$$

We note that: $\sum_{k \geq 1} \frac{J_k^{(m)}}{\gamma_m} = 1$.

We then define the random Dirichlet measure $P_{0,m}^{(1/G)}(dx)$ by the formula:

$$P_{0,m}^{(1/G)}(dx) := \sum_{k \geq 1} \frac{J_k^{(m)}}{\gamma_m} \, \delta_{1/G_k}(dx)$$

and we obtain, from (34) and (37):

$$D_m(G) \stackrel{(law)}{=} \int_0^\infty x \, P_{0,m}^{(1/G)}(dx) \quad (38)$$

This relation (38) justifies the denomination, for $D_m(G)$, of a Dirichlet means.

The study of the Dirichlet means can be traced back to an early work of Cifarelli and Regazzini [10] which culminates into the more recognized [11]. Additional early works on this topic include [14; 21; 26; 56].

See also Bertoin [3] for an example of (34) where he shows that a Cauchy random variable, $C_1$, may be represented as

$$C_1 \stackrel{(law)}{=} -\frac{1}{\gamma_1} \int_0^1 \cot(\pi s) d\gamma_s \stackrel{(law)}{=} \int_0^1 M'_s ds.$$

Where, $(M'_s, s \in [0, 1[)$ is the right derivative of the convex minorant of a Cauchy process.

**1.4.d** Multiplication by a beta variable.

In this section we discuss what happens when $D_m(G)$ is multiplied by certain independent beta random variables. This idea, and restatements of the results (1.4.d i-iii) below, first appear in James (see [28], Theorem 3.1, revised in [29]).

**1.4.d i)** We denote by $\mathcal{D}^{(m)}$ the set of the laws of r.v.'s of the form

$$\int_0^m h(u) d_u(D_u^{(m)}) := D_m(h), \quad \text{with} \int_0^m \log(1 + h(u)) du < \infty, \qquad h \geq 0.$$

If $D_m(h) \in \mathcal{D}^{(m)}$ and if $\beta_{m,m'-m}$, with $m' > m$, is a beta r.v. with parameters $(m, m' - m)$ then:

$$\beta_{m,m'-m} \cdot D_m(h) \in \mathcal{D}^{(m')} \tag{39}$$

Indeed, $\frac{\gamma_m}{\gamma_{m'}}$ is independent from $\gamma_{m'}$ and follows a beta law, with parameters $(m, m' - m)$. Hence:

$$\beta_{m,m'-m} \cdot D_m(h) \stackrel{(law)}{=} \frac{\gamma_m}{\gamma_{m'}} \cdot \int_0^m h(u) \, d_u(D_u^{(m)})$$

$$\stackrel{(law)}{=} \frac{\gamma_m}{\gamma_{m'}} \int_0^m h(u) \frac{d\gamma_u}{\gamma_m}$$

$$\stackrel{(law)}{=} \frac{1}{\gamma_{m'}} \int_0^{m'} h(u) 1_{[0,m]}(u) d_u(D_u^{(m')}) \in \mathcal{D}^{(m')}$$

**1.4.d ii)** In the same spirit as for the preceding point, we note that, if $G$ is a positive r.v. such that $E(\log^+ G) < \infty$ and if $m' > m$, then:

$$\Gamma_m\left(\frac{1}{G}\right) \stackrel{(law)}{=} \Gamma_{m'}\left(\frac{1}{G \cdot Y_p}\right) \tag{40}$$

where, on the RHS of (40), $G$ and $Y_p$ are independent and $Y_p$ is a Bernoulli r.v. with parameter $p = \frac{m}{m'}$:

$$P(Y_p = 1) = p = 1 - P(Y_p = 0) \tag{41}$$

Indeed, we deduce from (24):

$$\begin{aligned}
E(e^{-\lambda \Gamma_m(1/G)}) &= \exp\left\{-m \int_0^\infty (1-e^{-\lambda x})\frac{dx}{x} E(e^{-x/G})\right\} \\
&= \exp\left\{-m' \int_0^\infty (1-e^{-\lambda x})\frac{dx}{x}\frac{m}{m'} E(e^{-x/G})\right\} \\
&= \exp\left\{-m' \int_0^\infty (1-e^{-\lambda x})\frac{dx}{x} E(e^{-(\frac{x}{GY_p})})\right\} \\
&= E(e^{-\lambda \Gamma_{m'}(\frac{1}{GY_p})})
\end{aligned}$$

**1.4.d iii)** We now write the relation (40) in a slightly different manner with the introduction of the r.v.'s $D_m(1/G)$ and $D_m\bigl(\frac{1}{GY_p}\bigr)$. We have, from (40) and (35) for $m' > m$:

$$\Gamma_m(1/G) \stackrel{(law)}{=} \gamma_m D_m(1/G) \stackrel{(law)}{=} \Gamma_{m'}\left(\frac{1}{GY_p}\right) \stackrel{(law)}{=} \gamma_{m'} \cdot D_{m'}\left(\frac{1}{GY_p}\right)$$

Hence:

$$\frac{\gamma_m}{\gamma_{m'}} \gamma_{m'} D_m(1/G) \stackrel{(law)}{=} \gamma_{m'} D_{m'}\left(\frac{1}{GY_p}\right)$$

so that, from point **1.4.b**:

$$\frac{\gamma_m}{\gamma_{m'}} \cdot D_m(1/G) \stackrel{(law)}{=} D_{m'}\bigl(\frac{1}{GY_p}\bigr) \quad \text{i.e. :}$$

$$\beta_{m,m'-m} \cdot D_m(1/G) \stackrel{(law)}{=} D_{m'}\left(\frac{1}{GY_p}\right) \quad \left(p = \frac{m}{m'}\right) \tag{42}$$

In particular, for $m < 1$ and $m' = 1$:

$$\beta_{m,1-m} \cdot D_m(1/G) \stackrel{(law)}{=} D_1\left(\frac{1}{GY_m}\right) \tag{43}$$

**1.4.d iv)** Some elements of $\mathcal{D}^{(m)}$.

Let $T$ denote a positive r.v. which belongs to the Bondesson class $\mathbb{B}$ (see [5], p. 73, Th. 5.2.2) i.e. whose density $f_T$ writes:

$$f_T(x) = C\, x^{\beta-1} h_1(x) h_2(1/x) \qquad (\beta \in \mathbb{R},\ x > 0)$$

with

$$h_j(x) = \exp\left\{-b_j\, x + \int_0^\infty \log\left(\frac{1+y}{x+y}\right) \nu_j(dy)\right\} \qquad j = 1, 2$$

and

$$\int_0^\infty \frac{\nu_j(dy)}{1+y} < \infty$$

Then, from Bondesson ([5], Th. 5.2.2, p. 79), we know that:

$$\Gamma_m := \gamma_m \cdot T \tag{44}$$

is GGC (for every $m > 0$) with some, possibly unknown, associated Thorin measure we denote as $\mu_m$. Assuming furthermore that $E(T^{-m}) < \infty$, we get, from (30):

$$\mu_m(]0, \infty[) = \sup\left\{\alpha \geq 0 \; ; \; \lim_{x \searrow 0_+} \frac{f_{\Gamma_m}(x)}{x^{\alpha-1}} = 0\right\} \tag{45}$$

But, an elementary computation, starting from (44), shows that:

$$f_{\Gamma_m}(x) = \frac{1}{\Gamma(m)} \, x^{m-1} E\left(e^{-\frac{x}{T}} \frac{1}{T^m}\right) \tag{46}$$

Thus, we have:

$$\mu_m(]0, \infty[) = m$$

Hence, there exists, from Proposition 1.3, a positive r.v. $G_m$ such that $E\left(\log^+\left(\frac{1}{G_m}\right)\right) < \infty$ and also such that, from (33):

$$\Gamma_m \stackrel{(law)}{=} \gamma_m \cdot T = \Gamma_m(G_m) = \gamma_m \cdot D_m(G_m) \tag{47}$$

Thus, from point **1.4.b**

$$T \stackrel{(law)}{=} D_m(G_m) \in \mathcal{D}^{(m)} \tag{48}$$

We now summarize what we have just obtained:

**Proposition 1.5**
Let $T$ denote a positive r.v. which belongs to $\mathbb{B}$, such that:

$$E\left(T^{-m}\right) < \infty$$

Then:

$$1) \quad T \in \mathcal{D}^{(m)}$$

$$2) \quad \textit{For every } m' > m, \; \beta_{m,m'-m} \cdot T \in \mathcal{D}^{(m')} \quad \big(\text{from } (43)\big)$$

In particular, Proposition 1.5 may be applied in the following cases:

- If $T$ is a generalized inverse Gaussian r.v., i.e. its density is given by:

$$f_T(x) = C\, x^{\beta-1} \exp\left\{-\frac{1}{2}\left(c_1 x + \frac{c_2}{x}\right)\right\} \cdot 1_{[0,\infty[}(x)$$

  ($\beta \in \mathbb{R}$, $c_1, c_2 > 0$) then, for every $m > 0$, $T \in \mathcal{D}^{(m)}$ (although $T$ is a GGC variable with Thorin measure of infinite total mass (see [5], p. 59).
- If $T$ is a Gamma r.v. $\gamma_\theta$ with parameter $\theta$ then for every $m > 0$, $\gamma_\theta \in \mathcal{D}^{(m)}$. Indeed $\gamma_\theta \in \mathbb{B}$ and for every $m < \theta$: $E\left(\frac{1}{\gamma_\theta^m}\right) < \infty$. On the other hand, assuming still $m < \theta$, we have:

$$\beta_{m,\theta-m} \cdot \gamma_\theta \stackrel{(law)}{=} \gamma_m \in \mathcal{D}^{(\theta)} \quad \text{from (40)}$$

- If $T$ is a positive stable r.v. with index $\alpha$:

$$E(e^{-\lambda T}) = \exp\{-C\lambda^\alpha\} \qquad (\lambda, C > 0, \alpha \in\, ]0,1[)$$

then, for $\alpha \leq \frac{1}{2}$, $T \in \mathbb{B}$ (see [5], p. 85–88) and hence $T \in \mathcal{D}^{(m)}$ for every $m > 0$ since $E\left(\frac{1}{T^m}\right) < \infty$.

**Remark 1.5** Epifani, Guglielmi and Melilli ([16], section 4; see also [17]), posed the natural question of which kind of probability measures are the laws of Dirichlet means. They were able to find some examples in cases where those particular random variables possessed all finite moments. One sees that Proposition 1.5, in a rather simple way, identifies a large number of possible distributions.

**1.4.e** Another representation of $\Gamma_m(G)$.

We have been interested mainly until now in the distributions of the $(m, G)$ GGC r.v.'s. We shall now describe a realization of such a r.v. with the help of a compound Poisson process. Besides, this realization allows us to show that a $(m, G)$ GGC solves an "affine equation". For a nice survey of these equations, see Vervaat [51].

Let $m > 0$ and let $K$ be a positive r.v. We shall say that $(Y_t,\, t \geq 0)$ is a $(m, K)$ $\mathbb{R}_+$ valued compound Poisson process if:

$$Y_t := \sum_{i=1}^{N_t} K_i$$

where $(K_1, K_2, \cdots)$ is a sequence of i.i.d. variables, distributed as $K$, and with $(N_t,\, t \geq 0)$ a Poisson process with parameter $m$, independent of the sequence $(K_i, i = 1, 2 \cdots)$. In particular, $N_t$ is a Poisson r.v. with parameter $mt$.

**Proposition 1.6.**

*Let $\Gamma_m(G)$ a $(m, G)$ GGC r.v. $\left(\text{with } E\left(\log^+\left(\frac{1}{G}\right)\right) < \infty \text{ and } m > 0\right)$. Define $K$ by:*

$$K \stackrel{(law)}{=} \frac{\mathfrak{e}}{G} \qquad (\mathfrak{e},\ a\ standard\ exponential\ r.v.\ independent\ of\ G)$$

*Then*

$$\textbf{1)} \qquad \Gamma_m(G) \stackrel{(law)}{=} \int_0^\infty e^{-t} dY_t$$

*where $(Y_t,\, t \geq 0)$ is a $(m, K)$ compound Poisson process.*

**2)** *$\Gamma_m(G)$ satisfies the affine equation:*

$$\Gamma_m(G) \stackrel{(law)}{=} U^{1/m}\bigl[\Gamma_m(G) + K\bigr]$$

*$\bigl(\text{with } U \text{ uniformly distributed on } [0,1]\bigr)$.*

**Proof of Proposition 1.6.**
*i)* <u>We first prove point 1</u>. We consider $(Y_t, \ t \geq 0)$ a $(m, K)$ compound Poisson process. Then, approximating $\int_0^\infty e^{-t} dY_t$ by the Riemann sums

$$\sum (Y_{t_{i+1}} - Y_{t_i}) e^{-t_i},$$

we obtain:

$$E\left(\exp\left\{-\lambda \int_0^\infty e^{-t} dY_t\right\}\right) = \exp\left\{-\int_0^\infty (1 - e^{-\lambda v}) \mu([v, \infty]) \frac{dv}{v}\right\}$$

where $\mu$ is the Lévy measure of the subordinator $(Y_t, \ t \geq 0)$. Since this subordinator is a $(m, K)$ compound Poisson process, we have:

$$\mu([v, \infty[) = m P(K \geq v) = m P\left(\frac{\mathfrak{e}}{G} \geq v\right)$$
$$= m P(\mathfrak{e} \geq v G) = m E(e^{-vG}).$$

Hence:

$$E\left(\exp\left\{-\lambda \int_0^\infty e^{-t} dY_t\right\}\right) = \exp\left\{-m \int_0^\infty (1 - e^{-\lambda v}) \frac{dv}{v} E(e^{-vG})\right\}$$
$$= E(e^{-\lambda \Gamma_m(G)})$$

*ii)* <u>We now prove point 2</u>. We have:

$$\int_0^\infty e^{-t} dY_t = \int_0^{T_1} e^{-t} dY_t + \int_{T_1}^\infty e^{-t} dY_t$$

(where $T_1$ is the first jump time of $(N_t, \ t \geq 0)$)

$$= e^{-T_1} K_1 + e^{-T_1} \int_0^\infty e^{-t} d(Y_{T_1 + t})$$

and we observe that:

$$e^{-T_1} \stackrel{(law)}{=} U^{1/m}, \quad \int_0^\infty e^{-t} d(Y_{T_1 + t}) \stackrel{(law)}{=} \int_0^\infty e^{-t} dY_t,$$

and $\int_0^\infty e^{-t} d(Y_{T_1+t})$ is independent of $T_1$.

### 1.5. The subordinators $(\Gamma_t(G), \ t \geq 0)$

Until now, we have been interested uniquely in the "individual" GGC variables. However, to each GGC r.v. $\Gamma$ we may, since $\Gamma$ is infinitely divisible and positive, associate a unique subordinator $(\Gamma_t, \ t \geq 0)$ such that $\Gamma_1 \stackrel{(law)}{=} \Gamma$. It is this subordinator which we shall now define and describe.

**1.5.a** Let $\Gamma$ denote a GGC r.v. From Lévy-Khintchine formula, there exists a (unique) subordinator $(\Gamma_t,\ t \geq 0)$ such that:

$$\Gamma_1 \stackrel{(law)}{=} \Gamma \tag{49}$$

hence we have:

$$E(e^{-\lambda \Gamma_t}) = \exp\left\{-t \int_0^\infty (1 - e^{-\lambda x}) \frac{dx}{x} \int_0^\infty e^{-xz} \mu(dx)\right\} \tag{50}$$

where $\mu$ denotes the Thorin measure associated with $\Gamma$.

**1.5.b** Let $G$ denote a positive r.v. which satisfies (19). Then, there exists a subordinator $(\Gamma_t(G),\ t \geq 0)$ which is characterized by:

$$E(e^{-\lambda \Gamma_t(G)}) = \exp\left\{-t \int_0^\infty (1 - e^{-\lambda x}) \frac{dx}{x} E(e^{-xG})\right\} \tag{51}$$

$$= \exp\left\{-t E\left(\log\left(1 + \frac{\lambda}{G}\right)\right)\right\} \tag{52}$$

In particular, for every $t > 0$, $\Gamma_t(G)$ is a $(t, G)$ GGC r.v. Thus, there exists, following (33), a family of r.v.'s $D_t(G),\ t \geq 0$, whose laws are characterized by:

$$\Gamma_t(G) = \gamma_t \cdot D_t(G) \tag{53}$$

and, for every $t > 0$, from (34) and Proposition 1.3:

$$\Gamma_t(G) \stackrel{(law)}{=} \int_0^t \frac{1}{F_G^{-1}\left(\frac{u}{t}\right)} d\gamma_u \stackrel{(law)}{=} \int_0^t \frac{1}{F_G^{-1}(1 - \frac{u}{t})} d\gamma_u \tag{54}$$

$$D_t(G) \stackrel{(law)}{=} \int_0^t \frac{1}{F_G^{-1}\left(\frac{u}{t}\right)} d_u(D_u^{(t)}) \stackrel{(law)}{=} \int_0^t \frac{1}{F_G^{-1}\left(\frac{u}{t}\right)} d_u(D_u^{(t)}) \tag{55}$$

We note that the relations (54) and (55) are only true for fixed $t$, for any $t > 0$, but do not hold as equalities in law between processes. On the other hand, since: $E(e^{-\lambda \gamma_t}) = \frac{1}{(1+\lambda)^t}$, we deduce from (53) that:

$$E(e^{-\lambda \Gamma_t(G)}) = E\left(\frac{1}{\left(1 + \lambda D_t(G)\right)^t}\right) \tag{56}$$

**1.5.c** Some elementary properties of $(\Gamma_t(G),\ t \geq 0)$ and $(D_t(G),\ t \geq 0)$
Let $G$ denote a positive r.v. which satisfies (19).

**Proposition 1.7**
**1)** *The family of laws of the r.v.'s $D_t(G),\ t \geq 0$ solves the equation: for every $t, s \geq 0$:*

$$\gamma_{t+s} \cdot X_{t+s} \stackrel{(law)}{=} \gamma_t \cdot X_t + \gamma_s X_s \tag{57}$$

**2)** The family of laws of the r.v.'s $\big(D_t(G),\ t\geq 0\big)$ solves the equation: for every $s,t\geq 0$:
$$X_{t+s}(G) \stackrel{(law)}{=} \beta_{t,s}X_t + (1-\beta_{t,s})X_s \qquad (58)$$

**3)** We assume furthermore that $E\big(\frac{1}{G}\big) < \infty$ $\big($which implies $E\left(\log^+\left(\frac{1}{G}\right)\right) < \infty\big)$. Then:

$3\,i)$    for every $t>0$,    $E\big(\Gamma_t(G)\big) = t\,E(1/G)$ and $E\big(D_t(G)\big) = E(1/G)$  (59)

$3\,ii)$    $\dfrac{1}{t}\Gamma_t(G) \xrightarrow[t\to 0]{a.s.} E(1/G)$   and   $D_t(G) \xrightarrow[t\to\infty]{L^1} E(1/G)$  (60)

$3\,iii)$    $\big[\Gamma_t(G)\big]^t \xrightarrow[t\to 0]{(law)} U$   and   $\big[D_t(G)\big]^t \xrightarrow[t\to 0]{(law)} 1$  (61)

where $U$ is uniform on $[0,1]$.

Point 2 of this Proposition 1.7 is due to Hjort and Ongaro (see [27]), which can be seen as a consequence of Ethier and Griffiths [19, Lemma 1].

We note the following remarkable feature of points 1 and 2 of Proposition 1.7: the affine equations (where the unknowns are the laws of the $(X_t,\ t\geq 0)$:

$$\gamma_{t+s} \cdot X_{t+s} \stackrel{(law)}{=} \gamma_t \cdot X_t + \gamma_s \cdot X_s \qquad (s_i t \geq 0)$$
$$\text{and}\quad X_{t+s} \stackrel{(law)}{=} \beta_{t,s}X_t + (1-\beta_{t+s})X_s \qquad (s_i t \geq 0)$$

both admit infinitely many solutions: $X_u = D_u(G)$, $u \geq 0$ for every r.v. G which satisfy (1.19).

**Proof of Proposition 1.7**
<u>Point 1</u>:
$$\gamma_{t+s} \cdot D_{t+s}(G) \stackrel{(law)}{=} \gamma_t\, D_t(G) + \gamma_s\, D_s(G) \qquad (62)$$

follows from:
$$\Gamma_{t+s}(G) \stackrel{(law)}{=} \Gamma_t(G) + \Gamma_s(G)$$

since $\big(\Gamma_t(G),\ t\geq 0\big)$ is a subordinator and since, from the definition of $D_a(G)$:
$$\Gamma_a(G) \stackrel{(law)}{=} \gamma_a \cdot D_a(G).$$

We now show (58). From the beta-gamma algebra, we have:
$$(\gamma_t, \gamma_s) \stackrel{(law)}{=} \big(\beta_{t,s}\,\gamma_{t+s},\ (1-\beta_{t,s})\gamma_{t+s}\big) \qquad (63)$$

hence, plugging (63) in (62), we obtain:
$$\gamma_{t+s} \cdot D_{t+s}(G) \stackrel{(law)}{=} \beta_{t,s}\,\gamma_{t+s}\, D_t(G) + (1-\beta_{t,s})\gamma_{t+s}D_s(G)$$

which, using point **1.4.b**, implies relation (58). <u>We now prove 3 $i$)</u>. We have:

$$E[\Gamma_t(G)] = -\frac{\partial}{\partial \lambda} E(e^{-\lambda \Gamma_t})\big|_{\lambda=0} = -\frac{\partial}{\partial \lambda} \exp\left\{-t E\left(\log\left(1+\frac{\lambda}{G}\right)\right)\right\}\big|_{\lambda=0}$$
$$= t E(1/G)$$

and we deduce from (53), that:

$$E[\Gamma_t(G)] = E[\gamma_t \cdot D_t(G)] = t E(D_t(G)) = t E(1/G)$$

<u>We now prove 3 $ii$)</u>. The a.s. convergence (hence, the convergence in law) of $\frac{1}{t}\Gamma_t(G)$ as $t \to \infty$ towards $E(1/G)$ follows from the law of large numbers.

We prove that $D_t(G) \xrightarrow[t\to\infty]{L^1} E(1/G)$. Indeed, we may write:

$$\frac{\Gamma_t(G)}{t} \stackrel{(law)}{=} D_t(G)\left(\frac{\gamma_t}{t} - 1\right) + D_t(G) \tag{64}$$

and we shall prove that:

$$D_t(G)\left|\frac{\gamma_t}{t} - 1\right| \xrightarrow[t\to\infty]{L^1} 0$$

which implies that $D_t(G) \xrightarrow[t\to\infty]{P} E(1/G)$. Since $E(D_t(G)) = E\left(\frac{1}{G}\right)$, this implies that the family of laws of $\{D_t(G), t > 0\}$ is tight, hence convergence in probability implies convergence in $L^1$. It remains to prove: $E\left(D_t(G)\left|\frac{\gamma_t}{t} - 1\right|\right) \xrightarrow[t\to\infty]{} 0$. But we have:

$$E\left(D_t(G)\left|\frac{\gamma_t}{t} - 1\right|\right) = E\left(\frac{1}{G}\right) E\left(\left|\frac{\gamma_t}{t} - 1\right|\right)$$
$$\leq E\left(\frac{1}{G}\right)\left(E\left(\left(\frac{\gamma_t}{t}\right)^2 - 1\right)\right)^{\frac{1}{2}} = E\left(\frac{1}{G}\right) \frac{\sqrt{t+1}}{t} \xrightarrow[t\to\infty]{} 0$$

Hence the result.

<u>We now prove point 3 $iii$)</u>.

If we knew that $\bigl(D_t(G)\bigr)^t \xrightarrow[t\to 0]{(law)} 1$, then, combining this result with the classical one: $(\gamma_t)^t \xrightarrow[t\to 0]{(law)} U$, we would deduce:

$$\bigl(\Gamma_t(G)\bigr)^t \stackrel{(law)}{=} (\gamma_t)^t \cdot \bigl(D_t(G)\bigr)^t \xrightarrow[t\to 0]{(law)} U$$

In fact, we shall proceed in the other direction. We shall show further (see point 7 of Remark 2.2) that $\bigl(\Gamma_t(G)\bigr) \xrightarrow[t\to 0]{(law)} U$. The relation $\bigl(\Gamma_t(G)\bigr)^t = (\gamma_t)^t \bigl(D_t(G)\bigr)^t$ and:

$$(\gamma_t)^t \xrightarrow[t\to 0]{(law)} U, \quad \bigl(\Gamma_t(G)\bigr)^t \xrightarrow[t\to 0]{(law)} U \quad \text{then imply easily that:}$$

$$\bigl(D_t(G)\bigr)^t \xrightarrow[t\to 0]{(law)} 1.$$

**Remark 1.7** Lijoi and Reggazini [36] have shown that the support of the law of $D_t(G)$ is the closure of the convex hull of the support of the law of $1/G$. In particular:

$$\bullet \quad \text{If } G \leq a \text{ a.s., then } D_t(G) \geq \frac{1}{a} \text{ a.s.} \tag{65}$$

$$\bullet \quad \text{If } G \geq a \text{ a.s., then } D_t(G) \leq \frac{1}{a} \text{ a.s.} \tag{66}$$

In Section 3 of this survey, we shall verify this assertion on inspection of numerous examples.

### 1.6. Some examples of GGC subordinators

Let $\Gamma$ denote a GGC variable with associated Thorin measure $\mu$, and let $(\Gamma_t, t \geq 0)$ denote the subordinator such that $\Gamma_1 \stackrel{(law)}{=} \Gamma$. We have:

$$E(e^{-\lambda \Gamma_t}) = \exp\left\{ -t \int_0^\infty (1 - e^{-\lambda x}) \frac{dx}{x} \int_0^\infty e^{-xz} \mu(dz) \right\} \tag{67}$$

Such a subordinator is called a GGC subordinator.

Here are now some examples of such subordinators. They are lifted from a paper by H. Matsumoto, L. Nguyen and M. Yor [38] on one hand, and from the study of hyperbolic subordinators made in Pitman and Yor [41] on the other hand. The reader may refer to these papers for further information. In the study of these examples, we shall denote subordinators with curly letters, especially to avoid some possible confusion with the modified Bessel functions, which are traditionally written with ordinary capital letters.

**1.6.a** The hyperbolic subordinators (see [41] for a probabilistic description of these subordinators.)

*i)* The subordinator $(\mathcal{C}_t, t \geq 0)$ is characterized by:

$$E(e^{-\lambda \mathcal{C}_t}) = \exp\{-t \log \cosh \sqrt{2\lambda}\} = \left(\frac{1}{\cosh \sqrt{2\lambda}}\right)^t \qquad (\lambda, t \geq 0) \tag{68}$$

Its Lévy density $l_\mathcal{C}$ equals:

$$l_\mathcal{C}(x) = \frac{1}{x} \sum_{n=1}^\infty \exp\left(-\frac{\pi^2}{8}(2n-1)^2 x\right) \qquad (x \geq 0) \tag{69}$$

Hence, its associated Thorin measure, (i.e.: the Thorin measure of $\mathcal{C}_1$) equals:

$$\mu_\mathcal{C}(dx) = \sum_{n=1}^\infty \delta_{\frac{\pi^2}{8}(2n-1)^2}(dx) \tag{70}$$

It has infinite total mass.

*ii)* The subordinator $(\mathcal{S}_t,\ t \geq 0)$ is characterized by:

$$E(e^{-\lambda \mathcal{S}_t}) = \Big(\frac{\sqrt{2\lambda}}{\sinh\sqrt{2\lambda}}\Big)^t \qquad (\lambda, t \geq 0) \tag{71}$$

Its Lévy density equals:

$$l_\mathcal{S}(x) = \frac{1}{x} \sum_{n=1}^\infty \exp\Big(-\frac{\pi^2}{2} n^2 x\Big) \qquad (x \geq 0) \tag{72}$$

Hence, its Thorin measure equals:

$$\mu_\mathcal{S}(dx) = \sum_{n=1}^\infty \delta_{\frac{\pi^2 n^2}{2}}(dx) \qquad (x \geq 0) \tag{73}$$

It has infinite total mass.

We note that the subordinator $(\mathcal{T}_t,\ t \geq 0)$ which is characterized by:

$$E(e^{-\lambda \mathcal{T}_t}) = \Big(\frac{1}{\cosh\sqrt{2\lambda}}\Big)^t \cdot \Big(\frac{\sinh\sqrt{2\lambda}}{\sqrt{2\lambda}}\Big)^t = \Big(\frac{\tanh\sqrt{2\lambda}}{\sqrt{2\lambda}}\Big)^t \tag{74}$$

satisfies:

$$(\mathcal{C}_t,\ t \geq 0) \stackrel{(law)}{=} (\mathcal{S}_t + \mathcal{T}_t,\ t \geq 0) \tag{75}$$

and that its Lévy density equals:

$$l_\mathcal{T}(x) = l_\mathcal{C}(x) - l_\mathcal{S}(x)$$

However, this subordinator $(\mathcal{T}_t\ t \geq 0)$ is not GGC, as its would be 'Thorin measure' $\mu_\mathcal{T}$ is a signed measure:

$$\mu_\mathcal{T} = \mu_\mathcal{C} - \mu_\mathcal{S}.$$

**1.6.b** The subordinators $(\mathcal{J}_t^{(0)},\ t \geq 0)$ and $(\mathcal{K}_t^{(0)},\ t \geq 0)$

We denote by $I_\nu$ and $K_\nu$ the modified Bessel functions with index $\nu$ (see [34], p. 108).

*i)* The subordinator $(\mathcal{J}_t^{(0)},\ t \geq 0)$ is characterized by:

$$E(e^{-\lambda \mathcal{J}_t^{(0)}}) = \big(1+\lambda+\sqrt{(1+\lambda)^2-1}\big)^{-t} = \exp\Big\{-t\int_0^\infty (1-e^{-\lambda x})\frac{dx}{x} I_0(x)e^{-x}\Big\} \tag{76}$$

Its Lévy density $l_{\mathcal{J}^{(0)}}$ equals:

$$l_{\mathcal{J}^{(0)}}(x) = I_0(x)\frac{e^{-x}}{x} \qquad (x \geq 0) \tag{77}$$

and its Thorin measure $\mu_{\mathcal{J}^{(0)}}$, with total mass equal to 1, equals:

$$\mu_{\mathcal{J}^{(0)}}(dx) = \frac{1}{\pi} \frac{dx}{\sqrt{x(2-x)}} 1_{[0,2]}(x) \tag{78}$$

This formula was obtained by Bondesson ([7], Ex. 5.2). In Section 3, Theorem 3.1, we shall meet again this subordinator $(\mathcal{J}_t^{(0)},\ t \geq 0)$.

*ii)* The subordinator $(\mathcal{K}_t^{(0)},\ t \geq 0)$ is characterized by:

$$\begin{aligned} E(e^{-\lambda \mathcal{K}_t^{(0)}}) &= \exp\Big\{-\frac{t}{2}\big(\operatorname{arg cosh}(1+\lambda)\big)^2\Big\} \\ &= \exp\Big\{-t\int_0^\infty (1-e^{-\lambda x})\frac{dx}{x} K_0(x) e^{-x}\Big\} \end{aligned} \quad (79)$$

Hence, its Lévy density equals:

$$l_{\mathcal{K}^{(0)}}(x) = K_0(x)\frac{e^{-x}}{x} \qquad x \geq 0 \quad (80)$$

and its Thorin measure $\mu_{\mathcal{K}^{(0)}}$ equals:

$$\mu_{\mathcal{K}^{(0)}}(dx) = \frac{dx}{\sqrt{x(x-2)}}\, 1_{[2,\infty[}(x)dx \quad (81)$$

We note that this latter formula follows from:

$$\begin{aligned} e^{-x}K_0(x) &= e^{-x}\int_1^\infty e^{-x(\cosh t)}dt \qquad \text{(see [34], p. 119)} \\ &= e^{-x}\int_1^\infty e^{-xu}\frac{du}{\sqrt{u^2-1}} = \int_2^\infty e^{-xv}\frac{dv}{\sqrt{v(v-2)}} \end{aligned}$$

*iii)* The subordinators $(\mathcal{K}_t^{(0)},\ t \geq 0)$ and $(\mathcal{J}_t^{(0)},\ t \geq 0)$ are connected via the subordination relation:

$$(\mathcal{J}_t^{(0)},\ t \geq 0) \stackrel{(law)}{=} (\mathcal{K}_{S_{1/2}(t)}^{(0)},\ t \geq 0) \quad (82)$$

where, on the RHS from (82), the processes $(\mathcal{K}_u^{(0)},\ u \geq 0)$ and $(S_{1/2}(t),\ t \geq 0)$ are independent, and where $(S_{1/2}(t),\ t \geq 0)$ is the stable subordinator with index 1/2 characterized by:

$$E(e^{-\lambda S_{1/2}(t)}) = \exp\big(-t\sqrt{2\lambda}\big) \quad (83)$$

*iv)* The v.a. $\mathcal{K}_t^{(0)}$, for $t$ fixed, may be realized in the following manner: let $\big(b_u(s)\,;\ 0 \leq s \leq u\big)$ the brownian bridge with length $u$ (with $b_u(0) = b_u(u) = 0$) and let:

$$A(b_u) := \int_0^u \exp\big(2b_u(s)\big)ds$$

Then:

$$\mathcal{K}_t^{(0)} \stackrel{(law)}{=} \big[A(b_{1/t})\big]^{-1} \stackrel{(law)}{=} t\Big(\int_0^1 \exp\Big(\frac{2b(s)}{\sqrt{t}}\Big)ds\Big)^{-1} \quad (84)$$

where, in (84), $\big(b(s),\ 0 \leq s \leq 1\big)$ denotes the standard brownian bridge $b_1$.

*v)* In relation with the preceding, D. Dufresne and M. Yor [15], in a work in preparation [15], establish the following formula: let $\bigl(b_u^{(x)}(s),\ 0 \leq s \leq u\bigr)$ denote the brownian bridge with length $u$, starting from 0 and such that $b_u^{(x)}(u) = x$. Let:

$$A\bigl[b_u^{(x)}\bigr] := \int_0^u \exp\bigl(2b_u^{(x)}(s)\bigr)ds \tag{85}$$

Then, D. Dufresne and M. Yor obtain the representation as a Wiener-Gamma integral of the r.v. $\frac{1}{A[b_u^{(x)}]}$:

$$\frac{1}{A[b_u^{(x)}]} \stackrel{(law)}{=} \int_0^\infty \frac{2}{e^{2x} + 2e^x \cosh(su) + 1}\, d\gamma_s \tag{86}$$

and in particular:

$$\frac{1}{A[b_u]} \stackrel{(law)}{=} \int_0^\infty \frac{d\gamma_s}{1 + \cosh(su)} \tag{87}$$

**1.6.c** The subordinators $(\mathcal{J}_t^{(\nu)},\ t \geq 0)$ and $(\mathcal{K}_t^{(\nu)},\ t \geq 0)$

They are obtained by replacing in (76), resp. (79), $I_0$ by $I_\nu$, resp. $K_0$ by $K_\nu$.

*i)* The subordinator $(\mathcal{J}_t^{(\nu)},\ t \geq 0)$, which is defined for $\nu > -1$, is characterized by:

$$E(e^{-\lambda \mathcal{J}_t^{(\nu)}}) = \exp\Bigl\{ -t \int_0^\infty (1 - e^{-\lambda x})\, \frac{dx}{x}\, e^{-x} I_\nu(x) \Bigr\} \tag{88}$$

Thus, its Lévy density $l_{\mathcal{J}^{(\nu)}}$ equals:

$$l_{\mathcal{J}^{(\nu)}}(x) = \frac{e^{-x}}{x} I_\nu(x) \qquad x \geq 0 \tag{89}$$

Its Thorin measure, in the case $-1/2 < \nu < 0$, is obtained by the following calculus:

$$I_\nu(z) = \frac{\bigl(\frac{z}{2}\bigr)^\nu}{\sqrt{\pi}\,\Gamma(\nu + 1/2)} \cdot \int_0^1 (1-t^2)^{\nu-1/2}(e^{zt} + e^{-zt})dt \quad (Re\,\nu > -1/2)$$

(see [34], p. 119, form.5.10.22)

and:

$$e^{-z}I_\nu(z) = C_\nu\, z^\nu \int_0^1 (1-t^2)^{\nu-1/2}(e^{-z(1-t)} + e^{-z(1+t)})dt$$

But, for $\nu < 0$:

$$z^\nu e^{-za} = \frac{1}{\Gamma(-\nu)} \int_0^\infty s^{-\nu-1} e^{-z(a+s)}ds$$

Hence:

$$\begin{aligned}
e^{-z}I_\nu(z) &= C'_\nu \int_0^1 (1-t^2)^{\nu-1/2} dt \int_0^\infty s^{-\nu-1}(e^{-z(1-t+s)} + e^{-z(1+t+s)}) ds \\
&= C'_\nu \int_0^1 (1-t^2)^{\nu-1/2} dt \Big[\int_{1-t}^\infty e^{-zh}(h-1+t)^{-\nu-1} dh \\
&\quad + \int_{1+t}^\infty e^{-zh}(h-1-t)^{-\nu-1} dh\Big] \\
&= C'_\nu \int_0^1 e^{-zh} dh \int_0^1 dt (1-t^2)^{\nu-1/2}\big[(h-1+t)^{-\nu-1}\mathbf{1}_{h>1-t} \\
&\quad + (h-1-t)^{-\nu-1}\mathbf{1}_{h>1+t}\big]
\end{aligned}$$

Hence, in the case $-1/2 < \nu < 0$, the Thorin measure of $(\mathcal{J}_t^{(\nu)}, t \geq 0)$ admits a density equal to:

$$\mu_{\mathcal{J}^{(\nu)}}(dx)/dx = \Big[C'_\nu \int_0^1 (1-t^2)^{\nu-1/2} \big\{\big(x-(1-t)\big)^{-\nu-1}\mathbf{1}_{x>1-t} \\
+ \big(x-(1+t)\big)^{-\nu-1}\mathbf{1}_{x>1+t}\big\} dt\Big]$$

Hence, if $-1/2 < \nu < 0$, $(\mathcal{J}_t^\nu, t \geq 0)$ is a GGC subordinator.

*ii)* The subordinator $(\mathcal{K}_t^{(\nu)}, t \geq 0)$, defined for $\nu < 1$, is characterized by:

$$E(e^{-\lambda \mathcal{K}_t^{(\nu)}}) = \exp\Big\{-t \int_0^\infty (1-e^{-\lambda x}) \frac{dx}{x} e^{-x} K_\nu(x)\Big\} \tag{90}$$

Its Lévy density $l_{\mathcal{K}^{(\nu)}}$ equals:

$$l_{\mathcal{K}^{(\nu)}}(x) = \frac{e^{-x}}{x} K_\nu(x) \qquad (x \geq 0) \tag{91}$$

Its Thorin measure $\mu_{\mathcal{K}^{(\nu)}}$, obtained by using the formula ([34], p. 119)

$$K_\nu(z) = \int_0^\infty e^{-z\cosh u}(\cosh \nu u) du$$

equals:

$$\mu_{\mathcal{K}^{(\nu)}}(dx) = \cosh\big[\nu \arg\cosh(x-1)\big] \frac{dx}{\sqrt{x(x-2)}} \mathbf{1}_{[2,\infty[}(x) \tag{92}$$

The formula in (92) simplifies using,

$$\cosh\big[\nu \arg\cosh(x-1)\big] = \frac{1}{2}\Big\{\big((x-1)+\sqrt{(x-2)x}\big)^\nu + \big((x-1)+\sqrt{(x-2)x}\big)^{-\nu}\Big\},$$

since:

$$\arg\cosh(y) = \log\big(y + \sqrt{y^2-1}\big).$$

### 1.7. The subordinator $(\Gamma_t(S_\beta),\ t \geq 0)$ $(0 < \beta < 1)$

We give here a realization for the subordinator $(\Gamma_t(S_\beta),\ t \geq 0)$, where $S_\beta$ is a positive $\beta$-stable r.v.

Let $\alpha > 0$. We denote $\sigma_{t,n} := \{s_1 = 0 < s_2 \cdots < s_n = t\}$ a subdivision of the interval $[0,t]$ with mesh $\delta(\sigma_{t,n}) := \sup_i (s_{i+1} - s_i)$. The limit along the decreasing filtering set of sequences $\sigma_{t,n}$ whose meshes $\delta(\sigma_{t,n})$ tend to zero, of $\sum_{s_k \in \sigma_{t,n}} (\gamma_{s_{k+1}} - \gamma_{s_k})^\alpha$ exists a.s. Observe that: $E\big(\sum_{s_k \in \sigma_{t,n}} (\gamma_{s_{k+1}} - \gamma_{s_k})^\alpha\big) = \sum_{s_k \in \sigma_{t,n}} \frac{\Gamma(\alpha + s_{k+1} - s_k)}{\Gamma(s_{k+1} - s_k)} \xrightarrow[\delta(\sigma_{t,n}) \to 0]{} t\,\Gamma(\alpha)$. Let:

$$V^{(\alpha)}(t) := \lim_{\delta(\sigma_{t,n}) \to 0} \sum_{s_k \in \sigma_{t,n}} (\gamma_{s_{k+1}} - \gamma_{s_k})^\alpha = \sum_{0 < s \leq t} (\gamma_s - \gamma_{s^-})^\alpha \qquad (93)$$

Let

$$N_{t,x} := \#\{s \leq t\ ;\ \gamma_s - \gamma_{s^-} > x^{1/\alpha}\}$$

From Section 1.4, $N_{t,x}$ is a Poisson r.v. with parameter $t \cdot \int_{x^{1/\alpha}}^\infty \frac{e^{-u}}{u}\, du$. But, evidently:

$$N_{t,x} = \#\{s \leq t\ ;\ (\gamma_s - \gamma_{s^-})^\alpha > x\} \qquad (94)$$

We deduce from (94) that the Poisson point process

$$\{(\gamma_s - \gamma_{s^-})^\alpha,\ s \geq 0\}\ \text{admits}\ n^{(\alpha)}(dx) = \frac{1}{\alpha}\frac{dx}{x}\, e^{-x^{1/\alpha}} \qquad (x > 0)$$

as intensity measure. Hence $(V^{(\alpha)}(t),\ t \geq 0)$ is a subordinator whose Lévy measure equals $n^{(\alpha)}(dx)$:

$$E(e^{-\lambda V^{(\alpha)}(t)}) = \exp\Big\{ - \frac{t}{\alpha} \int_0^\infty (1 - e^{-\lambda x})\, \frac{dx}{x}\, e^{-x^{1/\alpha}} \Big\}$$

This subordinator, for every $\alpha > 0$, is self-decomposable but it is GGC only in the case $\alpha \geq 1$. In this case: $\alpha \geq 1$ we have:

$$E(e^{-\lambda V^{(\alpha)}(t)}) = \exp\Big( - \frac{t}{\alpha} \int_0^\infty (1 - e^{-\lambda x})\frac{dx}{x}\, E(e^{-x S_{1/\alpha}}) \Big\}$$

where $S_{1/\alpha}$ is a positive $\frac{1}{\alpha}$ stable r.v. In other terms, with the notation of Section **1.5.b**:

$$(\Gamma_t(S_\beta),\ t \geq 0) \stackrel{(law)}{=} (V^{(\alpha)}(t\alpha),\ t \geq 0) \qquad \Big(\alpha = \frac{1}{\beta} \geq 1\Big) \qquad (95)$$

### 1.8. Remark

We now end up this Introduction by indicating how the study of the GGC subordinators may be embedded in a more general one.

A subordinator $(N_t,\ t \geq 0)$ is said to belong to the Thorin class $T^{(\chi)}(\mathbb{R}_+)$, with $\chi > 0$ (see [Grig]), if its Lévy measure admits a density $l_N$ of the form:

$$l_N(x) = x^{\chi-2} k(x) \qquad x > 0 \tag{96}$$

where $k$ is a completely monotonic function, i.e. it may be represented as:

$$k(x) := \int_0^\infty e^{-xy}\, \mu(dy) \tag{97}$$

for a positive Radon measure $\mu$ carried by $\mathbb{R}_+^*$. Thus, the subordinators which we study in this work, i.e.: the GGC subordinators, belong to the class $T^{(1)}(\mathbb{R}_+)$. The class $T^{(2)}(\mathbb{R}_+)$ has been studied by Goldie [24], Steutel [47] and Bondesson [6]. The r.v's which belong to this class are the generalized convolutions of mixtures of exponential laws. Note that $T^{(\chi)}(\mathbb{R}_+) \subset T^{(\chi')}(\mathbb{R}_+)$ if $\chi < \chi'$. We also note that [BNMS] present extensions of these notions to $\mathbb{R}^d$.

In the same manner as condition (9) is necessary and sufficient for a measure $\mu$ to be the Thorin measure associated to a positive r.v., B. Grigelionis [25] obtains a necessary and sufficient analytical condition so that a measure $\mu$ defines, via (96) and (97), a subordinator $(N_t,\ t \geq 0)$ which belongs to the class $T^{(\chi)}(\mathbb{R}_+)$.

**1.9** We now detail the contents of the sequel of this paper:

• In Section 2, we present a duality result which connects on one hand the r.v.'s $\Gamma_t(G)$ and $D_t(G)$ to the r.v.'s $\Gamma_t(1/G)$ and $\Gamma_t(G)$ on the other hand.

• In Section 3, we study in depth the examples of subordinators $\Gamma_t(\mathbb{G}_\alpha)$ ($0 \leq \alpha \leq 1$) for which we know how to compute explicitly their Laplace transforms, i.e.: their Lévy exponents, as well as their densities at any time, and their Wiener-Gamma representations.

## 2. A duality principle

Throughout this section, $G$ denotes a positive r.v. such that $E(|\log G|) < \infty$. Thus, we have:

$$E\bigl(\log^+(1/G)\bigr) < \infty \quad \text{and} \quad E\bigl(\log^+(G)\bigr) < \infty \tag{98}$$

Consequently, the subordinators $\bigl(\Gamma_t(G),\ t \geq 0\bigr)$ and $\bigl(\Gamma_t(1/G),\ t \geq 0\bigr)$ are well defined. We denote by $\psi_G$ (resp. $\psi_{1/G}$) the characteristic exponent (i.e.: the Bernstein function or Lévy exponent) of the subordinator $\bigl(\Gamma_t(G),\ t \geq 0\bigr)$, resp. $\bigl(\Gamma_t(1/G),\ t \geq 0\bigr)$.

$$\psi_G(\lambda) := \int_0^\infty (1 - e^{-\lambda x})\, \frac{dx}{x}\, E(e^{-xG}) \tag{99}$$

and the same formula for $\psi_{1/G}$ obtained when replacing $G$ by $1/G$. We denote by $F_G$ the cumulative distribution function of $G$ and by $F_G^{-1}$ its right continuous inverse, in the sense of the composition of functions.

### 2.1. The duality theorem

**Theorem 2.1** (Duality)

$$1) \quad F_{1/G}^{-1}(y) \cdot F_G^{-1}(1-y) = 1 \quad a.e. \quad (y \in [0,1]) \tag{100}$$

*2)* For any $\lambda \geq 0$:

$$\psi_{1/G}(\lambda) - \psi_G(1/\lambda) = E(\log G) + \log \lambda \tag{101}$$

*3) i)* $f_{\Gamma_t(1/G)}$ the density of the r.v. $\Gamma_t(1/G)$, equals:

$$f_{\Gamma_t(1/G)}(x) = e^{-t E(\log G)} E\left\{ \left(\frac{\Gamma_t(G)}{x}\right)^{\frac{1-t}{2}} J_{t-1}\left(2\sqrt{x\,\Gamma_t(G)}\right) \right\} \tag{102}$$

where $J_\nu$ denotes the Bessel function with index $\nu$:

$$J_\nu(z) = \sum_{k=0}^{\infty} (-1)^k \frac{1}{\Gamma(k+1)\,\Gamma(k+\nu+1)} \left(\frac{z}{2}\right)^{\nu+2k}, \quad |z|<\infty,\ |\operatorname{Arg} z|<\pi \tag{103}$$

(see [34], p. 102)

$$ii)\ f_{\Gamma_t(1/G)}(x) = \frac{x^{t-1}}{\Gamma(t)}\, e^{-t E(\log G)}\, E(e^{-x D_t(G)}) \quad (x>0) \tag{104}$$

$$iii)\, E\left(\frac{1}{\left(D_t(G)\right)^t}\right) = e^{t E(\log G)} \tag{105}$$

*4) i)* $$E\left(\frac{1}{\left(\lambda + D_t(G)\right)^t}\right) = \exp\left\{ -t\left(\psi_{1/G}(\lambda) - E(\log G)\right) \right\} \tag{106}$$

*ii)* The densities $f_{D_t(G)}$ and $f_{D_t(1/G)}$ of, resp., $D_t(G)$ and $D_t(1/G)$ are related by:

$$f_{D_t(1/G)}(x) = x^{t-2}\, e^{-t E(\log G)}\, f_{D_t(G)}(1/x) \quad (x>0) \tag{107}$$

**Remark 2.2.**
**1)** The relation (101) is found, under a slightly different form, in Bondesson ([5], p. 48, III "Curious composition"), under the only hypothesis $E(\log^+(1/G)) < \infty$. We have chosen to work under the stronger hypothesis $E(|\log G|) < \infty$ since it seems that only under this hypothesis can we obtain the most interesting results.

**2)** Of course, all formulae of Theorem 2.1 are "involutive", in that $G$ may be replaced there by $1/G$.

**3)** Formula (104) agrees with, and makes more precise, the result of Bondesson which we recalled in (28) and (30). In particular:

$$f_{\Gamma_t(G)}(x) \underset{x \to 0}{\sim} \frac{x^{t-1}}{\Gamma(t)}\, e^{t E(\log G)} \tag{108}$$

**4)** Formula (101) may be generalized as follows:
Let $a, b, c, d \in \mathbb{R}^4$, with $ad - bc = \pm 1$ and let, for $x, \lambda \geq 0$:

$$\sigma(x) = \frac{ax + b}{cx + d}, \quad \widetilde{\sigma}(\lambda) = \frac{d\lambda + b}{c\lambda + a} \tag{109}$$

so that $\sigma(G)$ be a positive r.v. Then:

$$-\psi_G(\widetilde{\sigma}(\lambda)) + \psi_{\sigma(G)}(\lambda) = \log(c\lambda + a) + k$$

with $k := E\big(\log\big(\frac{G+\lambda_0}{\sigma(G)+\lambda_0} \cdot \frac{\sigma(G)}{G}\big)\big) - \log(c\lambda_0 + a)$ where $\lambda_0$ is a fixed point of $\widetilde{\sigma}$. The relation (101) corresponds to $a = d = 0$, $b = c = 1$. We shall study, in the Appendix, the case where $a = d = \sinh u$, $b = c = \cosh u$ $(u \geq 0)$.

**5)** Cifarelli and Regazzini [11], Cifarelli and Melilli [9] have obtained the density of $D_t(G)$ for $t \geq 1$ and James, Lijoi and Prünster [30] have obtained it for $t \leq 1$. For $t \leq 1$, they obtained:

$$f_{D_t(G)}(x) = \int_0^x (x-u)^{t-1} \frac{d}{du}\big(\theta_t(u)\big) du \tag{110}$$

with

$$\theta_t(u) := \frac{1}{\pi} \sin\big(t\pi F_{1/G}(u)\big) \exp\Big\{-tE\Big(\log\Big(u - \frac{1}{G}\Big) 1_{u \neq \frac{1}{G}}\Big)\Big\} \tag{111}$$

The proof of this formula (110) hinges upon the knowledge of the density of the r.v. $D_1\big(\frac{1}{G Y_t}\big)$ which is defined by: (see (41) and (42) with $m = t$ and $m' = 1$).

$$\beta_{t,1-t} \cdot D_t(1/G) \stackrel{(law)}{=} D_1\Big(\frac{1}{G \cdot Y_t}\Big) \qquad (t \leq 1)$$

This density equals:

$$f_{D_1}\Big(\frac{1}{G \cdot Y_t}\Big)(x) = \frac{\sin\big(t\pi F_G(1/x)\big)}{\pi} x^{t-1} \exp\Big\{-tE\Big(\log\Big(x - \frac{1}{G}\Big) \cdot 1_{G \neq \frac{1}{x}}\Big)\Big\} \tag{112}$$

which is obtained by inverting its Stieltjes transform.

Other formulae for densities of Dirichlet means may be found in Regazzini, Guglielmi and DiNunno [42].

**6)** Formula (105) is obtained from (104) by integrating between 0 and $\infty$:

$$1 = \int_0^\infty f_{\Gamma_t(1/G)}(x) dx = \int_0^\infty \frac{1}{\Gamma(t)} x^{t-1} e^{-tE(\log G)} E(e^{-x D_t(G)}) dx$$

$$= \frac{e^{-tE(\log G)}}{\Gamma(t)} E\Big[\frac{1}{(D_t(G))^t}\Big] \Gamma(t), \text{ hence}:$$

$$E\Big(\frac{1}{(D_t(G))^t}\Big) = e^{tE(\log G)} \tag{113}$$

The interest of this formula (113) is the following: it allows, in a situation where the law of $G$ is not known but when one knows the laws of $D_t(G)$ and $D_t(1/G)$ to show that $E(|\log G|) < \infty$ as soon as $E\big((D_t(G))^{-t}\big) < \infty$ and $E\big((D_t(1/G))^{-t}\big) < \infty$.

Formula (104), once multiplied by $x^\nu$, and integrated between 0 and $\infty$, leads to:

$$E\big[(\Gamma_t(1/G))^\nu\big] = \frac{\Gamma(\nu+t)}{\Gamma(t)} e^{-tE(\log G)} E\Big(\frac{1}{(D_t(G))^{\nu+t}}\Big) \tag{114}$$

The formula is also true even when the expectations which appear in this expression are infinite.

**7)** The result:

$$\big(\Gamma_t(G)\big)^t \xrightarrow[t\to 0]{(law)} U \tag{115}$$

with $U$ uniform on $[0,1]$ (see point 3 iii) of Proposition 1.7) is a consequence of (104). Indeed:

$$\begin{aligned}
P\big((\Gamma_t(G))^t \le e^{-a}\big) &= P\big(\Gamma_t(G) \le e^{-\frac{a}{t}}\big) \\
&= \int_0^{e^{-\frac{a}{t}}} f_{\Gamma_t(G)}(x) dx \\
&= \frac{e^{tE(\log G)}}{\Gamma(t)} \int_0^{e^{-\frac{a}{t}}} x^{t-1} E(e^{-xD_t(1/G)}) dx
\end{aligned}$$

(from (104) applied when replacing $G$ by $1/G$):

$$\underset{t\to 0}{\sim} \frac{1}{\Gamma(t)} \int_0^{e^{-\frac{a}{t}}} x^{t-1} dx = \frac{e^{-a}}{\Gamma(t+1)} \xrightarrow[t\downarrow 0]{} e^{-a}$$

We note that, from (64), the family of the laws of $D_t(1/G)$, $t \le 1$ is tight as soon as $E(G) < \infty$, since $E\big(D_t\big(\frac{1}{G}\big)\big) = E(G)$.

### 2.2. Proof of Theorem 2.1

**2.2.a** Point 1 is trivial. We show point 2:
Since, from Frullani's formula:

$$\begin{aligned}
\psi_G(\lambda) &= E\Big(\log\Big(1+\frac{\lambda}{G}\Big)\Big) = E\Big(\log\Big(\frac{\lambda}{G}\Big(1+\frac{G}{\lambda}\Big)\Big)\Big) \\
&= \log\lambda - E(\log G) + E\Big(\log\Big(1+\frac{G}{\lambda}\Big)\Big)
\end{aligned}$$

we have, by changing $\lambda$ in $1/\lambda$:

$$-\psi_G(1/\lambda) + \psi_{1/G}(\lambda) = -\log\lambda + E(\log G)$$

**2.2.b** We now show point 3 of Theorem 2.1:

*i)* From formula (101), after multiplying by $t$ and exponentiating, we obtain:

$$E(e^{-\lambda \Gamma_t(1/G)}) = \frac{e^{-tE(\log G)}}{\lambda^t} E(e^{-\frac{1}{\lambda} \Gamma_t(G)}) \tag{116}$$

Then, taking the Laplace transform of both sides of (116) in the variable $\lambda$, we obtain:

$$\begin{aligned}
\int_0^\infty E(e^{-\lambda(\Gamma_t(1/G))})e^{-\beta\lambda}d\lambda &= E\Big(\frac{1}{\beta + \Gamma_t(1/G)}\Big) \\
&= e^{-tE(\log G)} E\Big(\int_0^\infty e^{-\beta\lambda - \frac{1}{\lambda}\Gamma_t(G)}\frac{d\lambda}{\lambda^t}\Big) \\
&= 2e^{-tE(\log G)} E\Big\{\Big(\frac{\Gamma_t(G)}{\beta}\Big)^{\frac{1-t}{2}} K_{1-t}\big(2\sqrt{\beta\Gamma_t(G)}\big)\Big\}
\end{aligned} \tag{117}$$

where $K_\nu$ denotes the Bessel-McDonald function with index $\nu$ and where we have used formula 5.10.25 in [34], p. 108 and 119.

We now use (117) to compute $f_{\Gamma_t(1/G)}$ by inverting its Stieltjes transform [53].

$$\begin{aligned}
f_{\Gamma_t(1/G)}(u) &= \frac{1}{2i\pi} \lim_{\eta\downarrow 0} \Big\{ E\Big(\frac{1}{-u - i\eta + \Gamma_t(1/G)} - \frac{1}{-u + i\eta + \Gamma_t(1/G)}\Big)\Big\} \\
&= \frac{e^{-tE(\log G)}}{i\pi} E\Big(\Big(\frac{\Gamma_t(G)}{u}\Big)^{\frac{1-t}{2}} \Big[e^{\frac{i\pi(1-t)}{2}} K_{1-t}\big(2\sqrt{u\Gamma_t(G)}\,e^{-\frac{i\pi}{2}}\big) \\
&\quad - e^{-\frac{i\pi(1-t)}{2}} K_{1-t}\big(2\sqrt{u\Gamma_t(G)}\,e^{-\frac{i\pi}{2}}\big)\Big]\Big)
\end{aligned} \tag{118}$$

However, it is well known that (see [34], p. 108 and 109):

- $K_\nu = K_{-\nu}$
- $\pi < \text{Arg}\, z < \frac{\pi}{2} \Rightarrow K_\nu(z) = \frac{i\pi}{2} e^{\frac{i\nu\pi}{2}} H_\nu^{(1)}(z\, e^{\frac{i\pi}{2}})$
- $\frac{\pi}{2} < \text{Arg}\, z < \pi \Rightarrow K_\nu(z) = -\frac{i\pi}{2} e^{-\frac{i\nu\pi}{2}} H_\nu^{(2)}(z\, e^{-\frac{i\pi}{2}})$

where the Hankel functions $H_\nu^{(1)}$ and $H_\nu^{(2)}$ satisfy ([34], p. 108)

$$H_\nu^{(1)} + H_\nu^{(2)} = 2J_\nu$$

Plugging these relations into (118), we obtain:

$$f_{\Gamma_t(1/G)}(u) = e^{-tE(\log G)} E\Big\{\Big(\frac{\Gamma_t(G)}{u}\Big)^{\frac{1-t}{2}} J_{t-1}\big(2\sqrt{u\Gamma_t(G)}\big)\Big\}$$

*ii)* We now prove (104), which is a consequence of (102). Indeed, from (102) and (103):

$$\begin{aligned}
f_{\Gamma_t(1/G)}(u) &= e^{-tE(\log G)} E\Big\{\Big(\frac{\Gamma_t(G)}{u}\Big)^{\frac{1-t}{2}} \sum_{k=0}^\infty \frac{(-1)^k}{\Gamma(k+1)\Gamma(k+t)} \big(u\Gamma_t(G)\big)^{k+\frac{t-1}{2}}\Big\} \\
&= e^{-tE(\log G)} u^{t-1} E\Big\{\sum_{k=0}^\infty \frac{(-1)^k}{\Gamma(k+1)\Gamma(k+t)} \big(u\gamma_t D_t(G)\big)^k\Big\}
\end{aligned}$$

from (53). But, since $E(\gamma_t^k) = \frac{\Gamma(k+t)}{\Gamma(t)}$, we have:

$$\begin{aligned}
f_{\Gamma_t(1/G)}(u) &= \frac{e^{-tE(\log G)}}{\Gamma(t)} u^{t-1} E\Big(\sum_{k=0}^{\infty} \frac{(-1)^k}{k!} \big(uD_t(G)\big)^k\Big) \\
&= \frac{e^{-tE(\log G)}}{\Gamma(t)} u^{t-1} E(e^{-uD_t(G)}) \tag{119}
\end{aligned}$$

On the other hand, as we already noticed in point 6 of Remark 2.2, formula (105) may be obtained by integrating from 0 to $\infty$ in (119).

**2.2.c** <u>We now prove point 4 of Theorem 2.1.</u>

*i)* Formula (106) is immediate. Indeed, since $\Gamma_t(G) \stackrel{(law)}{=} \gamma_t D_t(G)$ we get:

$$E(e^{-\lambda \Gamma_t(G)}) = \exp\big(-t\psi_G(\lambda)\big) = E\Big(\frac{1}{\big(1+\lambda D_t(G)\big)^t}\Big)$$

Replacing then $\lambda$ by $1/\lambda$ in this latter formula, we obtain:

$$\begin{aligned}
E\Big[\frac{1}{\big(\lambda + D_t(G)\big)^t}\Big] &= \exp\big\{-t\big(\psi_G(1/\lambda) + \log\lambda\big)\big\} \\
&= \exp\big\{-t\big(\psi_{1/G}(\lambda) - E(\log G)\big)\big\} \tag{120}
\end{aligned}$$

from (101). We note that (105) may also be obtained by taking $\lambda = 0$ in (120).

*ii)* We now show (107). For this purpose, we shall prove that the two members of (107) admit the same Stieltjes transform with index $t$. Indeed:

$$\int_0^{\infty} \frac{e^{-tE(\log G)}}{(\lambda+x)^t} x^{t-2} f_{D_t(G)}\Big(\frac{1}{x}\Big) dx = e^{-tE(\log G)} \int_0^{\infty} \frac{1}{(1+\lambda x)^t} f_{D_t(G)}(x) dx$$

after making the change of variables $y = 1/x$

$$\begin{aligned}
&= e^{-tE(\log G)} E(e^{-\lambda \Gamma_t(G)}) \qquad \text{from (56)} \\
&= \exp -t\big(E(\log G) + \psi_G(\lambda)\big) \tag{121}
\end{aligned}$$

whereas:

$$\begin{aligned}
\int_0^{\infty} \frac{1}{(\lambda+x)^t} f_{D_t(1/G)}(x) dx &= \frac{1}{\lambda^t} \int_0^{\infty} \frac{1}{\big(1+\frac{x}{\lambda}\big)^t} f_{D_t(1/G)}(x) dx \\
&= \exp\big\{-t\big(\log\lambda + \psi_{1/G}(1/\lambda)\big)\big\} \quad \text{(from (56))} \\
&= \exp\big\{-t\big(E(\log G) + \psi_G(\lambda)\big)\big\} \quad \text{(from (101))}
\end{aligned} \tag{122}$$

The comparison of (121) and (122) and the injectivity of the Stieltjes transform of index $t$ imply (107).

*iii)* We now show (107), for $t = 1$, by using (110). For $t = 1$, (110) writes:

$$f_{D_1(G)}(x) = \frac{\sin(\pi F_{1/G}(x))}{\pi} e^{-E((\log|x - \frac{1}{G}|)1_{x \neq 1/G})} \tag{123}$$

Thus:

$$\begin{aligned}
f_{D_1(1/G)}(x) &= \frac{\sin(\pi F_G(x))}{\pi} e^{-E((\log|x - G|)1_{x \neq G})} \\
&= \frac{\sin(\pi F_G(x))}{\pi} \exp\left\{-E\left(\log\left(xG\left|\frac{1}{x} - \frac{1}{G}\right|\right) \cdot 1_{x \neq G}\right)\right\} \\
&= \frac{\sin(\pi F_G(x))}{\pi} \frac{e^{-E(\log G)}}{x} e^{-E(\log|\frac{1}{x} - \frac{1}{G}|1_{x \neq G})} \\
&= \frac{1}{x} e^{-E(\log G)} f_{D_1(G)}(1/x), \quad \text{from (123)}
\end{aligned}$$

since $\sin(\pi F_G(x)) = \sin(\pi(1 - F_G(x))) = \sin(\pi F_{1/G}(1/x))$. It is possible, by using (110), to extend this proof for all $t < 1$. We leave the details to the interested reader.

### 2.3. A complement to the duality theorem

Here again, $G$ denotes a r.v. such that $E(|\log G|) < \infty$, and we recall that, for any $t \leq 1$ $Y_t$ denotes a Bernoulli r.v. with parameter $t$ (see (41) and (42), with $m = t$ and $m' = 1$).

**Theorem 2.3.** *For all $t \in [0, 1[$:*
**1)** *The density $f_{D_1(\frac{G}{Y_t})}$ may be expressed in terms of $G$, as:*

$$\begin{aligned}
f_{D_1\left(\frac{G}{Y_t}\right)}(x) &= \frac{\sin(\pi t F_G(1/x))}{\pi} x^{t-1} \\
&\quad \times \exp\left\{-t E\left(\left(\log\left|x - \frac{1}{G}\right|\right) 1_{x \neq \frac{1}{G}}\right)\right\} \qquad x > 0, \quad (124)
\end{aligned}$$

*and when $1/x$ is in the support of $F_G$. Otherwise replace $\sin(\pi t F_G(1/x))$ by $\sin(\pi(1 - t))$ when $x > 0$.*
**2)** *The following duality formula holds:*

$$f_{D_1\left(\frac{G}{Y_t}\right)}(x) = \frac{\sin(\pi t F_G(1/x))}{\sin(\pi t(1 - F_G(1/x)))} \cdot x^{t-2} f_{D_1\left(\frac{1}{GY_t}\right)}\left(\frac{1}{x}\right) \cdot e^{t E(\log G)} \tag{125}$$

**3)** *Let $\wedge_t : [0, \infty[ \longrightarrow [0, 1]$ be defined by:*

$$\wedge_t(y) := 1 - \frac{1}{\pi t} \operatorname{arc tg}\left(\frac{\sin(\pi t)}{\cos(\pi t) + y}\right) \tag{126}$$

*Then:*

$$F_G(1/x) = \wedge_t \left( \frac{f_{D_1\left(\frac{G}{Y_t}\right)}(x)}{x^{t-2} f_{D_1\left(\frac{1}{GY_t}\right)}\left(\frac{1}{x}\right) e^{t\,E(\log G)}} \right) \quad (127)$$

**4)** Equivalently,

$$F_G(1/x) = \wedge_t \left( \frac{E[(D_t(G) - x)_+^{-t}]}{E[(x - D_t(G))_+^{-t}]} \right)$$

**Remark 2.4.**

**1)** We note that the right-hand side of (127) depends on $t$ $(t \in [0,1[)$, whereas the left-hand side does not depend on $t$.

**2)** Since (see (43)) $D_1\left(\frac{1}{GY_t}\right) \stackrel{(law)}{=} \beta_{t,1-t} \cdot D_t(1/G)$, the knowledge of the law of $D_t(1/G)$ and of that of $D_t(G)$, for one $t < 1$, allow to determine that of $G$. We shall exploit this fact, in points 2.5 and 2.6 below, to determine the Thorin measure of a Pareto distribution and of a power of a gamma variable.

**3)** We note that finding an explicit Thorin measure of an arbitrary GGC is akin to finding the Lévy measure of some infinitely divisible random variable. Bondesson([5], Theorem 4.3.2, p. 61), using inversion techniques, obtains an expression for the Thorin measure, but notes that it seldom yields explicit expressions. On the other hand the use of statement 3) of Theorem 2.3 will often lead to tractable expressions for the Thorin measure.

**4)** We shall prove (see Section 3.1.b) that $\wedge_t$ is the cumulative distribution function of a r.v. $Z_t$ which we shall describe. On the other hand, some trigonometric computations allow to see that $\wedge_t^{-1}$, the inverse of $\wedge_t$ in the sense of composition of functions, equals:

$$\wedge_t^{-1}(x) = \frac{\sin(\pi t x)}{\sin(\pi t(1-x))} \quad x \in [0,1] \quad (128)$$

**5)** Point 1 of Theorem 2.3 is due to James (see[28]).

### 2.4. Proof of Theorem 2.3

**2.4.a** <u>We first prove point 1.</u>
As this point has already been established by James (see [28]), we shall only give a quick proof.

We have, from (56):

$$E\left[\frac{1}{1+\lambda D_1\left(\frac{G}{Y_t}\right)}\right] = \exp\left(-t\psi_G(\lambda)\right) = \frac{1}{\lambda} E\left(\frac{1}{\frac{1}{\lambda} + D_1\left(\frac{G}{Y_t}\right)}\right)$$

Hence, changing $\lambda$ in $\frac{1}{\lambda}$ and using (101):

$$\begin{aligned}
E\left[\frac{1}{\lambda + D_1\left(\frac{G}{Y_t}\right)}\right] &= \exp\{-t\,\psi_G(1/\lambda) - \log\lambda\} \\
&= \exp\left((t-1)\log\lambda\right) \cdot e^{t\,E(\log)} \exp\left(-t\,\psi_{1/G}(\lambda)\right) \quad (129)
\end{aligned}$$

We then compute $f_{D_1\left(\frac{G}{Y_t}\right)}$ by inverting its Stieltjes transform:

$$f_{D_1\left(\frac{G}{Y_t}\right)}(u) =$$
$$\frac{e^{t\,E(\log G)}}{2i\pi} \lim_{\eta\downarrow 0} \left\{\exp\left[((t-1)\log(-u-i\eta)) - t\,\psi_{1/G}(-u-i\eta)\right]\right.$$
$$\left. -\exp\left[((t-1)\log(-u+i\eta)) - t\,\psi_{1/G}(-u+i\eta)\right]\right\} \qquad (u>0) \quad (130)$$

It then suffices to observe that:

$$\exp(t-1)\log(-u-i\eta) \xrightarrow[\eta\downarrow 0_+]{} \exp\{(t-1)\log|u| - i\pi(t-1)\} = u^{t-1}e^{-i\pi(t-1)}$$

and:

$$\lim_{\eta\downarrow 0_+} \exp\{(t-1)\log(-u+i\eta)\} = u^{t-1}e^{i\pi(t-1)}$$

as well as:

$$t\,\psi_{1/G}(-u-i\eta) = t\,E\big(\log(1-uG-i\eta G)\big)$$
$$= t\,E\big(\log(1-uG-i\eta G)1_{uG<1} + \log(1-uG-i\eta G)1_{uG>1}\big)$$
$$\xrightarrow[\eta\downarrow 0_+]{} t\,E\big(\log(|1-uG|)1_{uG\neq 1} - i\pi P(uG>1)\big)$$

whereas:

$$t\,\psi_{1/G}(-u+i\eta) \xrightarrow[\eta\downarrow 0_+]{} t\,E\big(\log(|1-uG|)1_{uG\neq 1} + i\pi P(uG>1)\big)$$

Then, plugging the values of these different limits in (130), we obtain point 1 of Theorem 2.3.

**2.4.b** We now prove point 2 of Theorem 2.3.
We deduce from (104) that:

$$f_{D_1\left(\frac{G}{Y_t}\right)}(x)$$
$$= \frac{\sin(\pi t\,F_G(1/x))}{\pi} x^{t-1}\exp\left\{-t\,E\left(\log\left(\frac{1}{G}|1-xG|\right)1_{xG\neq 1}\right)\right\}$$
$$= \frac{\sin(\pi t\,F_G(1/x))}{\pi} x^{t-1}e^{t\,E(\log G)}\exp\left\{-t\,E\big((\log|1-xG|)1_{xG\neq 1}\big)\right\} \quad (131)$$

whereas:

$$x^{t-2}f_{D_1\left(\frac{1}{GY_t}\right)}\left(\frac{1}{x}\right)e^{t\,E(\log G)} = x^{t-2}\frac{\sin(\pi t F_{1/G}(x))}{\pi} e^{t\,E(\log G)} \cdot x^{1-t}$$
$$\cdot \exp\left\{-t\,E\left(\log\left(\left|\frac{1}{x}-\frac{1}{G}\right|\right)1_{xG\neq 1}\right)\right\}$$
$$= \frac{1}{x}\frac{\sin(\pi t F_{1/G}(x))}{\pi} e^{tE(\log G)}\exp\left\{-t\,E\left(\log\left(\frac{1}{x}|1-xG|\right)1_{xG\neq 1}\right)\right\}$$
$$= \frac{\sin(\pi t\,F_{1/G}(x))}{\pi} x^{t-1}e^{tE(\log G)}\exp\left\{-t\,E\big(\log|1-xG|\big)1_{xG\neq 1}\right\} \quad (132)$$

hence, point 2 of Theorem 2.3 follows, by comparison of (132) and (131).

**2.4.c** We now prove point 3 of Theorem 2.3.
From (125), we obtain:

$$\frac{f_{D_1\left(\frac{G}{Y_t}\right)}(x)}{x^{t-2}f_{D_1\left(\frac{1}{GY_t}\right)}\left(\frac{1}{x}\right)e^{tE(\log G)}} = \frac{\sin\left(\pi t\, F_G(1/x)\right)}{\sin\left(\pi t(1-F_G(1/x))\right)} = \wedge_t^{-1}\left(F_G(1/x)\right)$$

$$\text{with} \qquad \wedge_t^{-1}(x) = \frac{\sin(\pi t x)}{\sin \pi t(1-x)} \qquad \left(\text{see }(128)\right)$$

Then, inverting this formula, since $\wedge_t \circ \wedge_t^{-1}(x) = x, \quad x \in [0,1]$ we obtain:

$$F_G(1/x) = \wedge_t\left(\frac{f_{D_1\left(\frac{G}{Y_t}\right)}(x)}{x^{t-2}e^{tE(\log G)}f_{D_1\left(\frac{1}{GY_t}\right)}\left(\frac{1}{x}\right)}\right)$$

**2.4.d** We now prove point 4 of Theorem 2.3. This result follows by using the generic form of the density of a $\beta_{t,1-t}$ random variable multiplied by an independent random variable in the particular cases of $\beta_{t,1-t}D_t(G)$ and $\beta_{t,1-t}D_t(1/G)$. The result is completed by applying the identity in (107).

## 2.5. Computation of the Thorin measure of a Pareto r.v.

Here is a first application of Theorem 2.3.

**2.5.a** Let $m > 0$ fixed and:

$$X_\theta := \frac{\gamma_\theta}{\gamma_m} \qquad (0 < \theta < 1) \tag{133}$$

with density:

$$f_{X_\theta}(x) = \frac{\Gamma(\theta+m)}{\Gamma(\theta)\,\Gamma(m)}\, x^{\theta-1}(1+x)^{-(\theta+m)} \tag{134}$$

$$= \frac{\Gamma(\theta+m)}{\Gamma(\theta)\cdot\Gamma(m)}\, x^{\theta-1}E(e^{-x\gamma_{\theta+m}}) \tag{135}$$

The r.v. $X_\theta$ is a GGC r.v. (see Bondesson [5], p. 59 ). The rationale of our work is now the following:
• We first compute the Thorin measure associated with $X_\theta$.
• Then, letting $\theta$ converge 1, we shall obtain - as the Thorin measure depends continuously (for the narrow topology) on the law of $X_\theta$ (see Bondesson, [5]) - the Thorin measure associated with the r.v. $X_1 = \frac{\gamma_1}{\gamma_m}$, i.e. to a Pareto r.v. with parameter $m(m > 0)$.

**2.5.b** Thorin measure associated with $X_\theta$, $\theta < 1$.
Since, from (134):

$$\sup\left\{\alpha > 0 \ ; \ \lim_{x \to 0} \frac{f_{X_\theta}(x)}{x^{\alpha-1}} = 0\right\} = \theta$$

we deduce from (30) the existence of a r.v. $G_\theta$, such that $E\left(\log^+\left(\frac{1}{G_\theta}\right)\right) < \infty$ and such that $X_\theta$ is a $(\theta, G_\theta)$ GGC r.v.:

$$X_\theta \stackrel{(law)}{=} \Gamma_\theta(G_\theta) \stackrel{(law)}{=} \gamma_\theta D_\theta(G_\theta) \stackrel{(law)}{=} \frac{\gamma_\theta}{\gamma_m} \tag{136}$$

hence:

$$D_\theta(G_\theta) \stackrel{(law)}{=} \frac{1}{\gamma_m} \tag{137}$$

On the other hand, from (104):

$$f_{\Gamma_\theta(G_\theta)}(x) = \frac{x^{\theta-1}}{\Gamma(\theta)} e^{\theta E(\log G_\theta)} E(e^{-x D_\theta(1/G_\theta)}) \tag{138}$$

Comparing (137) and (135) yields:

$$e^{\theta E(\log G_\theta)} = \frac{\Gamma(\theta+m)}{\Gamma(m)} \tag{139}$$

$$D_\theta(1/G_\theta) \stackrel{(law)}{=} \gamma_{\theta+m} \tag{140}$$

Since, on the other hand:

$$D_1\left(\frac{1}{G_\theta Y_\theta}\right) \stackrel{(law)}{=} \beta_{\theta,1-\theta} \cdot D_\theta(1/G_\theta) \stackrel{(law)}{=} \beta_{\theta,1-\theta} \cdot \gamma_{\theta+m} \tag{141}$$

$$D_1\left(\frac{G_\theta}{Y_\theta}\right) \stackrel{(law)}{=} \beta_{\theta,1-\theta} \cdot D_\theta(G_\theta) \stackrel{(law)}{=} \beta_{\theta,1-\theta} \cdot \frac{1}{\gamma_m} \tag{142}$$

we easily deduce from these formulae that:

$$f_{D_1\left(\frac{1}{G_\theta Y_\theta}\right)}(z) = \frac{\sin \pi\theta}{\pi \Gamma(\theta+m)} z^{\theta-1} \int_z^\infty e^{-y} \frac{y^{m+\theta-1}}{(y-z)^\theta} dy \tag{143}$$

$$f_{D_1\left(\frac{G_\theta}{Y_\theta}\right)}(z) = \frac{\sin \pi\theta}{\pi \Gamma(\theta)} z^{-m-1} \int_0^1 e^{-\frac{y}{z}} y^{m+\theta-1}(1-y)^{-\theta} dy \tag{144}$$

The Thorin measure of $X_\theta$, which equals: $\theta P_{G_\theta}(dx)$, where $P_{G_\theta}$ is the law of $G_\theta$, is then, by applying (127):

$$F_{G_\theta}\left(\frac{1}{z}\right) = \wedge_\theta \left(\frac{\int_z^\infty e^{-y} \frac{y^{m+\theta-1}}{(y-z)^\theta} dy}{z^m \int_0^1 e^{-yz} y^{m+\theta-1}(1-y)^{-\theta} dy}\right) \tag{145}$$

where, to obtain (145), we have used (139).

**2.5.c** Thorin measure of $\frac{\gamma_1}{\gamma_m}$ $(m > 0)$.

By continuity, the r.v. $\frac{\gamma_1}{\gamma_m}$ is a $(1, G)$ GGC r.v. and its Thorin measure is the law of $G$ whose cumulative distribution function is obtained by letting $\theta$ tend to 1 in (145). To obtain this limit, we shall develop several computations.

*2.5.i)* Development of $\int_z^\infty e^{-y} \frac{y^{m+\theta-1}}{(y-z)^\theta} dy$ as $\beta = 1 - \theta \to 0$

$$\int_z^\infty e^{-y} \frac{y^{m+\theta-1}}{(y-z)^\theta} dy$$

$$= e^{-z} \int_0^\infty e^{-u} \frac{(z+u)^{m+\theta-1}}{u^\theta} du \qquad (y = z+u)$$

$$= \frac{e^{-z}}{\beta} \int_0^\infty e^{-u} u^\beta \big[(z+u)^{m-\beta} - (m-\beta)(z+u)^{m-\beta-1}\big] du$$

with $(\beta = 1 - \theta)$ and after integrating by parts:

$$= \frac{e^{-z}}{\beta} \int_0^\infty e^{-u}(z+u)^m \left[\left(\frac{u}{z+u}\right)^\beta \left(1 - (m-\beta)\frac{1}{z+u}\right)\right] du$$

$$\underset{\beta \to 0}{\sim} \frac{e^{-z}}{\beta} \int_0^\infty e^{-u}(z+u)^m \left[\left(1 + \beta \log\left(\frac{u}{z+u}\right)\right)\right]\left(1 - (m-\beta)\frac{1}{z+u}\right) du$$

$$\underset{\beta \to 0}{=} \frac{e^{-z}}{\beta} \int_0^\infty e^{-u}(z+u)^m \left(1 - \frac{m}{z+u}\right) du$$

$$+ e^{-z} \int_0^\infty e^{-u}(z+u)^m \left(\log \frac{u}{1+u} + \frac{1}{z+u}\right) du + o(\beta)$$

$$= \frac{e^{-z}z^m}{\beta} + e^{-z}z^m \int_0^\infty e^{-u}\left(1 + \frac{u}{z}\right)^m \left(\log \frac{u}{1+u} + \frac{1}{z+u}\right) du + o(\beta)$$

since $\int_0^\infty e^{-u}(z+u)^m \left(1 - \frac{m}{z+u}\right) du = z^m$, by integration by parts. Finally:

$$\int_z^\infty e^{-y} \frac{y^{m+\theta-1}}{(y-z)^\theta} dy \underset{\beta \to 0}{=} \frac{e^{-z}z^m}{\beta} + e^{-z}z^m C_m(z) + o(\beta) \qquad (146)$$

with

$$C_m(z) = \int_0^\infty e^{-u}\left(1 + \frac{u}{z}\right)^m \left(\log \frac{u}{1+u} + \frac{1}{z+u}\right) du \qquad (147)$$

*2.5. ii)* Development of $z^m \int_0^1 e^{-xz} x^{m+\theta-1}(1-x)^{-\theta} dx$ as $\beta = 1 - \theta \longrightarrow 0$

$$z^m \int_0^1 e^{-xz} x^{m-\beta}(1-x)^{\beta-1} dx$$

$$= e^{-z} z^m \int_0^1 e^{zu}(1-u)^{m-\beta} u^\beta du \qquad (x = 1-u)$$

$$= \frac{e^{-z}z^m}{\beta} \int_0^1 e^{zu}\big((m-\beta)(1-u)^{m-\beta-1} - z(1-u)^{m-\beta}\big) u^\beta du$$

with $\beta = 1 - \theta$, and after integrating by parts:

$$= \frac{e^{-z}z^m}{\beta} \int_0^1 e^{zu}(1-u)^m \left(\frac{u}{1-u}\right)^\beta \left(\frac{m-\beta}{1-u} - z\right) du$$

$$\underset{\beta \to 0}{\sim} \frac{e^{-z}z^m}{\beta} \int_0^1 e^{zu}(1-u)^m \left(\frac{m-\beta}{1-u} - z\right)\left(1 + \beta \log \frac{u}{1-u}\right) du$$

$$\underset{\beta \to 0}{\sim} \frac{e^{-z}z^m}{\beta} \int_0^1 e^{zu}(1-u)^m \left(\frac{m}{1-u} - z\right) du + e^{-z}z^m \int_0^1 e^{zu}(1-u)^m \cdot$$

$$\left[-\frac{1}{1-u} + \left(\frac{m}{1-u} - z\right)\log \frac{u}{1-u}\right] du + o(\beta)$$

$$= \frac{e^{-z}z^m}{\beta} + e^{-z}z^m \widetilde{C}_m(z) + o(\beta) \qquad (148)$$

with

$$\widetilde{C}_m(z) = \int_0^1 e^{zu}(1-u)^m \left(-\frac{1}{1-u} + \left(\frac{m}{1-u} - z\right)\log \frac{u}{1-u}\right) du \qquad (149)$$

since: $\int_0^1 e^{zu}(1-u)^m \left(\frac{m}{1-u} - z\right) du = 1$, by integrating by parts.

2.5. *iii)* Let

$$q_{1-\beta}(z) := \frac{\int_z^\infty e^{-y} \frac{y^{m+\theta-1}}{(y-z)^\theta} dy}{z^m \int_0^1 e^{-xz} x^{m+\theta-1}(1-x)^\theta dx} \qquad (150)$$

Thus, we have:

$$q_{1-\beta}(z) \underset{\beta \to 0}{=} \frac{\dfrac{e^{-z}}{\beta}z^m + e^{-z}z^m C_m(z) + o(\beta)}{\dfrac{e^{-z}z^m}{\beta} + e^{-z}z^m \widetilde{C}_m(z) + o(\beta)} \qquad (151)$$

where $C_m(z)$ and $\widetilde{C}_m(z)$ are given by (147) and (149). Hence:

$$q_{1-\beta}(z) = 1 + \beta\big(C_m(z) - \widetilde{C}_m(z)\big) + o(\beta) \qquad (152)$$

Plugging this expression in (145), we obtain:

$$F_{G_{1-\beta}}(1/z)$$

$$= 1 - \frac{1}{\pi\theta} \operatorname{arc tg}\left(\frac{\sin(\pi\theta)}{\cos(\pi\theta) + 1 + \beta\big(C_m(z) - \widetilde{C}_m(z) + o(\beta)\big)}\right)$$

$$= 1 - \frac{1}{\pi(1-\beta)} \operatorname{arc tg}\left(\frac{\sin \pi\beta}{-\cos \pi\beta + 1 + \beta\big(C_m(z) - \widetilde{C}_m(z) + o(\beta)\big)}\right)$$

$$\underset{\beta \to 0}{\longrightarrow} 1 - \frac{1}{\pi} \operatorname{arc tg}\left(\frac{\pi}{C_m(z) - \widetilde{C}_m(z)}\right) = F(1/z) \qquad (153)$$

where $F$ is the cumulative distribution function of the Thorin measure associated to the Pareto r.v. $\frac{\gamma_1}{\gamma_m}$, with parameter $m > 0$.

## 2.6. Thorin measure associated with $\gamma_1^{1/\alpha}$, $0 < \alpha < 1$

**2.6.a** Let $X_\alpha \stackrel{(law)}{=} \gamma_1^{1/\alpha}$ with $0 < \alpha < 1$, with density:

$$f_{X_\alpha}(x) = \alpha\, x^{\alpha-1} e^{-x^\alpha} = \alpha\, x^{\alpha-1} E(e^{-xS_\alpha}) \tag{154}$$

where $S_\alpha$ is a positive stable r.v. with index $\alpha$:

$$E(e^{-\lambda S_\alpha}) = \exp(-\lambda^\alpha) \qquad (\lambda \geq 0) \tag{155}$$

From Bondesson ([5], p. 60), $X_\alpha$ is a GGC r.v. Since:

$$\alpha = \sup\left\{\nu > 0,\ \lim_{x \downarrow 0} \frac{f_{X_\alpha}(x)}{x^{\nu-1}} = 0\right\}$$

$X_\alpha$ is a $(\alpha, G_\alpha)$ GGC r.v., for a r.v. $G_\alpha$ such that: $E(\log^+(1/G_\alpha)) < \infty$. Thus:

$$X_\alpha \stackrel{(law)}{=} (\gamma_1)^{1/\alpha} = \Gamma_\alpha(G_\alpha) = \gamma_\alpha D_\alpha(G_\alpha) \tag{156}$$

We shall now devote some effort to finding the law of $G_\alpha$, making use in particular of formula (107)

From (104):

$$f_{\Gamma_\alpha(G_\alpha)}(x) = \frac{x^{\alpha-1}}{\Gamma(\alpha)}\, e^{\alpha E(\log G_\alpha)} E(e^{-xD_\alpha(1/G_\alpha)}) \tag{157}$$

we deduce, by comparison with (154):

$$e^{\alpha E(\log G_\alpha)} = \alpha\, \Gamma(\alpha) = \Gamma(\alpha+1) \tag{158}$$

and

$$D_\alpha(1/G_\alpha) \stackrel{(law)}{=} S_\alpha \tag{159}$$

Since:

$$D_1\!\left(\frac{1}{G_\alpha Y_\alpha}\right) \stackrel{(law)}{=} \beta_{\alpha,1-\alpha} \cdot D_\alpha(1/G_\alpha) \stackrel{(law)}{=} \beta_{\alpha,1-\alpha} \cdot S_\alpha$$

we deduce that the density of $D_1\!\left(\frac{1}{G_\alpha Y_\alpha}\right)$ equals:

$$f_{D_1(\frac{1}{G_\alpha Y_\alpha})}(y) = \frac{\sin(\pi\alpha)}{\pi}\, y^{\alpha-1} \int_y^\infty \frac{1}{(x-y)^\alpha}\, f_\alpha(x)\, dx \tag{160}$$

where $f_\alpha$ denotes the density of $S_\alpha$. On the other hand, from Chaumont-Yor ([9], p. 112):

$$\gamma_1^{1/\alpha} \stackrel{(law)}{=} \frac{\gamma_1}{S_\alpha} \tag{161}$$

hence, since:

$$\frac{\gamma_1}{S_\alpha} \stackrel{(law)}{=} \gamma_\alpha \cdot D_\alpha(G_\alpha) \stackrel{(law)}{=} \beta_{\alpha,1-\alpha} \cdot \gamma_1 D_\alpha(G_\alpha) \stackrel{(law)}{=} \gamma_1 \cdot D_1\!\left(\frac{G_\alpha}{Y_\alpha}\right) \tag{162}$$

we get:

$$D_1\left(\frac{G_\alpha}{Y_\alpha}\right) \stackrel{(law)}{=} \frac{1}{S_\alpha} \quad \text{and} \tag{163}$$

$$f_{D_1(\frac{G_\alpha}{Y_\alpha})}(x) = \frac{1}{x^2} f_\alpha\left(\frac{1}{x}\right) \tag{164}$$

Then, applying (127), we get:

$$F_{1/G_\alpha}\left(\frac{1}{y}\right) = \wedge_\alpha \left(\frac{\frac{\Gamma(\alpha+1)\sin(\pi\alpha)}{\pi} \int_y^\infty \frac{1}{(x-y)^\alpha} f_\alpha(x)dx}{yf_\alpha(y)}\right) \tag{165}$$

where, to obtain (165), we have used (157). This expression (165), provides us then with the explicit form of the Thorin measure of $(\gamma_1)^{1/\alpha}$, which is equal to $P_{G_\alpha}$.

**2.6.b.** We shall now give a more suitable expression of (165). From the relations (162) and (163), we deduce:

$$\beta_{\alpha,1-\alpha} \cdot D_\alpha(G_\alpha) \stackrel{(law)}{=} D_1\left(\frac{G_\alpha}{Y_\alpha}\right) = \frac{1}{S_\alpha} \tag{166}$$

hence:

$$\frac{1}{y^2} f_\alpha\left(\frac{1}{y}\right) = \frac{\sin \pi\alpha}{\pi} y^{\alpha-1} \int_y^\infty \frac{1}{(z-y)^\alpha} f_{D_\alpha(G_\alpha)}(z)dz \tag{167}$$

But, the density of $D_\alpha(G_\alpha)$ may be computed from that of $D_\alpha(1/G_\alpha)$, thanks to (107):

$$f_{D_\alpha(G_\alpha)}(x) = x^{\alpha-2} e^{\alpha E(\log G_\alpha)} f_{D_\alpha(G_\alpha)}(1/x).$$

Now from (159) and (158), we obtain:

$$f_{D_\alpha(G_\alpha)}(x) = x^{\alpha-2}\Gamma(\alpha+1)f_\alpha(1/x). \tag{168}$$

We now note the interesting relationship concerning the law of $S_\alpha$.

**Lemma 2.4**
*Let $0 < \alpha < 1$, and $S$ denote a positive random variable with density $(f(y), y > 0)$ such that:*

$$yf(y) = CE\left[\frac{1}{(y-S)^\alpha} 1_{\{S<y\}}\right] \tag{169}$$

*for some $C > 0$. Then, $S$ is a stable $(\alpha)$ variable; precisely:*

$$E\left(\exp\{-\lambda S\}\right) = \exp\left\{-\frac{C\Gamma(1-\alpha)}{\alpha}\lambda^\alpha\right\}.$$

We postpone the proof of the Lemma for the moment, and we note that, plugging (169) into (165), we obtain:

$$F_{1/G_\alpha}\left(\frac{1}{y}\right)$$
$$= \wedge_\alpha \left( \frac{E\left(\frac{1}{(S_\alpha-y)^\alpha} 1_{S_\alpha>y}\right)}{E\left(\frac{1}{(y-S_\alpha)^\alpha} 1_{S_\alpha<y}\right)} \right)$$
$$= 1 - \frac{1}{\pi\alpha} \operatorname{arc tg} \left[ \frac{\sin(\pi\alpha) \cdot E\left(\frac{1}{(y-S_\alpha)^\alpha} 1_{S_\alpha<y}\right)}{\cos(\pi\alpha) \cdot E\left(\frac{1}{(y-S_\alpha)^\alpha} 1_{S_\alpha<y}\right) + E\left(\frac{1}{(S_\alpha-y)^\alpha} 1_{S_\alpha>y}\right)} \right], \quad (170)$$

from (126).

**Proof of Lemma 2.4** From (169), we take the Laplace transform of both sides:

$$\int_0^\infty e^{-\lambda y} y f(y) dy = CE\left[\int_S^\infty \frac{e^{-\lambda y}}{(y-S)^\alpha} dy\right]$$
$$= CE\left[e^{-\lambda S} \int_0^\infty \frac{e^{-\lambda z}}{z^\alpha} dy\right]$$
$$= C\Gamma(1-\alpha)\lambda^{\alpha-1} E\left(e^{-\lambda S}\right).$$

Denoting $\phi(\lambda) = E\left(e^{-\lambda S}\right)$, we get:

$$-\phi'(\lambda) = C\Gamma(1-\alpha)\lambda^{\alpha-1}\phi(\lambda),$$

from which we deduce:

$$\phi(\lambda) = \exp\left\{-\frac{C\Gamma(1-\alpha)}{\alpha}\lambda^\alpha\right\}.$$

## 3. Explicit examples of GGC variables associated with the ($\mathbb{G}_\alpha$, $0 \leq \alpha \leq 1$) family

All the examples discussed in this Section are related to the r.v.'s ($\mathbb{G}_\alpha$, $0 \leq \alpha \leq 1$) introduced in [4]. Below, we indicate the properties of these r.v.'s which we shall use. We also recall our notation:

$$\Gamma_t(\mathbb{G}_\alpha) = \gamma_t \cdot D_t(\mathbb{G}_\alpha) \qquad \text{and} \qquad (171)$$

$$E(e^{-\lambda \Gamma_t(\mathbb{G}_\alpha)}) = \exp\left\{-t \int_0^\infty (1-e^{-\lambda x})\frac{dx}{x} E(e^{-x\mathbb{G}_\alpha})\right\} \qquad (\lambda, t \geq 0) \quad (172)$$

### 3.1. The family ($\mathbb{G}_\alpha$, $0 \leq \alpha \leq 1$)

(see [4])

**3.1.a.** For $0 < \alpha < 1$, the density $f_{\mathbb{G}_\alpha}$ of $\mathbb{G}_\alpha$ equals:

$$f_{\mathbb{G}_\alpha}(x) = \frac{\alpha \sin(\pi\alpha)}{(1-\alpha)\pi} \frac{x^{\alpha-1}(1-x)^{\alpha-1}}{(1-x)^{2\alpha} - 2(1-x)^\alpha x^\alpha \cos(\pi\alpha) + x^{2\alpha}} 1_{[0,1]}(x) \quad (173)$$

In particular:

- for $\alpha = 1/2$, $\mathbb{G}_{1/2}$ follows the arc sine law:

$$f_{\mathbb{G}_{1/2}}(x) = \frac{1}{\pi} \frac{1}{\sqrt{x(1-x)}} 1_{[0,1]}(x), \quad \mathbb{G}_{1/2} \stackrel{(law)}{=} \beta_{1/2,\,1/2} \tag{174}$$

- for $\alpha = 1$, $\mathbb{G}_1$ is uniform on $[0,1]$ (175)

- for $\alpha = 0$, $\mathbb{G}_0 \stackrel{(law)}{=} \dfrac{1}{1+\exp(\pi C)}$ (176)

where $C$ is a standard Cauchy variable.

In general, for $0 < \alpha < 1$ one has:

$$\begin{aligned} E[e^{-\lambda \Gamma_{1-\alpha}(\mathbb{G}_\alpha)}] &= \exp\Big\{-(1-\alpha)\int_0^\infty (1-e^{-\lambda x})\frac{dx}{x} E(e^{-x\mathbb{G}_\alpha})\Big\} \\ &= (1+\lambda)^\alpha - \lambda^\alpha \end{aligned} \tag{177}$$

The density $f_{\Gamma_{1-\alpha}(\mathbb{G}_\alpha)}$ of $\Gamma_{1-\alpha}(\mathbb{G}_\alpha)$ equals:

$$f_{\Gamma_{1-\alpha}(\mathbb{G}_\alpha)}(x) = \frac{\alpha}{\Gamma(1-\alpha)} \frac{1}{x^{1+\alpha}} (1-e^{-x}) 1_{[0,\infty]}(x) \tag{178}$$

which may be translated as the following identities in law:

$$\Gamma_{1-\alpha}(\mathbb{G}_\alpha) \stackrel{(law)}{=} \frac{\gamma_{1-\alpha}}{\beta_{\alpha,1}} \stackrel{(law)}{=} \frac{\gamma_{1-\alpha}}{U^{1/\alpha}} \tag{179}$$

**3.1.b.** We note, for $0 < \mu < 1$, $S_\mu$ and $S'_\mu$ two independent copies of positive stable $(\mu)$ r.v.'s, i.e.:

$$E(e^{-\lambda S_\mu}) = \exp(-\lambda^\mu) \qquad (\lambda \geq 0) \tag{180}$$

and, we let:

$$Z_\mu := \Big(\frac{S_\mu}{S'_\mu}\Big)^\mu \tag{181}$$

Then (see [33] or [8], p. 116), the density $f_{Z_\mu}$ of $Z_\mu$ equals:

$$f_{Z_\mu}(x) = \frac{\sin(\pi\mu)}{\pi\mu} \frac{1}{x^2 + 2x\cos(\pi\mu) + 1} 1_{[0,\infty]}(x) \tag{182}$$

and we have:

$$\mathbb{G}_\alpha \stackrel{(law)}{=} \frac{(Z_{1-\alpha})^{1/\alpha}}{1+(Z_{1-\alpha})^{1/\alpha}} \tag{183}$$

or equivalently:

$$\frac{1}{\mathbb{G}_\alpha} \stackrel{(law)}{=} 1 + \frac{1}{(Z_{1-\alpha})^{1/\alpha}} \stackrel{(law)}{=} 1 + \frac{S'_{1-\alpha}}{S_{1-\alpha}} \tag{184}$$

We note that the cumulative distribution function $F_{Z_\mu}$ of $Z_\mu$ equals:

$$F_{Z_\mu}(x) = 1 - \frac{1}{\pi\mu} \operatorname{arc tg} \left[\frac{\sin(\pi\mu)}{\cos(\pi\mu) + x}\right] \qquad x \geq 0 \qquad (185)$$

and that its inverse, in the sense of composition of functions, is given by:

$$F_{Z_\mu}^{-1}(x) = \frac{\sin(\pi\mu x)}{\sin(\pi\mu(1-x))} \qquad (0 \leq x \leq 1) \qquad (186)$$

(see Remark 2.4, point 3).

**3.1.c.** Although this will not be used in the sequel, we indicate a realization of the r.v. $\Gamma_{1-\alpha}(\mathbb{G}_\alpha)$ which has been at the start of [4]. Let $(R_t, t \geq 0)$ denote a Bessel process starting from 0, with dimension $d = 2(1-\alpha)$, with $0 < d < 2$, or equivalently $0 < \alpha < 1$. Let, for any $t > 0$:

$$g_t^{(\alpha)} := \sup\{s \leq t \,;\, R_s = 0\}, \quad d_t^{(\alpha)} := \inf\{s \geq t, \, R_s = 0\} \qquad (187)$$

and let $\mathfrak{e} \overset{(law)}{=} \gamma_1$ an exponentially distributed r.v., with mean 1, independent from $(R_t, t \geq 0)$. Then:

$$\Gamma_{1-\alpha}(\mathbb{G}_\alpha) \overset{(law)}{=} d_\mathfrak{e}^{(\alpha)} - g_\mathfrak{e}^{(\alpha)} \qquad (188)$$

A more general study of quantities such as the RHS of (188), has been developed by M. Winkel [54].

### 3.2. Study of the subordinators $\bigl(\Gamma_t(\mathbb{G}_{1/2}), t \geq 0\bigr)$ and $\bigl(\Gamma_t(1/\mathbb{G}_{1/2}), t \geq 0\bigr)$

In this Section, $\mathbb{G}_{1/2}$ is a beta $(1/2, 1/2)$ variable, i.e.: its law is the arc sine distribution. (see (174)).

**3.2.a. Theorem 3.1.**
Let $\bigl(\Gamma_t(\mathbb{G}_{1/2}), t \geq 0\bigr)$ denote the subordinator characterized by:

$$E(e^{-\lambda\Gamma_t(\mathbb{G}_{1/2})}) = \exp\Bigl\{-t \int_0^\infty (1 - e^{-\lambda x}) \frac{dx}{x} E(e^{-x\mathbb{G}_{1/2}})\Bigr\} \qquad (189)$$

*The following explicit formulae hold:*
**1. Laplace transform of $\Gamma_t(\mathbb{G}_{1/2})$.**

$$E(e^{-\lambda\Gamma_t(\mathbb{G}_{1/2})}) = \exp\Bigl\{-t \int_0^\infty (1 - e^{-\lambda x}) \frac{dx}{x} I_0\Bigl(\frac{x}{2}\Bigr) e^{-\frac{x}{2}}\Bigr\} \qquad (190)$$

$$= \bigl(\sqrt{1+\lambda} - \sqrt{\lambda}\bigr)^{2t} = \Bigl(\frac{1}{\sqrt{1+\lambda} + \sqrt{\lambda}}\Bigr)^{2t} = \bigl(1 + 2\lambda + 2\sqrt{\lambda(1+\lambda)}\bigr)^{-t} \qquad (191)$$

where, in (190), $I_0$ denotes the modified Bessel function with index 0. (see [34], p. 108)

**2. Laws of $\Gamma_t(\mathbb{G}_{1/2})$ and of $D_t(\mathbb{G}_{1/2})$.**

$$\Gamma_t(\mathbb{G}_{1/2}) \stackrel{(law)}{=} \frac{\gamma_t}{\beta_{1/2,1/2+t}}, \quad D_t(\mathbb{G}_{1/2}) \stackrel{(law)}{=} \frac{1}{\beta_{1/2,1/2+t}} \stackrel{(law)}{=} 1 + \frac{\gamma_{t+1/2}}{\gamma_{1/2}} \quad (192)$$

The density of $f_{\Gamma_t(\mathbb{G}_{1/2})}$ equals:

$$f_{\Gamma_t(\mathbb{G}_{1/2})}(x) = \frac{2^{2t}\Gamma(1+t)}{2\pi\,\Gamma(2t)}\, x^{t-1}\Big(\int_0^1 e^{-xy}\big(y(1-y)\big)^{t-1/2}dy\Big)1_{[0,\infty]}(x) \quad (193)$$

**3. Wiener-Gamma representation of $\Gamma_t(\mathbb{G}_{1/2})$.**

For any $t > 0$:

$$\Gamma_t(\mathbb{G}_{1/2}) \stackrel{(law)}{=} \int_0^t \frac{d\gamma_u}{\sin^2\big(\frac{\pi u}{2t}\big)} \stackrel{(law)}{=} \int_0^t \frac{d\gamma_u}{\cos^2\big(\frac{\pi u}{2t}\big)} \quad (194)$$

Here is the dual version of Theorem 3.1:

**Theorem 3.1\*.** (Cifarelli and Melilli, [9]) *Let $\big(\Gamma_t(1/\mathbb{G}_{1/2}),\, t \geq 0\big)$ denote the subordinator characterized by:*

$$E(e^{-\lambda \Gamma_t(1/\mathbb{G}_{1/2})}) = \exp\Big\{-t\int_0^\infty (1 - e^{-\lambda x})\frac{dx}{x}\, E(e^{-\frac{x}{\mathbb{G}_{1/2}}})\Big\} \quad (195)$$

*Then:*

**1. Laplace transform of $\Gamma_t(1/\mathbb{G}_{1/2})$.**

$$E(e^{-\lambda \Gamma_t(1/\mathbb{G}_{1/2})}) = \Big(\frac{2}{1+\sqrt{1+\lambda}}\Big)^{2t} \qquad (\lambda, t \geq 0) \quad (196)$$

**2. Laws of $\Gamma_t(1/\mathbb{G}_{1/2})$ and of $D_t(1/\mathbb{G}_{1/2})$.**

$$\Gamma_t(1/\mathbb{G}_{1/2}) \stackrel{(law)}{=} \gamma_t \cdot \beta_{t+1/2,\, t+1/2}, \quad D_t(1/\mathbb{G}_{1/2}) \stackrel{(law)}{=} \beta_{t+1/2,\, t+1/2} \quad (197)$$

The density $f_{\Gamma_t(1/\mathbb{G}_{1/2})}$ of $\Gamma_t\big(\frac{1}{\mathbb{G}_{1/2}}\big)$ equals:

$$f_{\Gamma_t(1/\mathbb{G}_{1/2})}(x) = \frac{t \cdot 2^{2t}}{\Gamma(1/2+t)\sqrt{\pi}}\, x^{t-1}\Big(\int_0^1 e^{-\frac{x}{y}}\frac{(1-y)^{t-1/2}}{\sqrt{y}}\,dy\Big)1_{[0,\infty]}(x) \quad (198)$$

**3. Wiener-Gamma representation of $\Gamma_t(1/\mathbb{G}_{1/2})$.**

For any $t \geq 0$:

$$\Gamma_t\Big(\frac{1}{\mathbb{G}_{1/2}}\Big) \stackrel{(law)}{=} \int_0^t \cos^2\Big(\frac{\pi u}{2t}\Big)d\gamma_u \stackrel{(law)}{=} \int_0^t \sin^2\Big(\frac{\pi u}{2t}\Big)d\gamma_u \quad (199)$$

**Remark 3.2.**

**1)** Theorem 3.1*. has been obtained by Cifarelli and Melilli [9]. It would be possible to prove Theorem 3.1 by first using Theorem 3.1* and then by applying the duality Theorem 2.1. In fact, we shall operate conversely, as we shall first prove Theorem 3.1, then we shall show that Theorem 3.1* may be deduced from it, due to the duality Theorem 2.1.

**2)** By comparing formula (190) with formula (76), we deduce:

$$(\mathcal{J}_t^{(0)},\ t\geq 0) \stackrel{(law)}{=} \left(\frac{1}{2}\Gamma_t(\mathbb{G}_{1/2}),\ t\geq 0\right)$$

**3)** We come back to (196) and we write:

$$\frac{2}{1+\sqrt{1+\lambda}} = \frac{1}{1+\frac{1}{2}(\sqrt{1+\lambda}-1)} = \int_0^\infty e^{-u-\frac{u}{2}(\sqrt{1+\lambda}-1)}du$$
$$= E\left(\exp -\lambda S_{1/2}^{(1)}\left(\frac{\mathfrak{e}}{2}\right)\right) \qquad (200)$$

where $(S_{1/2}^{(1)}(t),\ t\geq 0)$ denotes the stable (1/2) subordinator, Esscher transformed, with the Esscher transformation with parameter 1 (see [45]). In (200), $\mathfrak{e}$ denotes a standard exponential r.v. independent from $(S_{1/2}^{(1)}(t),\ t\geq 0)$. Hence:

$$\left(\frac{2}{1+\sqrt{1+\lambda}}\right)^{2t} = E\left(\exp\left(-\lambda S_{1/2}^{(1)}\left(\frac{\gamma_{2t}}{2}\right)\right)\right) \qquad (201)$$

where the subordinators $(S_{1/2}^{(1)}(t),\ t\geq 0)$ and $(\gamma_{2t},\ t\geq 0)$ featured in (201) are being assumed independent. This formula led James and Yor to consider, more generally, the subordinator $(S_{1/2}^{(\nu)}(\gamma_t),\ t\geq 0)$. In the article [31], the following is obtained:

$$\left(\Gamma_t(\mathbb{G}_{1/2}),\ t\geq 0\right) \stackrel{(law)}{=} \left(S_{1/2}(\gamma_{2t}) + \widetilde{S}_{1/2}^{(1)}(\gamma_{2t}),\ t\geq 0\right) \qquad (202)$$

where, on the right-hand side of (202), the three subordinators $S_{1/2}, \widetilde{S}_{1/2}^{(1)}$ and $\gamma$ are assumed independent, and $(S_{1/2}(t),\ t\geq 0) \stackrel{(law)}{=} (\widetilde{S}_{1/2}(t),\ t\geq 0)$.

**3.2.b. <u>Proof of Theorem 3.1.</u>**

*i)* We already prove (191). In fact, (191) follows immediately from (176), since:

$$\begin{aligned}
E[e^{-\lambda\Gamma_t(\mathbb{G}_{1/2})}] &= \left(E[e^{-\lambda\Gamma_{1/2}(\mathbb{G}_{1/2})}]\right)^{2t} \\
&= \left(\sqrt{1+\lambda}-\sqrt{\lambda}\right)^{2t} = \left(\frac{1}{\sqrt{1+\lambda}+\sqrt{\lambda}}\right)^{2t} \\
&= \left(\frac{1}{1+2\lambda+2\sqrt{\lambda(1+\lambda)}}\right)^t
\end{aligned}$$

*ii)* We now prove (190), which is equivalent to:

$$E(e^{-x\mathbb{G}_{1/2}}) = e^{-\frac{x}{2}} I_0\left(\frac{x}{2}\right), \qquad x \geq 0 \tag{203}$$

From the Lipschitz-Hankel formula (see [52] or [38], Th. 1.1):

$$\nu \int_0^\infty e^{-ax} I_\nu(x) \frac{dx}{x} = \left(a + \sqrt{a^2 - 1}\right)^{-\nu} \qquad (a \geq 1, \nu > 0) \tag{204}$$

we deduce:

$$\nu \int_0^\infty (1 - e^{-\lambda x}) e^{-\frac{x}{2}} I_\nu\left(\frac{x}{2}\right) \frac{dx}{x} = 1 - \left(2\lambda + 1 + 2\sqrt{\lambda(1+\lambda)}\right)^{-\nu} \tag{205}$$

then, letting $\nu \to 0$ in (205):

$$\int_0^\infty (1 - e^{-\lambda x}) e^{-\frac{x}{2}} I_0\left(\frac{x}{2}\right) \frac{dx}{x} = \log\left(2\lambda + 1 + 2\sqrt{\lambda(1+\lambda)}\right) \tag{206}$$

Hence:

$$\begin{aligned}
E(e^{-\lambda \Gamma_1(\mathbb{G}_{1/2})}) &= \exp\left\{-\int_0^\infty (1 - e^{-\lambda x}) \frac{dx}{x} E(e^{-x\mathbb{G}_{1/2}})\right\} \\
&= \left(\frac{1}{\sqrt{1+\lambda} + \sqrt{\lambda}}\right)^2 \\
&= \frac{1}{2\lambda + 1 + 2\sqrt{\lambda(1+\lambda)}} \\
&= \exp\left\{-\int_0^\infty (1 - e^{-\lambda x}) \frac{dx}{x} e^{-\frac{x}{2}} I_0\left(\frac{x}{2}\right)\right\}
\end{aligned}$$

This formula, which is of interest by itself, shall not be used in the sequel of this proof.

*iii)* We now prove point 2 of Theorem 3.1.
For this purpose, we shall use the following property of hypergeometric functions: (see [34], p. 238). Let, for $\alpha, \delta, \gamma$ reals (cf [34], p. 239):

$$F(\alpha, \delta, \gamma\,;\,z) = \frac{\Gamma(\gamma)}{\Gamma(\delta)\Gamma(\gamma - \delta)} \int_0^1 t^{\delta - 1} (1-t)^{\gamma - \delta - 1} (1 - tz)^{-\alpha} dt \tag{207}$$

$(\gamma > \delta > 0,\ |z| < 1)$

$$= E\left((1 - z\beta_{\delta, \gamma - \delta})^{-\alpha}\right) \tag{208}$$

**Lemma 3.3.** *Let $a, b, c$ three positive reals. Then:*

1) $\displaystyle E\left(e^{-\lambda \frac{\gamma_a}{\beta_{b,c}}}\right) = \frac{\Gamma(a+b)\,\Gamma(b+c)}{\Gamma(b)\,(\Gamma(a+b+c)} \frac{1}{\lambda^a} F\left(a, a+b, a+b+c\,;\, -\frac{1}{\lambda}\right) \quad (\lambda \geq 0)$ (209)

2) $\quad E(e^{-\lambda \gamma_a \beta_{b,c}}) = F(a, b, b+c\,;\,-\lambda)$ $\hspace{4cm} (\lambda \geq 0)$ (210)

**Proof of Lemma 3.3.**

Since $E(e^{-\lambda \gamma_a}) = \frac{1}{(1+\lambda)^a}$, we have:

$$
\begin{aligned}
E(e^{-\lambda \frac{\gamma_a}{\beta_{b,c}}}) &= E\left[\left(\frac{1}{1+\frac{\lambda}{\beta_{b,c}}}\right)^a\right] \\
&= \frac{1}{B(b,c)} \int_0^1 \frac{x^a x^{b-1}}{(x+\lambda)^a} (1-x)^{c-1} dx
\end{aligned}
$$

(by definition of the beta $(b, c)$ law)

$$
\begin{aligned}
&= \frac{B(a+b, c)}{B(b, c)} \frac{1}{B(a+b, c)} \frac{1}{\lambda^a} \int_0^1 \frac{x^{a+b-1}(1-x)^{c-1}}{\left(1+\frac{x}{\lambda}\right)^a} dx \\
&= \frac{\Gamma(a+b)\,\Gamma(b+c)}{\Gamma(b)\,\Gamma(a+b+c)} \frac{1}{\lambda^a} E\left(\left(1+\frac{\beta_{a+b,c}}{\lambda}\right)^{-a}\right),
\end{aligned}
$$

hence point 1 of Lemma 3.3, from (208). Point 2 of Lemma 3.3 may be obtained from similar arguments.

We now end the proof of point 2 of Theorem 3.1. From Lemma 3.3, we have:

$$
\begin{aligned}
E(e^{-\lambda \frac{\gamma_t}{\beta_{\frac{1}{2},\frac{1}{2}+t}}}) &= \frac{\Gamma\left(t+\frac{1}{2}\right)\Gamma(t+1)}{\Gamma(\frac{1}{2})\Gamma(2t+1)} \frac{1}{\lambda^t} F\left(t, t+\frac{1}{2}, 2t+1\,;\,-\frac{1}{\lambda}\right) \\
&= \frac{\Gamma\left(t+\frac{1}{2}\right)\Gamma(t+1)}{\Gamma(\frac{1}{2})\Gamma(2t+1)} \frac{1}{\lambda^t} \left(\frac{1+\sqrt{1+\frac{1}{\lambda}}}{2}\right)^{-2t}
\end{aligned}
$$
(211)

since (see [34], p. 259):

$$
F\left(\alpha, \alpha+\frac{1}{2}, 2\alpha+1\,;\,z\right) = \left(\frac{1+\sqrt{1-z}}{2}\right)^{-2\alpha}
$$

Using Legendre's duplication formula ([34], p. 4):

$$
\Gamma\left(t+\frac{1}{2}\right)\Gamma(t+1) = 2^{-2t}\sqrt{\pi}(2t+1) \quad \left(\text{with } \Gamma\left(\frac{1}{2}\right) = \sqrt{\pi}\right)
$$

we obtain:

$$
E(e^{-\lambda \frac{\gamma_t}{\beta_{\frac{1}{2},\frac{1}{2}+t}}}) = \left(\frac{1}{\sqrt{1+\lambda}+\sqrt{\lambda}}\right)^{2t} = E(e^{-\lambda \Gamma_t(\mathbb{G}_{1/2})})
$$

from (191). Hence:

$$
\frac{\gamma_t}{\beta_{\frac{1}{2},\frac{1}{2}+t}} \stackrel{(law)}{=} \Gamma_t(\mathbb{G}_{1/2})
$$
(212)

The form of the density of $\Gamma_t(\mathbb{G}_{1/2})$ follows easily from (212).

*iv)* From Proposition 1.3 and (18), showing point 3 amounts to compute the inverse of the cumulative distribution function of $\mathbb{G}_{1/2}$. Now, we have:

$$F_{\mathbb{G}_{1/2}}(x) := P(\mathbb{G}_{1/2} \leq x) = \frac{1}{\pi} \int_0^x \frac{du}{\sqrt{u(1-u)}}$$

$$= \frac{2}{\pi} \text{ Arc sin } (\sqrt{x}) \qquad x \geq 0$$

so that $F_{\mathbb{G}_{1/2}}^{-1}(y) = \sin^2\left(\frac{\pi y}{2}\right) \qquad (0 \leq y \leq 1).$

### 3.2.c Proof of Theorem 3.1*.

*i)* We prove point 1, in two different ways :

• A direct proof.
We deduce, as a particular case of the beta-gamma algebra, that:

$$\frac{N^2}{2} \stackrel{(law)}{=} \gamma_{1/2} \stackrel{(law)}{=} \mathfrak{e}\beta_{\frac{1}{2},\frac{1}{2}} \tag{213}$$

where $N$ is a centered Gaussian variable, with variance 1, and $\mathfrak{e}$ a standard exponential. Hence:

$$P\left(\frac{N^2}{2} > x\right) = P(\mathfrak{e}\beta_{\frac{1}{2},\frac{1}{2}} > x) = P\left(\mathfrak{e} > \frac{x}{\beta_{\frac{1}{2},\frac{1}{2}}}\right) = E(e^{-\frac{x}{\beta_{\frac{1}{2},\frac{1}{2}}}})$$

Thus:

$$\int_0^\infty (1 - e^{-\lambda x}) \frac{dx}{x} E(e^{-\frac{x}{\mathbb{G}_{1/2}}}) = E\left(\int_0^{\gamma_{1/2}} \frac{dx}{x} (1 - e^{-\lambda x}) dx\right)$$

$$= \int_0^1 \frac{du}{u} \left(1 - E(\exp(-\lambda u \gamma_{1/2}))\right)$$

(from the change of variable: $x = u\gamma_{1/2}$)

$$= \int_0^1 \frac{du}{u}\left(1 - \frac{1}{(1+\lambda u)^{1/2}}\right)$$

$$= 2\left(\log(\sqrt{1+\lambda}+1)\right) - \log 2$$

hence:

$$E(e^{-\lambda \Gamma_t(\frac{1}{\mathbb{G}_{1/2}})}) = \exp\left\{-2t \log\left(\frac{\sqrt{1+\lambda}+1}{2}\right)\right\}$$

$$= \left(\frac{2}{1+\sqrt{1+\lambda}}\right)^{2t} = \left(\frac{4}{2+\lambda+2\sqrt{1+\lambda}}\right)^t$$

• We may also prove (196) with the help of the duality Theorem 2.1. Indeed,

from (101):

$$
\begin{aligned}
E(e^{-\lambda \Gamma_t(\frac{1}{\mathbb{G}_{1/2}})}) &= E\Big(\exp -\frac{1}{\lambda}\Gamma_t(\mathbb{G}_{1/2})\Big) \cdot \frac{\exp -t\,E(\log \mathbb{G}_{1/2})}{\lambda^t} \\
&= \frac{4^t}{\big(1+\frac{2}{\lambda}+2\sqrt{\frac{1}{\lambda}(1+\frac{1}{\lambda})}\big)^t} \cdot \frac{1}{\lambda^t}
\end{aligned}
$$

from (191) and since $E(\log \mathbb{G}_{1/2}) = -\log 4$, from (104) and (192). Thus:

$$
E(e^{-\lambda \Gamma_t(\frac{1}{\mathbb{G}_{1/2}})}) = \frac{4^t}{\big(2+\lambda+2\sqrt{1+\lambda}\big)^t} = \Big(\frac{2}{1+\sqrt{1+\lambda}}\Big)^{2t}
$$

*ii)* We now prove point 2 of Theorem 3.1*.
Of course, it would be possible to use point 2 of Lemma 3.3 to make this proof, in the same manner that we have used point 1 of this Lemma 3.3 to show point 2 of Theorem 3.1. In fact, we prefer to use formula (104) of the duality Theorem 2.1. We have:

$$
\begin{aligned}
f_{\Gamma_t(1/\mathbb{G}_{1/2})}(u) &= e^{-t\,E(\log \mathbb{G}_{1/2})} \frac{u^{t-1}}{\Gamma(t)} E(e^{-u\,D_t(\mathbb{G}_{1/2})}) \\
&= \frac{4^t}{\Gamma(t)} u^{t-1} E(e^{-\frac{u}{\beta_{\frac{1}{2},\frac{1}{2}+t}}}) \quad \text{(from (192))} \qquad (214)
\end{aligned}
$$

Considering now the Laplace transform of the two sides of (214) we obtain:

$$
\begin{aligned}
E(e^{-\lambda \Gamma_t(\frac{1}{\mathbb{G}_{1/2}})}) &= \frac{4^t}{\Gamma(t)} E\Big(\int_0^\infty e^{-\lambda u} u^{t-1} e^{-\frac{u}{\beta_{\frac{1}{2},\frac{1}{2}+t}}} du\Big) \\
&= \frac{4^t}{\Gamma(t)} E\Big[\Big(\frac{\beta_{\frac{1}{2},\frac{1}{2}+t}}{1+\lambda \beta_{\frac{1}{2},\frac{1}{2}+t}}\Big)^t\Big] \\
&= \frac{4^t B(t+\frac{1}{2},t+\frac{1}{2})}{B(\frac{1}{2},\frac{1}{2}+t)} E\Big[\Big(\frac{1}{1+\lambda \beta_{t+\frac{1}{2},\frac{1}{2}+t}}\Big)^t\Big] \qquad (215)
\end{aligned}
$$

Note that, by taking $\lambda = 0$ in (215), we get:

$$
\frac{4^t\,B(t+\frac{1}{2},t+\frac{1}{2})}{B(\frac{1}{2},t+\frac{1}{2})} = 1
$$

(which may also be recovered from the duplication formula for the Gamma function). Hence:

$$
E(e^{-\lambda \Gamma_t(\frac{1}{\mathbb{G}_{1/2}})}) = E\Big[\Big(\frac{1}{1+\lambda \beta_{t+1/2,\,t+1/2}}\Big)^t\Big] = E(e^{-\lambda \gamma_t \beta_{t+1/2,\,t+1/2}})
$$

The other formulae of Theorem 3.1* are now easily obtained. In particular, since: $F^{-1}_{1/\mathbb{G}_{1/2}}(u) = \frac{1}{F^{-1}_{\mathbb{G}_{1/2}}(1/u)} = \frac{1}{\sin^2 \frac{\pi u}{2}}$, we have, from Proposition 1.3:

$$
\Gamma_t(1/\mathbb{G}_{1/2}) \stackrel{(law)}{=} \int_0^t \frac{d\gamma_u}{F^{-1}_{1/\mathbb{G}_{1/2}}(\frac{u}{t})} = \int_0^t \sin^2\Big(\frac{\pi u}{2t}\Big) d\gamma_u.
$$

### 3.3. Study of the r.v.'s $\Gamma_{1-\alpha}(\mathbb{G}_\alpha)$ and $\Gamma_{1-\alpha}(1/\mathbb{G}_\alpha)$, $0 < \alpha < 1$

**3.3.a. Theorem 3.4.** *Let, for $0 < \alpha < 1$, $\bigl(\Gamma_t(\mathbb{G}_\alpha), \ t \geq 0\bigr)$ the subordinator characterized by:*

$$E(e^{-\lambda \Gamma_t(\mathbb{G}_\alpha)}) = \exp\Bigl\{ -t \int_0^\infty (1 - e^{-\lambda x}) \frac{dx}{x} \, E(e^{-x \mathbb{G}_\alpha}) \Bigr\} \qquad (216)$$

*We then have the following explicit formulae:*

**1. Laplace transform of $\Gamma_t(\mathbb{G}_\alpha)$.**

$$E(e^{-\lambda \Gamma_t(\mathbb{G}_\alpha)}) = \bigl((1+\lambda)^\alpha - \lambda^\alpha\bigr)^{\frac{t}{1-\alpha}} \qquad (\lambda, t \geq 0) \qquad (217)$$

**2. Distributions of $\Gamma_{1-\alpha}(\mathbb{G}_\alpha)$ and of $D_{1-\alpha}(\mathbb{G}_\alpha)$.**

$$\Gamma_{1-\alpha}(\mathbb{G}_\alpha) \stackrel{(law)}{=} \frac{\gamma_{1-\alpha}}{\beta_{\alpha,1}} \stackrel{(law)}{=} \frac{\gamma_{1-\alpha}}{U^{1/\alpha}} \qquad (218)$$

$$D_{1-\alpha}(\mathbb{G}_\alpha) \stackrel{(law)}{=} \frac{1}{\beta_{\alpha,1}} \stackrel{(law)}{=} \frac{1}{U^{1/\alpha}} \stackrel{(law)}{=} 1 + \frac{\gamma_1}{\gamma_\alpha} \qquad (219)$$

$\bigl($where $\frac{\gamma_1}{\gamma_\alpha}$ is, by definition, a Pareto r.v. with parameter $\alpha\bigr)$. *The density $f_{\Gamma_{1-\alpha}(\mathbb{G}_\alpha)}$ of $\Gamma_{1-\alpha}(\mathbb{G}_\alpha)$ equals:*

$$f_{\Gamma_{1-\alpha}(\mathbb{G}_\alpha)}(x) = \frac{\alpha}{\Gamma(1-\alpha)} \, \frac{1}{x^{1+\alpha}} \, (1 - e^{-x}) \, 1_{[0,\infty[}(x) \qquad (220)$$

**3. Wiener-Gamma representation of $\Gamma_t(\mathbb{G}_\alpha)$.**
*For every $t \geq 0$ and $0 < \alpha < 1$:*

$$\Gamma_t(\mathbb{G}_\alpha) \stackrel{(law)}{=} \int_0^t \left[1 + \left(\frac{\sin\bigl(\pi(1-\alpha)(\frac{t-u}{t})\bigr)}{\sin\bigl(\pi(1-\alpha)\frac{u}{t}\bigr)}\right)^{1/\alpha}\right] d\gamma_u \qquad (221)$$

$$\stackrel{(law)}{=} \int_0^t \left[1 + \left(\frac{\sin\bigl(\pi(1-\alpha)\frac{u}{t}\bigr)}{\sin\bigl(\pi(1-\alpha)(\frac{t-u}{t})\bigr)}\right)^{1/\alpha}\right] d\gamma_u$$

*In particular, for $t = 1 - \alpha$:*

$$\Gamma_{1-\alpha}(\mathbb{G}_\alpha) \stackrel{(law)}{=} \int_0^{1-\alpha} \left[1 + \left(\frac{\sin(\pi u)}{\sin \pi(1-\alpha-u)}\right)^{1/\alpha}\right] d\gamma_u \qquad (222)$$

We note that, for $\alpha = 1/2$, formula (222) coïncides with (194). The dual version of Theorem 3.4 is:

**Theorem 3.4\*.** *Let $0 < \alpha < 1$ and let $\bigl(\Gamma_t(1/\mathbb{G}_\alpha), \ t \geq 0\bigr)$ the subordinator characterized by:*

$$E(e^{-\lambda \Gamma_t(1/\mathbb{G}_\alpha)}) = \exp\Bigl\{ -t \int_0^\infty (1 - e^{-\lambda x}) \frac{dx}{x} \, E(e^{-\frac{x}{\mathbb{G}_\alpha}}) \Bigr\} \qquad (223)$$

*Then:*

1. **Laplace transform of $\Gamma_t(1/\mathbb{G}_\alpha)$.**

$$E(e^{-\lambda \Gamma_t(1/\mathbb{G}_\alpha)}) = \Big(\frac{1}{\alpha}\frac{(1+\lambda)^\alpha - 1}{\lambda}\Big)^{\frac{t}{1-\alpha}} \quad (224)$$

2. **Laws of $\Gamma_{1-\alpha}(1/\mathbb{G}_\alpha)$ and of $D_{1-\alpha}(1/\mathbb{G}_\alpha)$.**

$$\Gamma_{1-\alpha}(1/\mathbb{G}_\alpha) \stackrel{(law)}{=} \gamma_{1-\alpha} \cdot U, \quad D_{1-\alpha}(1/\mathbb{G}_\alpha) \stackrel{(law)}{=} U \quad (225)$$

(we note that the law of $D_{1-\alpha}(1/\mathbb{G}_\alpha)$ does not depend on $\alpha$, and it may be compared with (197): $\beta_{1,1} = U$). The density $f_{\Gamma_{1-\alpha}(1/\mathbb{G}_\alpha)}$ of $\Gamma_{1-\alpha}(1/\mathbb{G}_\alpha)$ equals:

$$f_{\Gamma_{1-\alpha}(1/\mathbb{G}_\alpha)}(x) = \frac{1}{\Gamma(1-\alpha)}\frac{1}{x^\alpha}\Big(\int_1^\infty e^{-xy}\frac{dy}{y^{\alpha+1}}\Big)1_{[0,\infty[}(x) \quad (226)$$

3. **Wiener-Gamma representation of $\Gamma_t(1/\mathbb{G}_\alpha)$.**

For every $t \geq 0$:

$$\Gamma_t\Big(\frac{1}{\mathbb{G}_\alpha}\Big) \stackrel{(law)}{=} \int_0^t \frac{d\gamma_u}{1 + \Big(\frac{\sin(\pi(1-\alpha)(\frac{t-u}{t}))}{\sin(\pi(1-\alpha)\frac{u}{t})}\Big)^{1/\alpha}} \quad (227)$$

$$\stackrel{(law)}{=} \int_0^t \frac{d\gamma_u}{1 + \Big(\frac{\sin(\pi(1-\alpha)\frac{u}{t})}{\sin(\pi(1-\alpha)(\frac{t-u}{t}))}\Big)^{1/\alpha}} \quad (228)$$

**Remark 3.5.**

1) Formula (221) was originally obtained by T. Fujita and M. Yor (see [23]).

2) We deduce from Proposition 1.6:

$$E\Big(D_{1-\alpha}\Big(\frac{1}{\mathbb{G}_\alpha}\Big)\Big) = E(\mathbb{G}_\alpha) = E(U) = \frac{1}{2} \quad \text{(from (225))} \quad (229)$$

Of course, we may verify directly that:

$$E(\mathbb{G}_\alpha) = \frac{1}{2}$$

starting from the formula:

$$E\big(\exp(-\lambda \mathfrak{e}\,\mathbb{G}_\alpha)\big) = \frac{\alpha}{1-\alpha}\frac{1-(1+\lambda)^{\alpha-1}}{(1+\lambda)^\alpha - 1} \quad (230)$$

(see [4], formula 1.19), then taking the derivative in $\lambda = 0$.

3) In [4], Section 4.2.1, the positive variables $X_{\alpha,1}$ whose Laplace transforms equal:

$$E(e^{-\lambda X_{\alpha,1}}) = \frac{1}{\alpha}\frac{(1+\lambda)^\alpha - 1}{\lambda} \quad (0 < \alpha < 1)$$

have been introduced. We note that, from (224):

$$X_{\alpha,1} \stackrel{(law)}{=} \Gamma_{1-\alpha}\left(\frac{1}{\mathbb{G}_\alpha}\right) \tag{231}$$

4) Since $\Gamma_{1-\alpha}(1/\mathbb{G}_\alpha) \stackrel{(law)}{=} \gamma_{1-\alpha} \cdot U$ and $\Gamma_{1-\alpha}(\mathbb{G}_\alpha) \stackrel{(law)}{=} \frac{\gamma_{1-\alpha}}{U^{1/\alpha}}$ the relations (221) and (227) may also be written as:

$$\frac{1}{U^{1/\alpha}} - 1 \stackrel{(law)}{=} \int_0^{1-\alpha} \left(\frac{\sin(\pi u)}{\sin(\pi(1-\alpha-u))}\right)^{1/\alpha} d_u(D_u^{(1-\alpha)})$$

$$\stackrel{(law)}{=} \int_0^{1-\alpha} \left(\frac{\sin(\pi(1-\alpha-u))}{\sin(\pi u)}\right)^{1/\alpha} d_u(D_u^{(1-\alpha)})$$

and

$$U \stackrel{(law)}{=} \int_0^{1-\alpha} \frac{d_u(D_u^{(1-\alpha)})}{1 + \left(\frac{\sin(\pi(1-\alpha-u))}{\sin(\pi u)}\right)^{1/\alpha}}$$

$$\stackrel{(law)}{=} \int_0^{1-\alpha} \frac{d_u(D_u^{(1-\alpha)})}{1 + \left(\frac{\sin(\pi u)}{\sin \pi(1-\alpha-u)}\right)^{1/\alpha}} \tag{232}$$

In particular, for $\alpha = \frac{1}{2}$, one finds:

$$U \stackrel{(law)}{=} \int_0^{1/2} \cos^2(\pi u) \, d_u(D_u^{(1/2)}) \stackrel{(law)}{=} \int_0^{1/2} \sin^2(\pi u) d_u(D_u^{(1/2)})$$

### 3.3.b. Proof of Theorem 3.4.

*i)* Point 1) is an immediate consequence of (176) and point 2) of (177), as may be shown after some elementary computations.

*ii)* We now prove point 3. From Proposition 1.3, it amounts to compute $F_{\mathbb{G}_\alpha}$ and its inverse. However, from (182) and (181), we find, successively:

$$F_{\mathbb{G}_\alpha}(x) = P(\mathbb{G}_\alpha < x) = P\left(\frac{Z_{1-\alpha}^{1/\alpha}}{1 + Z_{1-\alpha}^{1/\alpha}} < x\right)$$

$$= P\left(Z_{1-\alpha} < \left(\frac{x}{1-x}\right)^\alpha\right) \quad (0 \le x \le 1)$$

i.e.:

$$F_{\mathbb{G}_\alpha}(x) = F_{Z_{1-\alpha}}\left(\left(\frac{x}{1-x}\right)^\alpha\right) \quad (0 \le x \le 1) \tag{233}$$

hence:

$$F_{\mathbb{G}_\alpha}^{-1}(x) = \frac{\left(F_{Z_{1-\alpha}}^{-1}(x)\right)^{1/\alpha}}{1 + \left(F_{Z_{1-\alpha}}^{-1}(x)\right)^{1/\alpha}} \tag{234}$$

then:

$$\begin{aligned}
F_{Z_{1-\alpha}}(x) &= P(Z_{1-\alpha} \leq x) = \frac{\sin \pi(1-\alpha)}{\pi(1-\alpha)} \int_0^x \frac{dy}{y^2 + 2y(\cos \pi(1-\alpha)) + 1} \\
&= 1 - \frac{1}{\pi(1-\alpha)} \text{ arc tg} \left( \frac{\sin \pi(1-\alpha)}{\cos \pi(1-\alpha) + x} \right) \qquad (235) \\
&= \wedge_{1-\alpha}(x) \qquad \text{with the notation of (126)}
\end{aligned}$$

and, from (185):

$$F_{Z_{1-\alpha}}^{-1}(x) = \frac{\sin(\pi(1-\alpha)x)}{\sin(\pi(1-\alpha)(1-x))} \qquad (236)$$

hence, plugging (236) in (234), we obtain:

$$\frac{1}{F_{\mathbb{G}_\alpha}^{-1}(x)} = 1 + \left( \frac{\sin(\pi(1-\alpha)(1-x))}{\sin(\pi(1-\alpha)x)} \right)^{1/\alpha} \qquad (237)$$

which, thanks to Proposition 1.3, proves point 3) of Theorem 3.4.

### 3.3.c. Proof of Theorem 3.4*.

We use the duality theorem. From formula (104):

$$f_{\Gamma_{1-\alpha}(\mathbb{G}_\alpha)}(x) \underset{x\to 0}{\sim} \frac{x^{-\alpha}}{\Gamma(1-\alpha)} e^{(1-\alpha) E(\log \mathbb{G}_\alpha)} \qquad (238)$$

which we compare with (220):

$$f_{\Gamma_{1-\alpha}(\mathbb{G}_\alpha)}(x) = \frac{\alpha}{\Gamma(1-\alpha)} \frac{1}{x^{1+\alpha}} (1 - e^{-x}) 1_{[0,\infty[}(x)$$

and we deduce:

$$E(\log \mathbb{G}_\alpha) = \frac{\log \alpha}{1-\alpha} \qquad (239)$$

Then, using this time (101):

$$\begin{aligned}
E(e^{-\lambda \Gamma_t(\frac{1}{\mathbb{G}_\alpha})}) &= \exp\left\{ -t\left( E(\log \mathbb{G}_\alpha) + \log \lambda + \psi_{\mathbb{G}_\alpha}\left(\frac{1}{\lambda}\right) \right) \right\} \\
&= \left[ \exp\left( -\frac{t}{1-\alpha} \log \alpha \right) \right] \frac{1}{\lambda^t} \left[ \left(1 + \frac{1}{\lambda}\right)^\alpha - \left(\frac{1}{\lambda}\right)^\alpha \right]^{\frac{t}{1-\alpha}}
\end{aligned}$$

from (239) and (217):

$$= \left[ \frac{1}{\alpha} \frac{(1+\lambda)^\alpha - 1}{\lambda} \right]^{\frac{t}{1-\alpha}} \qquad (240)$$

which establishes (224). The density of $\Gamma_{1-\alpha}\left(\frac{1}{\mathbb{G}_\alpha}\right)$ may be computed from the density formula (104) and the knowledge of the law of $D_{1-\alpha}(\mathbb{G}_\alpha) \overset{(law)}{=} \frac{1}{U^{1/\alpha}}$, given by (219). All the formulae in point 2 of Theorem 3.4* follow easily from the explicit expression of $f_{\Gamma_{1-\alpha}(\frac{1}{\mathbb{G}_\alpha})}$, and point 3 of Theorem 3.4* follows from:

$$F^{-1}_{1/\mathbb{G}_\alpha}(y) \cdot F^{-1}_{\mathbb{G}_\alpha}(1-y) = 1,$$

and from Proposition 1.3.

### 3.4. Study of the subordinators $\big(\Gamma_t(\mu + \mathbb{G}_0),\ t \geq 0\big)$ and $\big(\Gamma_t(\frac{1}{\mu+\mathbb{G}_0}),\ t \geq 0\big)$

We recall the following formula, which was established in [4] (see formulae (181) and (179) in [4]):

$$\exp\Big\{ - \int_0^\infty (1 - e^{-\lambda x})\frac{dx}{x}\, E(e^{-x(\mu+\mathbb{G}_0)})\Big\} = \frac{\log\big(1 + \frac{1}{\lambda+\mu}\big)}{\log\big(1 + \frac{1}{\mu}\big)} \qquad \lambda, \mu \geq 0 \quad (241)$$

This formula was obtained in [4] by "letting $\alpha$ tend to 0 in the study of the r.v.'s $\mathbb{G}_\alpha$. Here, $\mathbb{G}_0$ satisfies (see (175)):

$$\mathbb{G}_0 \stackrel{(law)}{=} \frac{1}{1 + \exp \pi C}, \quad \text{where } C \text{ denotes a standard Cauchy r.v.} \qquad (242)$$

**3.4.a. Theorem 3.6.** *Let $\mu > 0$ be fixed and let $\big(\Gamma_t(\mu + \mathbb{G}_0),\ t \geq 0\big)$ denote the subordinator which is characterized by:*

$$E(e^{-\lambda \Gamma_t(\mu+\mathbb{G}_0)}) = \exp\Big\{ - t\int_0^\infty (1 - e^{-\lambda x})\frac{dx}{x}\, E(e^{-x(\mu+\mathbb{G}_0)})\Big\} \qquad (243)$$

*Then:*

**1. Laplace transform of $\Gamma_t(\mu + \mathbb{G}_0)$.**

$$E(e^{-\lambda \Gamma_t(\mu+\mathbb{G}_0)}) = \bigg(\frac{\log\big(1 + \frac{1}{\lambda+\mu}\big)}{\log\big(1 + \frac{1}{\mu}\big)}\bigg)^t \qquad (244)$$

**2. Laws of $\Gamma_1(\mu + \mathbb{G}_0)$ and $D_1(\mu + \mathbb{G}_0)$.**

The densities of $\Gamma_1(\mu + \mathbb{G}_0)$ and $D_1(\mu + \mathbb{G}_0)$ equal:

$$f_{\Gamma_1(\mu+\mathbb{G}_0)}(x) = \frac{1}{\log\big(1 + \frac{1}{\mu}\big)}\, e^{-\mu x}\, \frac{1 - e^{-x}}{x}\, 1_{[0,\infty[}(x) \qquad (245)$$

$$f_{D_1(\mu+\mathbb{G}_0)}(x) = \frac{1}{\log\big(1 + \frac{1}{\mu}\big)}\, \frac{1}{x}\, 1_{[\frac{1}{\mu+1}, \frac{1}{\mu}]}(x) \qquad (246)$$

**3. Wiener-Gamma representation of $\Gamma_t(\mu + \mathbb{G}_0)$.**
For every $t \geq 0$:

$$\Gamma_t(\mu + \mathbb{G}_0) \stackrel{(law)}{=} \int_0^t \frac{1 + \exp\big(\mathrm{cotg}\big(\frac{\pi u}{t}\big)\big)}{1 + \mu\big(1 + \exp\big(\mathrm{cotg}\big(\frac{\pi u}{t}\big)\big)\big)}\, d\gamma_u \qquad (247)$$

Here is the dual version of Theorem 3.6.

**Theorem 3.6\*.** *Let $\mu \geq 0$ and $\left(\Gamma_t\left(\frac{1}{\mu+\mathbb{G}_0}\right),\ t \geq 0\right)$ denote the subordinator characterized by:*

$$E(e^{-\lambda \Gamma_t(\frac{1}{\mu+\mathbb{G}_0})}) = \exp\left\{ - t \int_0^\infty (1 - e^{-\lambda x})\frac{dx}{x} E(e^{-\frac{x}{\mu+\mathbb{G}_0}}) \right\} \qquad (248)$$

*Then:*
 1. **Laplace transform of $\Gamma_t\left(\frac{1}{\mu+\mathbb{G}_0}\right)$.**

$$E\left[e^{-\lambda \Gamma_t(\frac{1}{\mu+\mathbb{G}_0})}\right] = \left[\frac{1}{\lambda} \log\left(\frac{1+\lambda(1+\mu)}{1+\lambda\mu}\right)\right]^t \qquad (249)$$

 2. **Laws of $\Gamma_1\left(\frac{1}{\mu+\mathbb{G}_0}\right)$ and $D_1\left(\frac{1}{\mu+\mathbb{G}_0}\right)$.**

$$\Gamma_1\left(\frac{1}{\mu+\mathbb{G}_0}\right) \stackrel{(law)}{=} \mathfrak{e}(U+\mu), \quad D_1\left(\frac{1}{\mu+\mathbb{G}_0}\right) \stackrel{(law)}{=} U+\mu \qquad (250)$$

*In particular, for $\mu=0$:*

$$\Gamma_1\left(\frac{1}{\mathbb{G}_0}\right) \stackrel{(law)}{=} \mathfrak{e} \cdot U, \quad D_1\left(\frac{1}{\mathbb{G}_0}\right) \stackrel{(law)}{=} U \qquad (251)$$

*The density of $\Gamma_1\left(\frac{1}{\mu+\mathbb{G}_0}\right)$ equals:*

$$f_{\Gamma_1\left(\frac{1}{\mu+\mathbb{G}_0}\right)}(x) = \left(\int_0^1 \frac{1}{\mu+y} \exp\left(-\frac{x}{\mu+y}\right) dy\right) 1_{[0,\infty[}(x) \qquad (252)$$

 3. **Wiener-Gamma representation of $\Gamma_t\left(\frac{1}{\mu+\mathbb{G}_0}\right)$.**
*For every $\mu \geq 0$ and $t \geq 0$:*

$$\Gamma_t\left(\frac{1}{\mu+\mathbb{G}_0}\right) \stackrel{(law)}{=} \int_0^t \frac{1+\mu\left(1+\exp\left(-\cotg\frac{\pi u}{t}\right)\right)}{1+\exp\left(-\cotg\frac{\pi u}{t}\right)}\, d\gamma_u \qquad (253)$$

**Remark 3.7.**
1) We may take $\mu=0$ in the statement of Theorem 3.6*, since $E\left(\log^+(\mathbb{G}_0)\right) < \infty$, but not in that of Theorem 3.6.

2) Formula (249), with $\mu=0$, was obtained by Bondesson ([5], Ex. 3.3.1, p. 42) and the relation $D_1\left(\frac{1}{\mathbb{G}_0}\right) \stackrel{(law)}{=} U$ has been obtained, independently, by Diaconis and Freedman ([12], Ex. 7.4, p. 74) and by Cifarelli and Melilli ([9], Ex. 2, p. 1393).

3) The formula $D_1\left(\frac{1}{\mathbb{G}_0}\right) \stackrel{(law)}{=} U$ may be obtained by letting $\alpha \longrightarrow 0$, in (225) and formula $\Gamma_1\left(\frac{1}{\mathbb{G}_0}\right) = \mathfrak{e}.U$ is formula (225) with $\alpha = 0$.

**3.4.b. Proof of Theorem 3.6.**
Point 1 follows immediately from (241), whereas formula (245) is a consequence of:

$$\begin{aligned}
E(e^{-\lambda \Gamma_1(\mu+\mathbb{G}_0)}) &= \frac{\log\left(1+\frac{1}{\lambda+\mu}\right)}{\log\left(1+\frac{1}{\mu}\right)} = \frac{1}{\log\left(1+\frac{1}{\mu}\right)} \int_{\lambda+\mu-1}^{\lambda+\mu} \frac{dv}{1+v} \\
&= \frac{1}{\log\left(1+\frac{1}{\mu}\right)} \int_{\lambda+\mu-1}^{\lambda+\mu} dv \int_0^\infty e^{-x(1+v)} dx \\
&= \frac{1}{\log\left(1+\frac{1}{\mu}\right)} \int_0^\infty e^{-x} dx \int_{\lambda+\mu-1}^{\lambda+\mu} \exp(-xv) dv \\
&= \frac{1}{\log\left(1+\frac{1}{\mu}\right)} \int_0^\infty e^{-(\lambda+\mu)x} \frac{1-e^{-x}}{x} dx
\end{aligned}$$

On the other hand, formula (246), which yields the density of $f_{D_1(\mu+\mathbb{G}_0)}$ follows from:

$$E[e^{-\lambda \Gamma_1(\mu+\mathbb{G}_0)}] = \frac{\log\left(1+\frac{1}{\lambda+\mu}\right)}{\log\left(1+\frac{1}{\mu}\right)} = E\left(\frac{1}{1+\lambda D_1(\mu+\mathbb{G}_0)}\right)$$

$$\frac{1}{\log\left(1+\frac{1}{\mu}\right)} \int_{\frac{1}{\mu+1}}^{\frac{1}{\mu}} \frac{dx}{x} \cdot \frac{1}{1+\lambda x} \quad \left(\text{after writing } \frac{1}{x(1+\lambda x)} = \frac{1}{x} - \frac{\lambda}{1+\lambda x}\right)$$

Then, the last point in Theorem 3.6 follows, with the help of (242), from:

$$\begin{aligned}
F^{-1}_{\mu+\mathbb{G}_0}(y) &= \mu + \frac{1}{1+\exp \text{tg}\left(\frac{\pi}{2}(1-2y)\right)} \\
&= \frac{1+\mu(\exp(\cotg(\pi y))+1)}{1+\exp \cotg(\pi y)}
\end{aligned} \tag{254}$$

**Proof of Theorem 3.6*.**
It may be proven by using the duality theorem. From (104) and (245):

$$f_{\Gamma_1(\mu+\mathbb{G}_0)}(0) = e^{E(\log(\mu+\mathbb{G}_0))} = \frac{1}{\log\left(1+\frac{1}{\mu}\right)}$$

i.e.:
$$E\left(\log(\mu+\mathbb{G}_0)\right) = -\log\left(\log\left(1+\frac{1}{\mu}\right)\right) \tag{255}$$

We then apply the duality theorem:

$$\begin{aligned}
E(e^{-\lambda \Gamma_t(\frac{1}{\mu+\mathbb{G}_0})}) &= E(e^{-\frac{1}{\lambda} \Gamma_t(\mu+\mathbb{G}_0)}) \cdot \frac{1}{\lambda^t} \cdot e^{-tE\left(\log(\mu+\mathbb{G}_0)\right)} \\
&= \left(\frac{\log\left(1+\frac{1}{\mu}\right)}{\lambda}\right)^t \cdot \left(\frac{\log\left(1+\frac{1}{\frac{1}{\lambda}+\mu}\right)}{\log\left(1+\frac{1}{\mu}\right)}\right)^t
\end{aligned}$$

(from (255) and (249))

$$= \left[\frac{1}{\lambda} \log\left(\frac{1+\lambda(1+\mu)}{1+\lambda\mu}\right)\right]^t$$

We now compute the density of $\Gamma_1\left(\frac{1}{\mu+\mathbb{G}_0}\right)$ by using (104):

$$\begin{aligned}
f_{\Gamma_1\left(\frac{1}{\mu+\mathbb{G}_0}\right)}(x) &= e^{-E(\log(\mu+\mathbb{G}_0))} E(e^{-x\,D_1(\mu+\mathbb{G}_0)}) \\
&= \log\left(1+\frac{1}{\mu}\right) \cdot \int_{\frac{1}{\mu+1}}^{\frac{1}{\mu}} \frac{1}{y\log\left(1+\frac{1}{\mu}\right)} e^{-xy} dy
\end{aligned}$$

(from (255) and (246))

$$= \int_0^1 \frac{1}{\mu+y} e^{-\frac{x}{\mu+y}} dy \qquad (256)$$

Formula (250) now follows easily from (256).

We deduce, from Proposition 1.6, that:

$$\begin{aligned}
E\left[\Gamma_1\left(\frac{1}{\mu+\mathbb{G}_0}\right)\right] &= E(\mu+\mathbb{G}_0) = E\left(D_1(\mu+\mathbb{G}_0)\right) \\
&= E(U+\mu) = \frac{1}{2}+\mu.
\end{aligned}$$

In particular, we have:

$$E(\mathbb{G}_0) = E(U) = \frac{1}{2}.$$

Finally, the last point of Theorem 3.6* follows from Proposition 1.3 and from (254).

### 3.5. Study of the subordinators $\left(\Gamma_t(\mathbb{G}_1),\ t \geq 0\right)$ and $\left(\Gamma_t(1/\mathbb{G}_1),\ t \geq 0\right)$

In this Section (see (175)), $\mathbb{G}_1$ denotes a uniform r.v. on $[0,1]$.

**Theorem 3.8.** Let $\left(\Gamma_t(\mathbb{G}_1),\ t \geq 0\right)$ denote the subordinator characterized by:

$$E(e^{-\lambda \Gamma_t(\mathbb{G}_1)}) = \exp\left\{-t\int_0^\infty (1-e^{-\lambda x})\frac{dx}{x} E(e^{-x\mathbb{G}_1})\right\} \qquad (257)$$

Then:

**1. Laplace transform of $\Gamma_t(\mathbb{G}_1)$.**

$$E(e^{-\lambda\Gamma_t(\mathbb{G}_1)}) = \left(\frac{\lambda^\lambda}{(1+\lambda)^{1+\lambda}}\right)^t \qquad \lambda, t \geq 0 \qquad (258)$$

**2. Laws of $\Gamma_1(\mathbb{G}_1)$ and $D_1(\mathbb{G}_1)$.**

The expressions of the densities of $\Gamma_1(\mathbb{G}_1)$ and $D_1(\mathbb{G}_1)$ are:

$$f_{\Gamma_1(\mathbb{G}_1)}(x) = \frac{1}{\pi}\left(\int_0^1 e^{-xy} \frac{\sin(\pi y)dy}{y^y(1-y)^{1-y}}\right) 1_{[0,\infty[}(x) \qquad (259)$$

$$f_{D_1(\mathbb{G}_1)}(x) = \frac{\sin\left(\frac{\pi}{x}\right)}{\pi(x-1)^{1-\frac{1}{x}}} 1_{[1,\infty[}(x) \tag{260}$$

**3. Wiener-Gamma representation of $\Gamma_t(\mathbb{G}_1)$.**
For every $t \geq 0$:
$$\Gamma_t(\mathbb{G}_1) \stackrel{(law)}{=} t \int_0^t \frac{d\gamma_u}{u} \stackrel{(law)}{=} t \int_0^t \frac{d\gamma_u}{t-u} \tag{261}$$

Here is now the dual version of Theorem 3.8, which is due to Diaconis and Kemperman [13].

**Theorem 3.8*.** Let $\left(\Gamma_t\left(\frac{1}{\mathbb{G}_1}\right), t \geq 0\right)$ denote the subordinator characterized by:
$$E(e^{-\lambda \Gamma_t(\frac{1}{\mathbb{G}_1})}) = \exp\left\{-t \int_0^\infty (1-e^{-\lambda x}) \frac{dx}{x} E(e^{-\frac{x}{\mathbb{G}_1}})\right\} \tag{262}$$

Then:

**1. Laplace transform of $\Gamma_t(1/\mathbb{G}_1)$.**
$$E(e^{-\lambda \Gamma_t(1/\mathbb{G}_1)}) = (e(1+\lambda)^{-\frac{1+\lambda}{\lambda}})^t \quad t, \lambda \geq 0 \tag{263}$$

**2. Laws of $\Gamma_1(1/\mathbb{G}_1)$ and $D_1(1/\mathbb{G}_1)$.**
The densities of $\Gamma_1(1/\mathbb{G}_1)$ and $D_1(1/\mathbb{G}_1)$ equal:
$$f_{\Gamma_1(1/\mathbb{G}_1)}(x) = \frac{e}{\pi} \left( \int_0^1 e^{-\frac{x}{y}} \frac{\sin(\pi y)dy}{y^{y+1}(1-y)^{1-y}} \right) 1_{[0,\infty[}(x) \tag{264}$$

$$f_{D_1(1/\mathbb{G}_1)}(x) = \frac{e\sin(\pi x)}{\pi} \frac{1}{x^x(1-x)^{1-x}} 1_{[0,1]}(x) \tag{265}$$

**3. Wiener-Gamma representation of $\Gamma_t(1/\mathbb{G}_1)$.**
For every $t \geq 0$:
$$\Gamma_t(1/\mathbb{G}_1) \stackrel{(law)}{=} \frac{1}{t} \int_0^t (t-u)\, d\gamma_u \stackrel{(law)}{=} \frac{1}{t} \int_0^t u\, d\gamma_u \tag{266}$$

**Remark 3.9.**
1) Theorem 3.8*, and in particular formula (265), are due to P. Diaconis and J. Kemperman [13]. Formula (265) may also be found in ([8], Ex. 4.4, p. 98) corrected with multiplication by $e$ in that formula.

2) *i)* A r.v. $Z$, which takes values on $\mathbb{R}$, is said to be a Luria-Delbrück r.v. if it satisfies:
$$E(e^{-\lambda Z}) = \lambda^\lambda \quad (\lambda \geq 0) \tag{267}$$

and M. Möhle [39] determined the density $f_Z$ of this r.v.:
$$\begin{aligned} f_Z(x) &= \frac{1}{\pi} \int_0^\infty e^{-\frac{xt}{2}} \cos(xt + t\log t)\, dt \quad (x \in \mathbb{R}) \\ &= \frac{1}{\pi} \int_0^\infty e^{-xt - t\log t}(\sin \pi t) dt \end{aligned} \tag{268}$$

*ii)* On the other hand, it is proven in $\big([44], \S \text{ III}, 1.3, \text{p. } 1251, \text{with } \alpha = k = 1\big)$ that there exists a Wald couple $(X, H)$, which is infinitely divisible and such that:

- $H$ is positive and $X$ and $H$ are infinitely divisible.
  - $E(e^{-\frac{\lambda^2}{2}H}) \cdot E(e^{\lambda X}) = 1 \qquad \lambda \geq 0$
  - $E(e^{-\frac{\lambda^2}{2}H}) = (1+\lambda)^{-(1+\lambda)}e^{\lambda}$,
  - $E(e^{\lambda X}) = (1+\lambda)^{1+\lambda}e^{-\lambda}$ (269)

Thus:
$$E(e^{\lambda(X+1)}) = (1+\lambda)^{1+\lambda} \tag{270}$$

hence, from (258), with $t = 1$:
$$\Gamma_1(\mathbb{G}_1) - (1+X) \stackrel{(law)}{=} Z \tag{271}$$

*iii)* We denote by $\widetilde{Z}$ a r.v. defined from $Z$ via the following Esscher transform:
$$E\big(\varphi(\widetilde{Z})\big) := E\big(\varphi(Z) \cdot e^{-Z}\big) \qquad (\varphi \text{ Borel and bounded})$$

$\big(\text{we note that } E(e^{-Z}) = 1, \text{ from } (267)\big)$. Thus, we have:
$$E(e^{-\lambda \widetilde{Z}}) = E(e^{-(1+\lambda)Z}) = (1+\lambda)^{1+\lambda}$$

so that:
$$\widetilde{Z} \stackrel{(law)}{=} -(1+X) \tag{272}$$

hence, from (271):
$$\Gamma_1(\mathbb{G}_1) + \widetilde{Z} \stackrel{(law)}{=} Z \tag{273}$$

We also note that Letac [35] characterized the law defined by (272).

**3.5.b. Proofs of Theorems 3.8 and 3.8\*.**

*i)* We prove (258), which it suffices to obtain for $t = 1$. We have:
$$\frac{\partial}{\partial \lambda} \log\left(\frac{\lambda^{\lambda}}{(1+\lambda)^{1+\lambda}}\right) = \log\left(\frac{\lambda}{1+\lambda}\right) = -\log\left(1+\frac{1}{\lambda}\right) \tag{274}$$

whereas:
$$\int_0^{\infty} (1-e^{-\lambda x})\frac{dx}{x} \, E(e^{-x\,\mathbb{G}_1}) = \int_0^{\infty} (1-e^{-\lambda x})(1-e^{-x})\frac{dx}{x}$$

Hence:
$$-\frac{\partial}{\partial \lambda} \int_0^{\infty} (1-e^{-\lambda x})\frac{dx}{x} \, E(e^{-x\,\mathbb{G}_1}) = -\int_0^{\infty} e^{-\lambda x}\frac{dx}{x}(1-e^{-x})$$
$$= -\log\left(1+\frac{1}{\lambda}\right) \tag{275}$$

from Frullani's integral. (258) now follows from the comparison of (274) and (275).

*ii)* We now prove (263):

$$
\begin{aligned}
\int_0^\infty (1-e^{-\lambda x})\frac{dx}{x} E(e^{-\frac{x}{\mathbb{G}_1}}) &= \int_0^\infty (1-e^{-\lambda x})\frac{dx}{x} \int_1^\infty e^{-ux}\frac{du}{u^2} \\
&= \int_1^\infty \frac{du}{u^2} \int_0^\infty \frac{1-e^{-\lambda x}}{x} e^{-ux} dx \quad \text{(Fubini)} \\
&= \int_1^\infty \frac{du}{u^2} \log\left(1+\frac{\lambda}{u}\right)
\end{aligned}
$$

(from Frullani's formula (see [34], p. 6, and (3) above )).

$$
= \int_0^1 \log(1+\lambda w)dw = \frac{1}{\lambda}\left((1+\lambda)\log(1+\lambda)-\lambda\right)
$$

Thus:

$$
\begin{aligned}
E(e^{-\lambda \Gamma_1(1/\mathbb{G}_1)}) &= \exp\left\{-\int_0^\infty (1-e^{-\lambda x})\frac{dx}{x} E(e^{-\frac{x}{\mathbb{G}_1}})\right\} \\
&= \exp\left\{-\frac{1}{\lambda}\left((1+\lambda)\log(1+\lambda)-\lambda\right)\right\} = \frac{e}{(1+\lambda)^{\frac{1+\lambda}{\lambda}}}.
\end{aligned}
$$

We may also show (263) with the help of the duality Theorem:

$$
\begin{aligned}
E(e^{-\lambda \Gamma_1(1/\mathbb{G}_1)}) &= \frac{e^{-E(\log \mathbb{G}_1)}}{\lambda} E\left(\exp -\frac{1}{\lambda}\Gamma_1(\mathbb{G}_1)\right) \\
&= \frac{e}{\lambda} \frac{\left(\frac{1}{\lambda}\right)^{\frac{1}{\lambda}}}{\left(1+\frac{1}{\lambda}\right)^{1+\frac{1}{\lambda}}} = e(1+\lambda)^{-\frac{1+\lambda}{\lambda}}
\end{aligned}
$$

from (258) and since: $E(\log \mathbb{G}_1) = -E(\gamma_1) = -1$.

*iii)* Point 3) of Theorem 3.8 and 3.8* follows easily from Proposition 1.3 and from the fact that $F_{\mathbb{G}_1}(u) = u \quad (0 \leq u \leq 1)$.

*iv)* We now prove (264). For this purpose, we use the result of Diaconis and Kemperman [13], see also [42],

$$
f_{D_1(1/\mathbb{G}_1)}(x) = \frac{e \sin \pi x}{\pi} \frac{1}{x^x(1-x)^{1-x}} 1_{[0,1]}(x) \tag{276}
$$

The relation:

$$
\Gamma_1(1/\mathbb{G}_1) = \gamma_1 \cdot D_1(1/\mathbb{G}_1)
$$

now implies easily (264).

*v)* We show (259) with the help of formula (104) in the duality Theorem:

$$\begin{aligned} f_{\Gamma_1(\mathbb{G}_1)}(x) &= e^{E(\log \mathbb{G}_1)} \, E(e^{-x\, D_1(\frac{1}{\mathbb{G}_1})}) \\ &= \frac{1}{e} \, E(e^{-x\, D_1(\frac{1}{\mathbb{G}_1})}) \\ &= \int_0^1 e^{-xy} \frac{\sin \pi y}{\pi} \frac{dy}{y^y (1-y)^{1-y}} \end{aligned} \quad (277)$$

from (276). Finally, in order to show (260), we apply (107):

$$\begin{aligned} f_{D_1(\mathbb{G}_1)}(x) &= x^{-1} \, e^{E(\log \mathbb{G}_1)} \, f_{D_1(1/\mathbb{G}_1)}\left(\frac{1}{x}\right) \\ &= \frac{1}{e\,x} \, \frac{e \sin\left(\frac{\pi}{x}\right)}{\pi} \, 1_{[0,1]}\left(\frac{1}{x}\right) \, \frac{1}{\left(\frac{1}{x}\right)^{\frac{1}{x}} \left(1-\frac{1}{x}\right)^{1-\frac{1}{x}}} \quad \text{(from (276))} \\ &= \frac{\sin\left(\frac{\pi}{x}\right)}{\pi (x-1)^{1-\frac{1}{x}}} \, 1_{[1,\infty[}(x) \end{aligned}$$

**Appendix**

**Interpolation between the subordinators $\bigl(\Gamma_t(1/G),\ t \geq 0\bigr)$ and $(\gamma_t,\ t \geq 0)$**

We denote, for every $u \geq 0$, by $\sigma_u : \mathbb{R} \longrightarrow \mathbb{R}$ the decreasing function defined by:

$$\sigma_u(x) := \frac{x \sinh u + \cosh u}{x \cosh u + \sinh u} = \frac{x \tanh u + 1}{x + \tanh u} \quad (278)$$

Since the image of $\mathbb{R}_+$ by $\sigma_u$ is equal to $]\tanh u, \coth u]$, then for every positive r.v. $G$, and every $u > 0$, we have: $E\bigl(|\log \sigma_u(G)|\bigr) < \infty$.

Let $\bigl(\Gamma_t(\sigma_u(G)),\ t \geq 0\bigr)$ denote the subordinator defined by:

$$\begin{aligned} E\bigl(\exp - \lambda \Gamma_t(\sigma_u(G))\bigr) &= \exp\bigl\{-t\psi_{\sigma_u(G)}(\lambda)\bigr\} \\ &= \exp\Bigl\{-t \int_0^\infty (1 - e^{-\lambda x}) \frac{dx}{x} E(e^{-x\sigma_u(G)})\Bigr\} \end{aligned} \quad (279)$$

Since $\sigma_0(G) = \frac{1}{G}$ and $\sigma_\infty(G) = 1$, we have:

$$\bigl(\Gamma_t(\sigma_0(G)),\ t \geq 0\bigr) = \bigl(\Gamma_t(1/G),\ t \geq 0\bigr)$$
$$\bigl(\Gamma_t(\sigma_\infty(G)),\ t \geq 0\bigr) = (\gamma_t,\ t \geq 0)$$

Thus, the family of subordinators $\bigl(\Gamma_t(\sigma_u(G)),\ t \geq 0\bigr)$ interpolates, as $u$ describes $\mathbb{R}_+$, between $\bigl(\Gamma_t(1/G),\ t \geq 0\bigr)$ and $(\gamma_t,\ t \geq 0)$.

The aim of this appendix is to show that one may compute "explicitly" the Laplace transform, the density, and the Wiener-Gamma representation of the r.v.'s $\Gamma_t(\sigma_u(G))$ in terms of those of $\bigl(\Gamma_t(G),\ t \geq 0\bigr)$, for every $u \geq 0$.

**Theorem A.1.** *We have, for every $u \geq 0$:*

1) $\quad -\psi_G(\sigma_u(\lambda)) + \psi_{\sigma_u(G)}(\lambda) = \log\left(\lambda \cosh u + \sinh u\right) + k_u \quad (280)$

$$\text{with} \quad k_u = E\left[\log\left(\frac{G}{G \sinh u + \cosh u}\right)\right] \quad (281)$$

2) $\quad F^{-1}_{\sigma_u(G)}(y) = \sigma_u\left(F_G^{-1}(1-y)\right) \quad (0 \leq y \leq 1) \quad (282)$

*hence:*

$$\Gamma_t\bigl(\sigma_u(G)\bigr) \stackrel{(law)}{=} \int_0^t \frac{d\gamma_s}{\sigma_u\bigl[F_G^{-1}(1-\frac{s}{t})\bigr]} \stackrel{(law)}{=} \int_0^t \frac{d\gamma_s}{\sigma_u\bigl[F_G^{-1}(\frac{s}{t})\bigr]} \quad (283)$$

3) *The density $f_{\Gamma_t(\sigma_u(G))}$ of $\Gamma_t\bigl(\sigma_u(G)\bigr)$, equals:*

$$f_{\Gamma_t(\sigma_u(G))}(x)$$
$$= \frac{e^{-x \tanh u - t k_u}}{\cosh u} E\left\{e^{-\Gamma_t(G)\tanh u}\left(\frac{\Gamma_t(G)}{x}\right)^{\frac{1-t}{2}} \cdot J_{t-1}\left(\frac{2\sqrt{y\Gamma_t(G)}}{\cosh u}\right)\right\} 1_{x>0}$$
$$= x^{t-1} \frac{e^{-x \tanh u - t k_u}}{\Gamma(t)} E\left[(\cosh u + D_t(G)\sinh u)^{-t} e^{-\frac{x D_t(G)}{\cosh u(\cosh u + D_t(G)\sinh u)}}\right]$$

(284)

**Remark A.2.**
**1.** For $u = 0$, we have: $\sigma_0(G) = \frac{1}{G}$, $k_0 = E(\log G)$ and the preceding formulae and indeed those of Theorem 2.1.

**2.** As $u \longrightarrow \infty$, the relation (283) becomes, by passage to the limit:

$$\Gamma_t\bigl(\sigma_\infty(G)\bigr) \stackrel{(law)}{=} \int_0^t d\gamma_s = \gamma_t \quad (285)$$

and it is not difficult to see, starting from (284), that:

$$f_{\Gamma_t(\sigma_u(G))}(x) \xrightarrow[u \to \infty]{} \frac{1}{\Gamma(t)} e^{-x} x^{t-1} 1_{[0,\infty[}(x) = f_{\gamma_t}(x)$$

**3.** More generally than (280), if we take:

$$\sigma(x) := \frac{ax+b}{cx+d} \quad \text{and} \quad \widetilde{\sigma}(\lambda) := \frac{dx+b}{cx+a}, \text{ we have :}$$
$$-\psi_G\left(\widetilde{\sigma}(\lambda)\right) + \psi_{\sigma(G)}(\lambda) = \log\left(c\lambda + a\right) + k(\sigma) \quad (286)$$
$$\text{with} \quad k(\sigma) := E\left(\frac{G}{aG+b}\right) \quad (287)$$

We observe that (286) and (287) follow immediately from:

$$-\psi_G(\widetilde{\sigma}(\lambda)) + \psi_{\sigma(G)}(\lambda) = -E\left[\log\left(1 + \frac{\widetilde{\sigma}(\lambda)}{G}\right)\right] + E\left[\log\left(1 + \frac{\lambda}{\sigma(G)}\right)\right]$$

$$= E\left[\log\left(1 + \frac{\lambda(cG+d)}{aG+b}\right)\right] - E\left[\log\left(1 + \frac{d\lambda+b}{c\lambda+a}\frac{1}{G}\right)\right]$$

$$= E\left[\log\left(\frac{(aG+b) + \lambda(cG+d)}{aG+b} \cdot \frac{(c\lambda+a)\,G}{(c\lambda+a)\,G + d\lambda+b}\right)\right]$$

$$= E\left[\log\frac{(c\lambda+a)G}{aG+d}\right].$$

**Proof of Theorem A.1.**
Point 1) is a particular case of point 3 of Remark A.2 and point 2) of Theorem A.1 is trivial. Let us prove point 3. From (280), after multiplying by $t$, and exponentiating, we deduce:

$$E\left[e^{-\lambda\Gamma_t(\sigma_u(G))}\right] = e^{-tk_u}\, E\left[e^{-\Gamma_t(G)\sigma_u(\lambda)}\, \frac{1}{(\lambda\cosh u + \sinh u)^t}\right] \quad (288)$$

then, multiplying each side of (288) by $\exp(-\beta\lambda)$ and integrating in $\lambda$ from $-\tanh u$ to $+\infty$, we obtain:

$$E\left(\int_{-\tanh u}^{\infty} e^{-\lambda\Gamma_t(\sigma_u(G)) - \beta\lambda}\, d\lambda\right)$$
$$= e^{-tk_u}\, E\left(\int_{-\tanh u}^{\infty} e^{-\Gamma_t(G)\sigma_u(\lambda) - \beta\lambda}\, \frac{d\lambda}{(\lambda\cosh u + \sinh u)^t}\right) \quad (289)$$

Making, on the LHS of (289) the change of variable: $\lambda\cosh u + \sinh u = \mu$, we obtain:

$$E\left[\frac{1}{(\beta + \Gamma_t(\sigma_u(G)))}\, e^{(\tanh u)(\Gamma_t(\sigma_u(G)) + \beta)}\right]$$
$$= \frac{e^{-tk_u}}{\cosh u}\, E\left[e^{-\Gamma_t(G)\tanh u + \beta\tanh u}\int_0^\infty \frac{e^{-\frac{1}{\cosh u}\left(\frac{\Gamma_t(G)}{\mu} + \beta\mu\right)}}{\mu^t}\, d\mu\right]$$

Thus, after simplifying by $\exp(\beta\tanh u)$ and using formula 5.10.25, p. 109 in [34], we obtain:

$$E\left[\frac{1}{(\beta + \Gamma_t(\sigma_u(G)))}\, e^{\tanh u\,\Gamma_t(\sigma_u(G))}\right]$$
$$= \frac{2e^{-tk_u}}{\cosh u}\, E\left[e^{-\Gamma_t(G)\tanh u}\left(\frac{\Gamma_t(G)}{\beta}\right)^{\frac{1-t}{2}} \cdot K_{1-t}\left(\frac{2\sqrt{\beta\,\Gamma_t(G)}}{\cosh u}\right)\right]$$

We then use the inversion formula of the Stieltjes transform, and thanks to a very similar computation to the one made in section 2.2.b, we arrive to the

identity:

$$e^{x\tanh u}f_{\Gamma_t(\sigma_u(G))}(x) = \frac{e^{-tk_u}}{\cosh u}E\left[e^{-\Gamma_t(G)\tanh u}\left(\frac{\Gamma_t(G)}{x}\right)^{\frac{1-t}{2}}J_{t-1}\left(\frac{2\sqrt{x\,\Gamma_t(G)}}{\cosh u}\right)\right]$$

then to (286), by using the series development of $J_{t-1}$:

$$f_{\Gamma_t(\sigma_u(G))}(x)$$
$$= \frac{e^{-x\tanh u - tk_u}x^{t-1}}{\Gamma(t)}\,E\left\{\left(\cosh u + D_t(G)\sinh u\right)^{-t}e^{-\frac{xD_t(G)}{\cosh u(\cosh u + D_t(G)\tanh u)}}\right\}$$

**Remark A.3.** Formulae (284) and (285), in the more general situation when: $\sigma(G) = \frac{aG+b}{cG+d}$ become:

$$f_{\Gamma_t(\sigma(G))}(x) = \frac{e^{-x\frac{a}{c}-tk(\sigma)}}{c}\,E\left(e^{-\frac{d}{c}\Gamma_t(G)}\left(\frac{\Gamma_t(G)}{x}\right)^{\frac{1-t}{2}}J_{t-1}\left(\frac{2\sqrt{x\,\Gamma_t(G)}}{c}\right)\right) \tag{290}$$

where $k(\sigma)$ is given by (287), and:

$$f_{\Gamma_t(\sigma(G))}(x) = \frac{x^{t-1}}{\Gamma(t)}\,e^{-x\frac{a}{c}-tk(\sigma)}\,E\left\{(c + dD_t(G))^{-t}\,e^{-\frac{x\,D_t(G)}{c(c+dD_t(G))}}\right\} \tag{291}$$

**List of notations**

$G, X, Z$ denote r.v.'s.

$\mu_G$: law of $G$; $f_G$: density of $G$; $F_G$: cumulative distribution function of $G$; $F_G^{-1}$: inverse function of $F_G$.

$\Gamma$: GGC r.v. with associated Thorin measure $\mu$.

$m$: total mass of $\mu$.

$\gamma_a$: Gamma r.v. with parameter $a$ $(a > 0)$ ; $\gamma_1 \stackrel{(law)}{=} \mathfrak{e}$, a standard exponential r.v.

$\beta_{a,b}$: a beta r.v. with parameters $a, b$ $(a, b > 0)$.

$U$: a uniform r.v. on $[0, 1]$.

$C$: a standard Cauchy r.v.

$Y_p$: a Bernoulli r.v. $(0 \leq p \leq 1)$. $\big(P(Y_p = 1) = p \,;\, P(Y_p = 0) = 1 - p\big)$

$S_\alpha$: a standard positive stable r.v., with index $\alpha$ $(0 < \alpha < 1)$, and density $f_\alpha$.

$\frac{\gamma_1}{\gamma_\delta}$: a Pareto r.v. of index $\delta$ ($\gamma_1$ and $\gamma_\delta$ independent).

$G_\theta$ generic notation for random variables considered in sections 2.5 and 2.6.

$(\mathbb{G}_\alpha, \, 0 \leq \alpha \leq 1)$ the family of r.v.'s introduced in 3.1. In particular:

$$\mathbb{G}_{1/2} \stackrel{(law)}{=} \beta_{\frac{1}{2}, \frac{1}{2}}, \quad \mathbb{G}_0 \stackrel{(law)}{=} \frac{1}{1 + \exp \pi C}, \quad \mathbb{G}_1 \stackrel{(law)}{=} U$$

$\Gamma_m(G)$ a GGC r.v. with associated Thorin measure $m \cdot \mu_G$.

$(\gamma_t, \, t \geq 0)$: standard gamma process (i.e.: subordinator with Lévy measure: $\frac{dx}{x} e^{-x}$).

$(D_t^{(m)}, \, 0 \leq t \leq m)$ Dirichlet process with parameter $m$.

$\big(\Gamma_t(G), \, t \geq 0\big)$: GGC subordinator with Thorin measure $\mu_G$, and Lévy density: $\frac{1}{x} E(e^{-xG})$.

$\psi_G$: the Bernstein function of the subordinator $\big(\Gamma_t(G), \, t \geq 0\big)$:

$$\psi_G(\lambda) = \int_0^\infty (1 - e^{-\lambda x}) \frac{dx}{x} E(e^{-xG})$$

$\widetilde{\Gamma}_t(h) = \int_0^t h(s) \, d\gamma_s$ the Wiener-Gamma integral of $h$.

$D_t(G)$ r.v. whose law is characterized by $\Gamma_t(G) \stackrel{(law)}{=} \gamma_t \, D_t(G)$

$D_t(G)$ satisfies:

$$D_t(G) \stackrel{(law)}{=} \int_0^t \frac{1}{F_G^{-1}(u/t)} \, d_u \, (D_u^{(t)})$$

$(\mathcal{C}_t, \, t \geq 0), (\mathcal{S}_t, \, t \geq 0), (\mathcal{T}_t, \, t \geq 0)$: hyperbolic subordinators.

$(\mathcal{J}_t^{(\nu)}, \, t \geq 0), (\mathcal{K}_t^{(\nu)}, \, t \geq 0)$: subordinators associated to the functions $I_\nu$ and $K_\nu$.

$K_\nu, I_\nu, J_\nu$: modified and unmodified Bessel functions with index $\nu$.

$\mathcal{S}^{(t)}(X)$: Stieltjes transform with index $t$ of the positive r.v. $X$:

$$\mathcal{S}^{(t)}(X)(\lambda) := E\left(\frac{1}{(\lambda + X)^t}\right) \qquad (\lambda \geq 0)$$

$(S_\alpha(t),\ t \geq 0)$ standard stable subordinator, with index $\alpha$ ($0 < \alpha < 1$).
$f_\alpha(x)$: density of $S_\alpha =^{(law)} S_\alpha(1)$.

# References

[1] BARNDORFF-NIELSEN, O. E., MAEJIMA, M. and SATO, K.I. (2006). Some classes of multivariate infinitely divisible distribution admitting stochastic integral representations *Bernoulli*, **12**, p. 1-33. MR2202318
[2] BERTOIN, J. (1996). *Lévy processes.* Cambridge Tracts in Mathematics, 121, Cambridge University Press. MR1406564
[3] BERTOIN, J. (2000). The convex minorant of the Cauchy process. *Electron. Comm. Probab.* **5** 51-55. MR1747095
[4] BERTOIN, J., FUJITA, T., ROYNETTE, B. and YOR, M. (2006). On a particular class of self-decomposable random variables: the duration of a Bessel excursion straddling an independent exponential time. *Prob. Math. Stat.* **26**, 315-366. MR2325310
[5] BONDESSON, L. (1992). *Generalized gamma convolutions and related classes of distributions and densities.* Lecture Notes in Statistics, **76**. Springer-Verlag, New York. MR1224674
[6] BONDESSON, L. (1981). Classes of infinitely divisible distributions and densities. *Z. Wahrsch. Verw. Gebiete* **57** 39-71. MR0623454
[7] BONDESSON, L. (1990). Generalized gamma convolutions and complete monotonicity. *Probab. Theory Related Fields* **85** 181–194. MR1050743
[8] CHAUMONT, L. and YOR, M. (2003). *Exercises in probability. A guided tour from measure theory to random processes, via conditioning.* Cambridge Series in Statistical and Probabilistic Mathematics, 13, Cambridge University Press. MR2016344
[9] CIFARELLI, D.M. and MELLILI, E. (2000). Some new results for Dirichlet priors. *Ann. Statist.* **28**, 1390–1413. MR1805789
[10] CIFARELLI, D. M. and REGAZZINI, E. (1979). Considerazioni generali sull'impostazione bayesiana di problemi non parametrici. Le medie associative nel contesto del processo aleatorio di Dirichlet I, II. *Riv. Mat. Sci. Econom. Social* **2**, 39–52. MR0573686
[11] CIFARELLI, D.M. and REGAZZINI, E. (1990). Distribution functions of means of a Dirichlet process. *Ann. Statist.* **18**, 429–442 (Correction in *Ann. Statist.* (1994) **22**, 1633–1634). MR1041402
[12] DIACONIS, P. and FREEDMAN, D. (1999). Iterated random functions. *SIAM Rev.*, **41**, 45-76. MR1669737
[13] DIACONIS, P. and KEMPERMAN, J. (1996). Some new tools for Dirichlet priors. In *Bayesian Statistics 5* (Bernardo, J.M., Berger, J.O., Dawid, A.P. and Smith, A.F.M., Eds.), 97–106. Oxford University Press, New York. MR1425401
[14] DOSS, H. and SELLKE, T. (1982). The tails of probabilities chosen from a Dirichlet prior. *Ann. Statist.* **10**, 1302–1305. MR0673666
[15] DUFRESNE, D. and YOR, M. (2007). In preparation.

[16] EPIFANI, I., GUGLIELMI, A. and MELILLI, E. (2004). Some new results on random Dirichlet variances. Technical Report IMATI 2004-15-MI (2004).

[17] EPIFANI, I., GUGLIELMI, A. and MELILLI, E. (2006). A stochastic equation for the law of the random Dirichlet variance. *Statist. Probab. Lett.* **76** 495–502. MR2266603

[18] EMERY, M. and YOR, M. (2004). A parallel between Brownian bridges and gamma bridges. Publ. Res. Inst. Math. Sci., Kyoto University, **40** 669-688. MR2074696

[19] ETHIER, S. N. and GRIFFITHS, R. C. (1993). The transition function of a Fleming-Viot process *Ann. Probab.* **21** 1571–1590. MR1235429

[20] FERGUSON, T.S. (1973). A Bayesian analysis of some nonparametric problems. *Ann. Statist.* **1** 209-230. MR0350949

[21] FEIGIN, P.D. and TWEEDIE, R.L. (1989). Linear functionals and Markov chains associated with Dirichlet processes. *Math. Proc. Philos. Soc.* **105**, p. 579-585. MR0985694

[22] FREEDMAN, D.A. On the asymptotic behavior of Bayes' estimates in the discrete case. (1963). *Ann. Math. Statist.*, **34**, 1386-1403. MR0158483

[23] FUJITA. T. and YOR, M. (2006). An interpretation of the results of the BFRY paper in terms of certain means of Dirichlet Processes.

[24] GOLDIE, C. A class of infinitely divisible random variables. (1967). *Proc. Cambridge Philos. Soc.* **63** 1141-1143. MR0215332

[25] GRIGELIONIS, B. (2007). Extended Thorin classes and stochastic integrals. *Lithuanian Math. Journal* **47**, n° 4, p. 406-411. MR2392715

[26] HANNUM, R.C., HOLLANDER, M. and LANGBERG, N.A. (1981). Distributional results for random functionals of a Dirichlet process *Ann. Probab.* **9** 665-670. MR0630318

[27] HJORT, N.L. and ONGARO, A. (2005). Exact inference for random Dirichlet means. *Stat. Inference Stoch. Process.* **8** 227-254. MR2177313

[28] JAMES, L.F. (2006). Gamma tilting calculus for GGC and Dirichlet means via applications to Linnik processes and occupation time laws for randomly skewed Bessel processes and bridges. http://arxiv.org/abs/math.PR/0610218. MR2275248

[29] JAMES, L.F. (2007). New Dirichlet mean identities. http://arxiv.org/abs/0708.0614.

[30] JAMES, L. F., LIJOI, A. and Prünster, I. (2008). Distributions of functionals of the two parameter Poisson-Dirichlet process. *Ann. Appl. Probab.* **18**, 521–551. MR2398765

[31] JAMES, L.F. and YOR, M. (2006). Tilted stable subordinators, Gamma time changes and Occupation Time of rays by Bessel Spiders. http://arxiv.org/abs/math.PR/0701049.

[32] JEANBLANC, M., PITMAN, J. and YOR, M. (2002). Self-similar processes with independent increments associated with Lévy and Bessel processes. *Stochastic Process. Appl.* **100** 223–231. MR1919614

[33] LAMPERTI, J. (1958). An occupation time theorem for a class of stochastic processes. Trans. Amer. Math. Soc., 88, p. 380-387 (1958). MR0094863

[34] LEBEDEV, N.N. (1972). *Special functions and their applications.* Revised

edition, translated from the Russian and edited by Richard A. Silverman.
[35] LETAC, G. (1985). A characterization of the Gamma distribution. *Adv. in Appl. Probab.* **17** 911-912. MR0809436
[36] LIJOI, A. and REGAZZINI, E. (2004). Means of a Dirichlet process and multiple hypergeometric functions. *Ann. Probab.* **32** 1469-1495. MR2060305
[37] LUKACS, E. (1970). *Characteristic functions.* Second edition, revised and enlarged. Hafner Publishing Co., New York. MR0346874
[38] MATSUMOTO, H., NGUYEN, L. and YOR, M. (2002). Subordinators related to the exponential functionals of Brownian bridges and explicit formulae for the semigroups of hyperbolic Brownian motions. Stochastic processes and related topics (Siegmundsburg, 2000), p. 213-235, Stochastic Monogr., 12, Taylor and Francis, London. MR1987318
[39] MÖHLE, M. (2005). Convergence results for compound Poisson distributions and applications to the standard Luria-Delbruck distribution. *J.Applied Prob.* **42** 620-631. MR2157509
[40] PATTERSON, S.J. (1998). *An introduction to the theory of the Riemann zeta-function.* Cambridge Studies in Advanced Mathematics, 14, Cambridge University Press. MR0933558
[41] PITMAN, J. and YOR, M. (2003). Infinitely divisible laws associated with hyperbolic functions. *Canad. J. Math.* **55** 292-330. MR1969794
[42] REGAZZINI, E., A. GUGLIELMI, A. and DI NUNNO, G. (2002). Theory and numerical analysis for exact distribution of functionals of a Dirichlet process. *Ann. Statist.* **30**, 1376-1411. MR1936323
[43] REVUZ, D. and YOR, M. (1999). *Continuous martingales and Brownian motion.* Third edition, Grundlehren der Mathematischen Wissenschaften (Fundamental Principles of Mathematical Sciences), 293, Springer-Verlag, Berlin. MR1725357
[44] ROYNETTE, B. and YOR, M. (2005). Couples de Wald indéfiniment divisibles. Exemples liés à la fonction gamma d'Euler et à la fonction zêta de Riemann. (Infinitely divisible Wald pairs: examples associated with the Euler gamma function and the Riemann zeta function.) *Ann. Inst. Fourier* (Grenoble) **55** 1219-1283. MR2157168
[45] SATO, K. I. (1999). *Lévy processes and infinitely divisible distributions.* Translated from the 1990 Japanese original. Revised by the author. Cambridge Studies in Advanced Mathematics, 68, Cambridge University Press. MR1739520
[46] SONG, R. and VONDRACEK, Z. (2006). Potential theory of special subordinators and subordinate killed stable processes. *J. Theor. Prob.* **19** 807-847. MR2279605
[47] STEUTEL, F. W. (1967). Note on the infinite divisibility of exponential mixtures. *Ann. Math. Statist.* **38** 1303-1305. MR0215339
[48] STEUTEL, F. W. and VAN HARN, K. (2004). *Infinite divisibility of probability distributions on the real line.* Monographs and Textbooks in Pure and Applied Mathematics, 259. Marcel Dekker. MR2011862
[49] THORIN, O. (1977). On the infinite divisibility of the lognormal distribution. *Scand. Actuar. J.* **3** 121-148. MR0552135

[50] VERSHIK, A., YOR, M. and TSILEVICH, N. (2004). On the Markov–Krein identity and quasi–invariance of the gamma process. *J. Math. Sci.* **121** 2303-2310. MR1879060
[51] VERVAAT, W. (1979). On a stochastic difference equation and a representation of non-negative infinitely divisible random variables. *Adv. in Appl. Probab.* **11**, 780-783. MR0544194
[52] WATSON, G. N. (1944). *A treatise on the theory of Bessel functions.* Cambridge University Press, Cambridge, England; The Macmillan Company, New York. MR0010746
[53] WIDDER, D.V. (1941). *The Laplace Transform.* Princeton Mathematical Series, v. 6. Princeton University Press, Princeton, N. J.. MR0005923
[54] WINKEL, M. (2005). Electronic foreign-exchange markets and passage events of independent subordinators. *J. Appl. Probab.* **42**, 138-152. MR2144899
[55] YOR, M. (2007). Some remarkable properties of the Gamma process. Festschrift for Dilip Madan, Advances in Mathematical Finance. p. 37-47, Applied Numer. Harmon And., eds M. Fu, R. Jarrow, R. Elliott, Birkhäuser Boston. MR2359361
[56] YAMATO, H. (1984). Characteristic functions of means of distributions chosen from a Dirichlet process. *Ann. Probab.* **12** 262-267. MR0723745